\documentclass[psamsfonts,11pt]{article}
\usepackage{graphicx}
\usepackage{amsmath}
\usepackage{amsfonts}
\usepackage{amssymb}

\newtheorem{theorem}{Theorem}

\newtheorem{corollary}[theorem]{Corollary}

\newtheorem{definition}[theorem]{Definition}
\newtheorem{example}[theorem]{Example}

\newtheorem{lemma}[theorem]{Lemma}

\newtheorem{proposition}[theorem]{Proposition}
\newtheorem{remark}[theorem]{Remark}

\newcommand{\R}{{\bf R}}
\newcommand{\Z}{{\bf Z}}
\newcommand{\Q}{{\bf Q}}
\newcommand{\C}{{\bf C}}
\newcommand{\Spec}{{\rm Spec}}
\newcommand{\Pf}{{\rm Pf}}

\begin{document}

\title{Pfaffians, the $G$-Signature Theorem and Galois de Rham discriminants}
\author{T. Chinburg\thanks{Supported by NSF Grant  DMS00-70433}, G.
Pappas\thanks{Supported by NSF Grant DMS02-01140.} and
M. J. Taylor\thanks{EPSRC Senior Research Fellow and a Royal Society Wolson
Merit award holder.}.}
\maketitle

\section{ INTRODUCTION.}

The principal objects of study in this article are the various bilinear forms
on the de Rham cohomology of an arithmetic variety which arise from duality
theory. S.~Bloch provided a fundamental new insight into de Rham discriminants
when he showed that, for a suitable arithmetic surface, the square root of the
de Rham discriminant is equal to the conductor of the surface (see [Bl]). This
striking result, which has been a fundamental influence on our work, may be
seen as extending to surfaces the important fact that the discriminant of a
ring of algebraic integers can be expressed in terms of Artin conductors. We
remark that Bloch's result has recently been extended to arithmetic varieties of
higher dimension by K.~Kato and T.~Saito in [KS].

We consider a regular scheme $\mathcal{X}$  which is projective and flat over
${\rm Spec} (\mathbf{Z} )$ of constant fibral dimension $d$, and which
supports an action by a finite group $G$.  We write $\mathcal{Y}$ for the
quotient scheme $\mathcal{X}/G$.  To each complex character $\theta$ of
$G$ we can associate the Artin-Hasse-Weil L-function $L (  \mathcal{Y}%
,\theta,s )$ which conjecturally satisfies a functional equation
\[
L (  \mathcal{Y},\theta,s )  =\varepsilon (  \mathcal{Y}%
,\theta )  A (  \mathcal{Y},\theta )  ^{-s}L (
\mathcal{Y},\overline{\theta},d+1-s )
\]
where $A (  \mathcal{Y},\theta )  $ denotes the conductor at
$\theta.$  The constant $\varepsilon(\mathcal{Y},\theta)$ is defined
independently of the conjecture.  By a theorem of Deligne and Langlands, after making certain
choices which we suppress in our notation, there is a product formula
\[
\varepsilon(  \mathcal{Y},\theta)  =%
{\textstyle\prod_{v}}
\varepsilon_{v}(  \mathcal{Y},\theta)
\]
where the product extends over the places of the rational field $\mathbf{Q}$
and almost all terms are 1. For details of the construction of the local constants
$\varepsilon
_{v} (  \mathcal{Y},\theta )  $ see [D2].  A discussion of local constants
which is well suited to the context of this paper is given in  [CEPT2].

Suppose now that the group $G$ acts tamely on $\mathcal{X}$ (for details see
condition (T1) below). A character of $G$ is
said to be symplectic when it is the character of representation of $G$ which
supports a non-degenerate $G$-invariant alternating form.  We shall show that for each symplectic character
$\theta$ of $G$, the $\theta$-component of the de Rham discriminant has a
natural square root called the Pfaffian.  One of the main
goals of this paper is to demonstrate an intimate relationship between the
sign of such Pfaffians and the archimedean constant $\varepsilon_{\infty
} (\mathcal{Y},\theta )$. In this regard we are most grateful to
Michael Atiyah who, early in our study, pointed out to us an analogous
phenomenon for Selberg $\zeta$-functions.  He showed in [Ms] that under suitable
circumstances, the constant in the functional equation of such a
$\zeta$-function is a signature invariant given by the regularised Pfaffian of
the Laplacian. This observation provided us with valuable insight into the
arithmetic situation which we consider here.

To describe a classical instance of our results, let  $N\subset\mathbf{C}$
be a number field 
which is tame and Galois over $\mathbf{Q}$ with Galois group $G$.
The simplest non-trivial case concerns quadratic $N$ over the rationals.
Let  $\sigma$ denote the
non-trivial automorphism of $N$. Since 2 is non-ramified in $N$, the ring
of integers of $N$ is equal to $\mathbf{Z}[  \alpha]  $ where
$\alpha=\frac{1}{2}\left(  1+\sqrt{d_{N}}\right)  $ and $d_{N}$ is the
discriminant of $N$. Let $\phi$ denote the non-trivial abelian character of
$G$; since $\phi$ is real-valued the character $2\phi$ is symplectic. Using
the fact that $\alpha-\sigma\left(  \alpha\right)  $ is a Pfaffian for the
trace form associated to $\phi$ (see Theorem 18, page 200 in [F1]), it follows
that the Pfaffian of the trace form associated to the symplectic character
$2\phi$ is equal to $\left(  \alpha-\sigma\left(  \alpha\right)  \right)
^{2}=d_{N}$, and this is of course positive or negative according as $N$ is
real or imaginary quadratic. We then remark that this sign coincides exactly
with the archimedean constant $\varepsilon_{\infty} (  \text{Spec
}(\mathbf{Z}),2\phi )  $, since $\varepsilon_{\infty} (\text{Spec}
(\mathbf{Z}),2\phi )  =\phi\left(  c\right)  $ where $c\in G$ is the image
of complex conjugation. This argument extends readily to arbitrary $G$ (see
Theorem 22 in [F2]).   The result is easiest to state when $N$ has
an integral normal basis.  Under this condition,   the
discriminant $d_{N}$ factors as a product of Pfaffians, and that for each
symplectic character $\psi$ of $G$ the sign of the Pfaffian associated to
$\psi$ is equal to $\varepsilon_{\infty}(  \text{Spec}(\mathbf{Z})%
,\psi )  $. Presently, we shall see in the main theorem that our work
on the signatures of de Rham discriminants provides a higher dimensional
version of this phenomenon.

The results of this paper are best seen as a complement to the work in [CPT2].
There we considered certain Euler characteristics of bounded complexes of
$G$-bundles which support metrics on the determinants of the isotypic
components of their cohomology.   We will refer to these Euler characteristics as
equivariant Arakelov classes.  In [CPT2]  we considered the equivariant Arakelov
classes associated to the de Rham complex, when the equivariant determinant of
cohomology was endowed with Quillen metrics arising from a choice of
K\"{a}hler metric.  In this approach all signature information was lost.
In [CPT1] we considered arithmetic surfaces in such a way as to retain the
signature information.  This was achieved by replacing the symplectic Arakelov
group by the somewhat sharper hermitian class group of Fr\"{o}hlich (see (5)
in 3.1.2 for a comparison of these two groups). The main result of this
article is to describe the de Rham discriminant of a suitable arithmetic
variety, in the hermitian class group, and thereby describe in full the
equivariant signature information of de Rham discriminant. In [CPT1] we were
also able to determine the non-archimedean local epsilon constant
$\varepsilon_{p} (  \mathcal{Y},\theta )  $ for a symplectic
character $\theta$ of $G$ for an arithmetic surface $\mathcal{Y}$, in terms of
the hermitian structure of the local de Rham complex. Here we do not attempt
to extend this result for non-archimedean local epsilon constants to higher
dimensions; however, we do remark that the work of T. Saito in [Sa] appears to
provide a good approach to this problem.

Prior to stating our main result, we first briefly comment on the methods and
tools used in this paper. One elementary tool, which is nonetheless crucial,
is the Pfaffian.  Our reformulation of the Pfaffian extends the description of
the discriminant of a non-degenerate bilinear form on a vector space as an
element of the dual  of the tensor square of the determinant of the vector
space. This approach extends readily to complexes, and thus one may then view
the discriminant as an element of the dual of the tensor square of the
determinant of the cohomology of the complex. For an alternating form on a
complex, the Pfaffian is a natural square root of the discriminant which is a
functional on the determinant of cohomology. We then recast Fr\"{o}hlich's
theory of hermitian classes in terms of the Pfaffian on determinants of
cohomology, and this approach turns out to be extremely well suited to the
calculation of  the de Rham discriminant.

We also wish to highlight the crucial use of the Atiyah-Singer $G$-signature
theorem. Recall that, when there is no group action, Hirzebruch gave a
beautiful formula for the signature of the cup-product form on the middle
dimensional cohomology group for a compact oriented real manifold $Z$ whose
dimension can be written $2d$ with $d$ an \textit{even} integer (see for
instance Theorem 6.6 in [AS]). Using the Index formula in the case where a
finite group $G$ acts on $Z,$ Atiyah and Singer obtained a formula for the
$G$-signature character of the virtual $\mathbf{R} [  G ]  $-module
\[
\text{H}_{B}^{d} (  Z,\mathbf{R} )  ^{+}-\text{H}_{B}^{d} (
Z,\mathbf{R} )  ^{-}%
\]
where ${\rm H}_{B}^{d} (  Z,\mathbf{R} )  ^{\pm}$ denotes a maximal
$\mathbf{R} [  G ]  $-module of the Betti cohomology group ${\rm H}_{B}%
^{d} (  Z,\mathbf{R} )  $ on which the cup-form product is positive
definite resp. negative definite.  In fact we shall need a modified
version of such invariants: namely, we let $\chi^{\pm}\left(  Z\right)  $
denote the dimension of ${\rm H}_{B}^{\bullet}\left(  Z,\mathbf{R}\right)  ^{\pm}$,
the (virtual) maximal space of the total Betti cohomology, on which the cup
product form is positive definite resp. negative definite (see 4.3 for
details).  Thus $\chi^{+}\left(  Z\right)  +\chi^{-}\left(  Z\right)  $  is
the Euler characteristic $\chi\left(  Z\right)  $. Then, for an arbitrary
compact oriented real manifold $Z$ of even dimension $2d$, we put
\[
\delta\left(  Z\right)  =\left\{
\begin{array}
[c]{c}%
 \chi\left(  Z\right)/2, \ \ \text{ if }d\text{ is odd},\\
\chi^{+}\left(  Z\right),  \text{ if }d\equiv2\operatorname{mod}4,\\
\chi^{-}\left(  Z\right),  \text{ if }d\equiv0\operatorname{mod}4.
\end{array}
\right.
\]

In order to present our main theorem we first need to introduce some notation.
Throughout $\mathcal{X}$ is a flat projective scheme over  ${\rm Spec} (
\mathbf{Z} )  $ which supports the action of a finite group $G;$ we let
$\mathcal{Y}$ denote the quotient $\mathcal{X}/G$, and we further assume that
the following two conditions are always satisfied:\smallskip

(T1) \ the action of $G$ on $\mathcal{X}$ is ``tame'' (for every point $x$ of
$\mathcal{X}$ the order of the inertia group $I_{x}\subset G$ is prime to the
residual characteristic of $x$).  Since $\mathcal{X}$ maps onto  ${\rm Spec}\left(
\mathbf{Z}\right)  $, it follows that   the locus of ramification
locus of the action of $G$ is fibral, and so, writing $X$ for the generic
fibre $\mathcal{X}\times_{\text{Spec}\left(  \mathbf{Z}\right)  }$Spec$\left(
\mathbf{Q}\right)  $, $G$ acts freely on $X$ and so the cover $X\rightarrow Y$
is etale;

(T2) \ both schemes $\mathcal{X}$ and $\mathcal{Y}\ $are regular and ``tame''
(i.e. they are regular and all their special fibres are divisors with normal
crossings with multiplicities prime to the residue characteristic).\medskip

Now let $\Omega_{\mathcal{X}/\mathbf{Z}}^{1}$ denote the coherent sheaf of
differentials of $\mathcal{X}\to {\rm Spec}(\mathbf{Z})$. Since $\mathcal{X}$ is regular, we
may choose a resolution of $\Omega_{\mathcal{X}/\mathbf{Z}}^{1}$ by a length
two complex $K^{\bullet}$ of $G$-equivariant locally free $O_{\mathcal{X}%
}$-sheaves; for
$i\geq0$ we let $L\wedge^{i} $ denote the $i$-th left derived exterior power
functor of Dold-Puppe on perfect complexes of $G$-equivariant $O_{\mathcal{X}%
}$-sheaves (that is to say, $O_{\mathcal{X}}$-sheaves with a $G$-action which
is compatible with the $G$-action on $O_{\mathcal{X}}$). Thus $L\wedge^{i}$
$K^{\bullet}$ denotes the complex arising from the application of $L\wedge
^{i}$ to $K^{\bullet} $ and we define $L\wedge^{\bullet}\Omega_{\mathcal{X}%
/\mathbf{Z}}^{1}$  to be the direct sum of the complexes $L\wedge^{i}$
$K^{\bullet} [  -i ]  $ for $0\leq i\leq d$. For details of the Dold-Puppe exterior power
functor the reader is referred to [DP], [I] and to 5.4-5.9 in [So]. We recall
from [CEPT1] that, because $G$ acts tamely, $R\Gamma (  \mathcal{X}%
, L\wedge^{\bullet}\Omega_{\mathcal{X}/\mathbf{Z}}^{1} )  $ may
represented by a perfect $\mathbf{Z}\left[  G\right]  $-complex. Note for
future reference that on the generic fibre $X$ of $\mathcal{X}$ each $\left(
L\wedge^{i}K^{\bullet}\right)  \otimes_{\mathbf{Z}}\mathbf{Q}$ is
quasi-isomorphic to the sheaf of differentials $\Omega_{X\mathbf{/Q}}^{i}$
viewed as a complex concentrated in degree zero.

In Sect. 4 we shall recall in detail from [CPT1] the symmetric $G$-invariant
pairings on the cohomology groups
\[
\sigma_{X}^{t}:\text{H}^{t}R\Gamma (  X,L\wedge
^{\bullet}\Omega_{X/\mathbf{Q}}^{1} )   [  d ]  \times
\text{H}^{-t}R\Gamma ( X, L\wedge^{\bullet}\Omega
_{X/\mathbf{Q}}^{1} )   [  d ]  \rightarrow\mathbf{Q}%
\]
arising from Serre duality. In [F1] Fr\"{o}hlich showed how to use the notion
of Pfaffian to construct a refined discriminant, or hermitian class, for any
locally free $\mathbf{Z} [  G ]  $-module which supports a
non-degenerate $G$-invariant symmetric form over $\mathbf{Q}$. In [CPT1] we
extended this construction to perfect $\mathbf{Z} [  G ]  $-complexes
with non-degenerate $G$-invariant symmetric forms on the cohomology of the
complex tensored by $\mathbf{Q}$. Thus, to the pair \hbox{$(  R\Gamma(
\mathcal{X},L\wedge^{\bullet}\Omega_{\mathcal{X}/\mathbf{Z}}^{1})
 [  d ]  ,\sigma_X)  $}, we may associate a
so-called hermitian Euler characteristic \hbox{$\chi_{\text{H}}^{\text{s}}(
R\Gamma(  \mathcal{X},L\wedge^{\bullet}\Omega_{\mathcal{X}/\mathbf{Z}%
}^{1})   [  d ]  ,\sigma_{X})  $} which takes
values in the hermitian class group ${\rm H}^{\text{s}} (  \mathbf{Z} [
G ]   )  $.

The hermitian Euler characteristic $\chi_{\text{H}}^{\text{s}} (
R\Gamma (  \mathcal{X},L\wedge^{\bullet}\Omega_{\mathcal{X}/\mathbf{Z}%
}^{1} )   [  d ]  ,\sigma_X )  $ was completely
determined in [CPT1] when $\mathcal{X}$  is an arithmetic surface. In this
paper  we shall essentially determine the hermitian Euler
characteristic  $\chi_{\text{H}}^{\text{s}} (  R\Gamma (
\mathcal{X},L\wedge^{\bullet}\Omega_{\mathcal{X}/\mathbf{Z}}^{1} )
 [  d ]  ,\sigma_X ) $ for arbitrary fibral
dimension $d $. To be a little more precise, we shall show that symplectic
hermitian Euler characteristics decompose into the product of a metric
invariant and a signature invariant. Writing $R_{G}^{\text{s}}$ for the group
of virtual symplectic characters of $G$, we shall see that, under the above
mentioned decomposition into metric and signature invariants (after tame
extension of coefficients), the image of  $\chi_{\text{H}}^{\text{s}} (
R\Gamma (  \mathcal{X},L\wedge^{\bullet}\Omega_{\mathcal{X}/\mathbf{Z}%
}^{1} )  [  d ]  ,\sigma_X )  $ lies in a
group which is naturally isomorphic to
\[
\text{Hom}_{\text{Gal}(  \overline{\mathbf{Q}}/\mathbf{Q})
} (  R_{G}^{\text{s}},\mathbf{Q}^{\times} )  \times\text{Hom}\left(
R_{G}^{\text{s}},\mathbf{\pm}1\right)  .
\]
We let $\chi_{1}$ resp. $\chi_{2}$ denote the image of $\chi_{\text{H}%
}^{\text{s}}(  R\Gamma(  \mathcal{X},L\wedge^{\bullet}%
\Omega_{\mathcal{X}/\mathbf{Z}}^{1})  \left[  d\right]  ,\sigma
_{\mathcal{X}})  $ in the first resp. second component. We shall see
that $ \chi_{1}$ coincides with the equivariant Arakelov class, which was
fully described in [CPT2]; there we saw that this class is given by the
$\varepsilon$-constant homomorphism which, for virtual characters $\theta$ of
degree zero, maps $\theta$ to $\varepsilon (  \mathcal{Y},\theta )
$. This product decomposition is of fundamental importance, for we note there
are two natural sign invariants: as indicated above, the first sign invariant
\ determines (conjecturally) the symmetry or skew-symmetry of the functional
equation of the Artin-Hasse-Weil L-function; whereas the second sign invariant
should be thought of as the archimedean signature. Such a double appearance of
sign invariants was apparent in the work of  Fr\"{o}hlich  (see for instance
Corollary 3 page 192 in [F1]). The essential contribution of this article is
the following evaluation of the signature class $\chi_{2}$.

\begin{theorem}
For \ a symplectic character $\theta$ of $G$
\[
\chi_{2}\left(  \theta\right)  =\left(  -1\right)  ^{\delta\left(  Y\right)
\theta\left(  1\right)  /2}\varepsilon_{\infty}\left(  \mathcal{Y}%
,\theta\right)
\]
where $\varepsilon_{\infty}\left(  \mathcal{Y},\theta\right)  $ is the
archimedean constant described in the first part of the Introduction.
\end{theorem}

Note that it is a remarkable fact that the signatures of such equivariant de
Rham discriminants, which come from de Rham cohomology, are determined by the
archimedean $\varepsilon$-constants (at least for virtual characters of degree
zero) which derive from the Hodge realisation of the real Artin motives
$X_{\mathbf{R}}\otimes_{G}V_{\theta}$ (see 5.3 in [D2], and see also Section 5
of [CEPT2]).

We conclude our introduction by providing a brief overview of the structure of
this paper. The basic definitions and results on Pfaffians are all presented
in Sect. 2. Then in Sect. 3 we introduce the symplectic hermitian classgroup
${\rm H}^{\text{s}}\left(  \mathbf{Z}\left[  G\right]  \right)  $ where we define
our hermitian Euler characteristics. It should be noted that this classgroup
is in fact slightly different from the hermitian classgroup defined by
Fr\"{o}hlich; we prefer to work with this version, because it contains the
symplectic equivariant Arakelov class group as a subgroup in a natural way. In
the Appendix we describe the natural map from Fr\"{o}hlich's classgroup to our
classgroup ${\rm H}^{\text{s}} (  \mathbf{Z} [  G ]  )  $. In the
Appendix, we shall also show how the hermitian Euler characteristics defined
in [CPT1], via lifts of pairings on cohomology to the whole perfect complex,
agree with the hermitian Euler characteristics that we use this paper, which
are defined via the Pfaffian on cohomology.

Finally in Sect. 4 we apply the foregoing theory to our arithmetic situation,
and we consider the de Rham discriminants: namely, the hermitian Euler
characteristics of the de Rham complex with forms on cohomology arising from
the natural duality pairings. Here we recall some basic results on archimedean
$\varepsilon$-constants from [CEPT2] and we study the signature properties of
de Rham cohomology. This then provides us with all the tools we need to
complete the proof of the main theorem in the final sub-section.

\section{PFAFFIANS.}

In this section we work over an arbitrary field $K$ of characteristic zero.
All vector spaces are assumed to be finite dimensional and all bilinear forms
are assumed to be non-degenerate.

\subsection{ DISCRIMINANTS AND PFAFFIANS.}

\subsubsection{Determinants.}

For a $K$-vector space $V$ we let $V^{D}$ denote the $K$-linear dual
Hom$_{K}\left(  V,K\right)  $, and if $V$ has dimension $d,$ we write
$\det\left(  V\right)  =%
{\textstyle\bigwedge^{d}}
V$ and let $\widehat{\det}\left(  V\right)$ denote the graded line
$\left(  \det\left(  V\right)  ,\dim V\right)  $. For a graded $K$-line
$\left(  L,n\right)  $, we put $\left(  L,n\right)  ^{D}=\left(
L^{D},-n\right)  $. In the sequel we shall often write $L^{-1}$ for $L^{D}$
and $\left(  L,n\right)  ^{-1}$ for $\left(  L,n\right)  ^{D}$.

Throughout this article we shall adopt the following convention in our use of
exterior products: we follow Deligne and we normalise the ``twist''
isomorphism between the tensor product of the determinants of two vector
spaces as follows: given two finite dimensional vector spaces $V,W$ over $K
$, the tensor product of the graded lines $\widehat{\det}\left(  V\right)  $
and $\widehat{\det}\left(  W\right)  $ is
\[
\widehat{\det}\left(  V\right)  \otimes\widehat{\det}\left(  W\right)
=\left(  \det\left(  V\right)  \otimes_{K}\det\left(  W\right)  ,\dim\left(
V\right)  +\dim\left(  W\right)  \right)
\]
and we twist the standard isomorphism $\det\left(  V\right)  \otimes
\det\left(  W\right)  \cong\det\left(  W\right)  \otimes\det\left(  V\right)
$ according to the Koszul rule of signs, i.e. by the factor $\left(
-1\right)  ^{\dim\left(  V\right)  \dim\left(  W\right)  }$; thus, to be
absolutely explicit, under the new isomorphism
\[
\left(  v_{1}\wedge\cdot\cdot\cdot\wedge v_{n}\right)  \otimes\left(
w_{1}\wedge\cdot\cdot\cdot\wedge w_{m}\right)  \leftrightarrow\left(
-1\right)  ^{\dim\left(  V\right)  \dim\left(  W\right)  }\left(  w_{1}%
\wedge\cdot\cdot\cdot\wedge w_{m}\right)  \otimes\left(  v_{1}\wedge\cdot
\cdot\cdot\wedge v_{n}\right)  .
\]
With this convention we then see that both the following diagram
\[%
\begin{array}
[c]{ccc}%
\widehat{\det}\left(  V\right)  \otimes\widehat{\det}\left(  W\right)  &
\rightarrow & \widehat{\det}\left(  V\oplus W\right) \\
\downarrow &  & \downarrow\\
\widehat{\det}\left(  W\right)  \otimes\widehat{\det}\left(  V\right)  &
\rightarrow & \widehat{\det}\left(  W\oplus V\right)
\end{array}
\]
and the corresponding diagram where we then forget the grading commute. Here
the horizontal maps are the maps induced by the isomorphisms $\det\left(
V\right)  \otimes\det\left(  W\right)  \cong\det\left(  V\oplus W\right)  $
and the above description of $\widehat{\det}\left(  V\right)  \otimes
\widehat{\det}\left(  W\right)  $; the right-hand vertical arrow is induced by
the natural isomorphism
\[
V\oplus W\cong W\oplus V\bigskip,\;\;v\oplus w\longmapsto w\oplus v\ ;
\]
the left-hand vertical arrow is the above ``twisted '' isomorphism. This
convention will help us avoid what Deligne calls the ``nightmare of signs''.

\bigskip

\subsubsection{\bigskip Pfaffians.}

We begin by recalling the notion of discriminant for a non-degenerate bilinear
form $h$ on $V$. Thus such a form $h$ affords an isomorphism $h:V\rightarrow
V^{D}$, via the rule $h\left(  x\right)  \left(  y\right)  =h\left(
y,x\right)  $. The discriminant $d_{h}$ is then defined to be the linear
isomorphism of one dimensional $K$ vector spaces
\[
d_{h}:\det\left(  V\right)  ^{\otimes2}\underset{1\otimes\det\left(  h\right)
}{\rightarrow}\det\left(  V\right)  \otimes\det\left(  V^{D}\right)
\rightarrow K
\]
given by using $\det(  V^{D})  \cong\det(  V)  ^{D}$
and contraction.

\medskip Suppose now that $h$ is an \textit{alternating form} and let
$\dim\left(  V\right)  =2n$; recall that $V$ has a hyperbolic basis $\left\{
u_{1},u_{1}^{\prime},u_{2},u_{2}^{\prime},\ldots,u_{n},u_{n}^{\prime
}\right\}  $ where
\[
h\left(  u_{i},u_{j}\right)  =0=h\left(  u_{i}^{\prime},u_{j}^{\prime}\right)
\text{ for all }i,j,\text{ and }h\left(  u_{i},u_{j}^{\prime}\right)
=\delta_{ij}.
\]
Thus in particular $Ku_{i}^{\prime}$ identifies, via $h$, as the dual line of
$Ku_{i}$, and so we can define
\begin{align*}
\text{Pf}_{h}  & :\det (  V )  =\otimes_{i=1}^{n}\det (
Ku_{i}\oplus Ku_{i}^{\prime} ) \\
& =\otimes_{i=1}^{n}\det (  Ku_{i}\oplus (  Ku_{i} )
^{D} )  \rightarrow\otimes_{i=1}^{n}K=K\ .
\end{align*}
Alternatively we see that ${\rm Pf}_{h}\;$is the unique $K$-linear functional on
$\det (  V )  $ such that
\[
\text{Pf}_{h}\left(  u_{1}\wedge u_{1}^{\prime}\wedge u_{2}\wedge
u_{2}^{\prime}\wedge\cdots\wedge u_{n}\wedge u_{n}^{\prime}\right)  =1.
\]
(Note that had we defined $h\left(  x\right)  \left(  y\right)  =h\left(
x,y\right)  $, then we would have $1={\rm Pf}_{h}\left(  u_{1}^{\prime
}\wedge u_{1}\cdots\right)  $.)

This notation suggests that in fact ${\rm Pf}_{h}$ does not depend on the choice of
particular hyperbolic basis. That this is indeed the case follows from the
first of the following two lemmas, both of whose proofs are routine.

\begin{lemma}
If $\left\{  v_{1},v_{1}^{\prime},v_{2},v_{2}^{\prime},\ldots,v_{n}%
,v_{n}^{\prime}\right\}  $ is a further hyperbolic basis of $V$, with respect
to $h$ then, since a symplectic automorphism has determinant 1,
\[
\text{\rm Pf}_{h}\left(  v_{1}\wedge v_{1}^{\prime}\wedge v_{2}\wedge
v_{2}^{\prime}\cdots\wedge v_{n}\wedge v_{n}^{\prime}\right)  =1.
\]
\end{lemma}

\begin{lemma}
{\textit For $i=1,2$ let $h_{i}$ be an alternating-form on the
vector space $V_{i}$. Let $h_{1}\oplus h_{2}$  denote the
orthogonal sum form on $V_{1}\oplus V_{2}$.  Then, with the above
convention, ${\rm Pf}_{h_{1}\oplus h_{2}} ={\rm Pf}_{h_{1}}\otimes{\rm Pf}_{h_{2}}%
$  under the identification $\det\left(  V_{1}\oplus V_{2}\right)
=\det\left(  V_{1}\right)  \otimes\det\left(  V_{2}\right)  $.}
\end{lemma}

The following lemmas describe the functorial properties of the Pfaffian which
we shall require. The proofs are all completely routine and follow from the
standard properties of determinants.

\begin{lemma}
\textit{For an alternating form }$h$\textit{\ on }$V$\textit{\ and for a given
isomorphism of }$K$-vector spaces \textit{\ }$\phi:V\rightarrow W,$%
\textit{\ let }$\phi^{\ast}h$\textit{\ denote the form on }$W$\textit{\ given
by the rule }
\[
\phi^{\ast}h\left(  x,y\right)  =h\left(  \phi^{-1}x,\phi^{-1}y\right)  .
\]
\textit{Then the following diagram commutes }
\[%
\begin{array}
[c]{lll}%
\det\left(  V\right)  & \overset{\text{\rm Pf}_{h}}{\longrightarrow} & K\\
\downarrow\det\left(  \phi\right)  &  & \downarrow\\
\det\left(  W\right)  & \overset{\text{\rm Pf}_{\phi^{\ast}h}}{\longrightarrow} & K.
\end{array}
\]
\textit{In particular if }$V=W$\textit{\medskip}$,$\textit{\ then }%
${\rm Pf}_{\phi^{\ast}h}=\det\left(  \phi\right)  ^{-1}{\rm Pf}_{h}$.
\end{lemma}

\begin{proposition}
\textit{For a given alternating form $h$ on $V$  and for an
automorphism $A$ of $V$, let $\widehat{A}$
 denote the adjoint of $A$ with respect to $h$; that is to say
$h\left(  Ax,y\right)  =h (  x,\widehat{A}y )$.  Suppose
$A$ is self-adjoint, so that $A=\widehat{A}$, and define
$h'\left(x,y\right ) = h\left (Ax,y\right )$.  Then there
is an automorphism $B$ of $V$ such that $h=B^{\ast}h^{\prime}$. This implies
that $A=\widehat{B}B$ and by the above ${\rm Pf}_{h^{\prime}}=\det (  B )
{\rm Pf}_{h}$.  The value $\det (  B )  $
 therefore depends only on $A$ and we call it the Pfaffian
of $A$, denoted $\mathbf{pf} (  A )  ,$ so that we have
\[
{\rm Pf}_{h^{\prime}}=\mathbf{pf} (  A )  {\rm Pf}_{h}.
\]}
\end{proposition}

\bigskip

\begin{remark}
{\rm In the sequel Pf will denote a functional on a $K$-line, whereas $\mathbf{pf}
$ will denote the Pfaffian of a matrix.}
\end{remark}

\subsection{ EXTENSION TO COMPLEXES.}

Let $C^{\bullet}$ denote a bounded complex of vector spaces over a field $K$.
We put
\[
C^{\text{ev}}=C^{0}%
{\textstyle\bigoplus_{i>0}}
 (  C^{2i}\oplus C^{-2i} )  \text{ \ \ and \ \ \ }C^{\text{odd}}=%
{\textstyle\bigoplus_{i\geq0}}
 (  C^{2i+1}\oplus C^{-2i-1} )  .
\]
and we recall that $\det (  C^{\bullet} )  =\otimes\det (
C^{i} )  ^{ (  -1 )  ^{i}}.$

There is a natural map (given by reordering)
\[
\upsilon_{C^{\bullet}}:\det (  C^{\bullet} )  \rightarrow\det (
C^{\text{ev}} )  \otimes\det (  C^{\text{odd}} )  ^{-1}%
\]
where in full the latter line is
\[
\det (  C^{0} )  \otimes\det (  C^{2} )  \otimes\det (
C^{-2} )  \otimes\cdots\otimes\det (  C^{1} )
^{-1}\otimes\det   (C^{-1} )  ^{-1}\otimes\cdots
\]

\begin{remark}{\rm 
\bigskip(a) If $D^{\bullet}$ is a further $K$-complex and if all the
terms of $C^{\bullet}$ and $D^{\bullet}$ have even dimension, then the map
\[
\det(  C^{\bullet}\oplus D^{\bullet})  \cong\det  (C^{\bullet
})  \otimes\det(  D^{\bullet})
\]
given by using the Koszul-twist isomorphisms coincides with the naive map
given by the reordering of terms.

(b) If  again all the terms $C^{i}$ have even dimension, then the map
\[
\det(  C^{\bullet})  \rightarrow\det(  C^{\text{ev}} )
\otimes\det(  C^{\text{odd}})  ^{-1}
\]
given by using the Koszul-twist isomorphisms coincides with the naive map
$\upsilon_{C^{\bullet}}$ given by the reordering of terms.}
\end{remark}

We shall write ${\rm H}^{\bullet} (  C^{\bullet} )  $ for the complex
$\{ {\rm H}^{i}\left(  C^{\bullet}\right)\}_i  $, with zero boundary maps. As
above we write
\[
\text{H}^{\text{ev}}=\text{H}^{\text{ev}}   (C^{\bullet} )
=\text{H}^{0} (  C^{\bullet} )
{\textstyle\bigoplus_{\dot{i}>0}}
(\text{H}^{2i} (  C^{\bullet} )  \oplus\text{H}^{-2i} (
C^{\bullet} ))
\]
\ and \
\[
\text{H}^{\text{odd}}=\text{H}^{\text{odd}} (  C^{\bullet} )  =%
{\textstyle\bigoplus_{i\geq0}}
(\text{H}^{2i+1} (  C^{\bullet} )  \oplus\text{H}^{-2i-1} (
C^{\bullet} ) ) .
\]
From [KM] we recall that there is a canonical isomorphism of $K$-lines
\[
\xi:\det (  C^{\bullet} )  \cong\det (  \text{H}^{\bullet} (
C^{\bullet} )   )  .
\]

\begin{definition}
Suppose we are given \textit{alternating} forms $h^{\text{\rm ev}}$ on
${\rm H}^{\text{\rm ev}}$ and $h^{\text{\rm odd}}$ on ${\rm H}^{\text{\rm odd}}$. Define
 ${\rm Pf}_{h }$ to be the element of the dual of the line
$\det\left(  C^{\bullet}\right)  ,\;$%
\[
\text{\Pf}_{h}:\det\left(  C^{\bullet}\right)  \rightarrow K\mathbf{,}%
\]
given by composing
\[
\text{\rm Pf}_{h}=\text{\rm Pf}_{h^{\text{\rm ev}}}\otimes\text{\rm Pf}_{h^{\text{\rm odd}}}%
^{-1}:\det (  \text{\rm H}^{\text{\rm ev}} (  C^{\bullet} )   )
\otimes\det (  \text{\rm H}^{\text{\rm odd}} (  C^{\bullet} )   )
^{-1}\rightarrow K
\]
with the isomorphism $\upsilon_{\text{H}^{\bullet}\left(  C^{\bullet}\right)
}\circ\xi.$ Note that in the sequel, for brevity, we shall usually write $h$
for the pair $\left\{  h^{\text{\rm ev}},h^{\text{\rm odd}}\right\}  .$
\end{definition}

\subsection{EQUIVARIANT PFAFFIANS.}

Suppose now that $G$ is a finite group, $K$ is a subfield of the real numbers
$\mathbf{R}$, and that $W$ is a symplectic complex representation of $G$ with
character $\theta$; thus, by definition, $W$ supports a non-degenerate
$G$-invariant alternating-form $\kappa$.

\begin{lemma}
\textit{If }$\kappa^{\prime}$\textit{\ is a further such form on }$W$\textit{,
then, since every pair of non-degenerate alternating forms on }$W$ are
isomorphic, there is an automorphism $B$ of $W\;$such that $\kappa^{\prime
}=B^{\ast}\kappa.$ Since both forms are $G$-invariant,
the\textit{\ self-adjoint automorphism }$A=\widehat{B}B$\textit{\ \ is a }%
$G$-automorphism of $W.$
\end{lemma}

\begin{definition}
A symmetric $K [  G ]  $-complex is a pair $   (C^{\bullet
},\sigma )  $ where $C^{\bullet}$ is a perfect $K [  G ]
$-complex and where $\sigma^{\text{\rm ev}}$ and $\sigma^{\text{\rm odd}}$ are
non-degenerate real-valued $G$-invariant \textit{symmetric} forms on
${\rm H}^{\text{\rm ev}} (  C^{\bullet} )  $ and ${\rm H}^{\text{\rm odd}} (
C^{\bullet} )  $ respectively.
\end{definition}

For a given symmetric complex $ (  C^{\bullet},\sigma )  $ and for
$W$ and $\kappa$ as above, we define $\det (  C_{W}^{\bullet} )  $ to
be the line $\det (   (  C^{\bullet}\otimes_{\mathbf{R}}W )
^{G} )  $; thus we have the canonical isomorphism
\[
\xi_{W}:\det (  C_{W}^{\bullet} )  \cong\det (  \text{\rm H}%
^{\bullet} (  C^{\bullet} )  _{W} )  .
\]
By restricting $\sigma^{\text{\rm ev}}\otimes\kappa$ to $ (  \text{\rm H}%
^{\text{\rm ev}}\otimes W )  ^{G}$ we obtain a non-degenerate alternating
form which we denote by $ (  \sigma^{\text{ev}}\otimes\kappa )
^{G};$ similarly we obtain a form $ (  \sigma^{\text{odd}}\otimes
\kappa )  ^{G}$ on $ (  \text{\rm H}^{\text{\rm odd}}\otimes W )  ^{G}.$
Thus we obtain the composite map
\[
\det (  C_{W}^{\bullet} )  \underset{\xi_{W}}{\cong}\det (
\text{\rm H}^{\bullet} (  C^{\bullet} )  _{W} )  \underset
{\upsilon_{H^{\bullet}}}{\cong}\det (  \text{\rm H}^{\text{\rm ev}} (Ä
C^{\bullet} )  _{W} )  \otimes\det (  \text{\rm H}^{\text{\rm odd}%
} (  C^{\bullet} )  _{W} )  ^{-1}\rightarrow K
\]
where the right hand arrow is ${\rm Pf}_{ (  \sigma\otimes\kappa )  ^{G}}
$.

\section{CLASS GROUPS.}

\subsection{HERMITIAN AND ARAKELOV CLASSGROUPS.}

In this subsection we give the definition of the symplectic hermitian class
group, and we also briefly recall the definition of the equivariant Arakelov
classgroup - for full details on the latter see [CPT2].

\subsubsection{Definition of classgroups.}

Let $R_{G}$ denote the group of complex virtual characters of $G$, and let
$R_{G}^{\text{s}}$ be the subgroup of virtual symplectic characters.  Let
$\overline{\mathbf{Q}}$ be an algebraic closure of $\mathbf{Q}$ in
$\mathbf{C}$, and define $\Omega={\rm Gal} (  \overline{\mathbf{Q}
}/\mathbf{Q} )  $.  Define $J_{f}$ (resp. $J_{\infty}$) to be the group of
finite ideles (resp. the archimedean ideles) of $\overline{\mathbf{Q}}$.  Thus
$J_{f}$ is the direct limit of the finite idele groups of all algebraic
number fields $E$ in $\overline{\mathbf{Q}}$, and
\[
J_{\infty}=\lim_{E\subset\overline{\mathbf{Q}}}\left(  E\otimes_{\mathbf{Q}%
}\mathbf{R}\right)  ^{\times}.
\]
The idele group of $\overline
{\mathbf{Q}}$ is  $J=J_{f}\times J_{\infty}$.

Let $\widehat{\mathbf{Z}}=\prod_{p}\mathbf{Z}_{p}$ denote the ring of integral
finite ideles of $\mathbf{Z}$. For $x\in$ $\widehat{\mathbf{Z}} [
G ]  ^{\times}$, the element ${\rm Det} (  x )  \in
{\rm Hom}_{\Omega} (  R_{G},J_{f} )  $ is defined by the rule that for
a representation $T$ of $G$ with character $\psi$
\[
\text{Det}\left(  x\right)  \left(  \psi\right)  =\det\left(  T\left(
x\right)  \right)  ;
\]
the group of all such homomorphisms is denoted by
\[
\text{Det} (  \widehat{\mathbf{Z}} [  G ]  ^{\times} )
\subseteq\text{Hom}_{\Omega}\left(  R_{G},J_{f}\right)  .
\]
More generally, for $n>1$ we can form the group ${\rm Det} (  GL_{n} (
\widehat{\mathbf{Z}} [  G ]   )  )  ;$ as each group ring
$\mathbf{Z}_{p} [  G ]  $ is semi-local, we have the equality
${\rm Det} (  GL_{n} (  \widehat{\mathbf{Z}} [  G ]   )
 )  ={\rm Det} (  \widehat{\mathbf{Z}} [  G ]  ^{\times} )
$ (see 1.2.6 in [T]).

Recall that by the Hasse-Schilling norm theorem
\begin{equation}
\text{Det} (  \mathbf{Q} [  G ]  ^{\times} )  =\text{Hom}%
_{\Omega}^{+} (  R_{G},\overline{\mathbf{Q}}^{\times} )
\end{equation}
where the right-hand expression denotes Galois equivariant homomorphisms whose
values on $R_{G}^{s}$ are all totally positive. We then have a diagonal map
\[
\Delta:\text{Hom}_{\Omega}^{+} (  R_{G},\overline{\mathbf{Q}}^{\times
} )  \rightarrow\text{Hom}_{\Omega} (  R_{G},J_{f} )
\times\text{Hom} (  R_{G},\mathbf{R}_{>0} )
\]
where $\Delta (  f )  =f\times\left|  f\right|  .$ Given a
homomorphism $f$ on $R_{G},$ we shall write $f^{\text{s}}$ for the restriction
of $f$ to $R_{G}^{\text{s}}$; in particular\ we write
\[
\Delta^{\text{s}}:\text{Hom}_{\Omega}^{+} (  R_{G}^{\text{s}}
,\overline{\mathbf{Q}}^{\times} )  \rightarrow\text{Hom}_{\Omega} (
R_{G}^{\text{s}},J_{f} )  \times\text{Hom} (  R_{G}^{\text{s}
},\mathbf{R}_{>0} )
\]
for the restriction of $\Delta$ to $R_{G}^{\text{s}}$, so  that
\[
\Delta^{\text{s}} (  f^{\prime} )  =f^{\prime}\times\left|
f^{\prime}\right|  =f^{\prime}\times f^{\prime}.\medskip
\]

\begin{definition}
The group of \textit{symplectic hermitian classes} ${\rm H}^{\text{\rm s}} (
\mathbf{Z}[G]   )  $ is defined to be the quotient group
\begin{equation}
\text{\rm H}^{\text{\rm s}}\left(  \mathbf{Z}[G]  \right)  =\frac
{\text{\rm Hom}_{\Omega}\left(  R_{G}^{\text{s}},J_{f}\right)  \times
\text{\rm Hom}\left(  R_{G}^{\text{s}},\mathbf{R}^{\times}\right)  }
{\operatorname{Im}(\Delta^{\text{\rm s}})\cdot  (  \text{\rm Det}^{\text{\rm s}} (
\widehat{\mathbf{Z}} [G]  ^{\times} )  \times1 )  }
\end{equation}
where ${\rm Det}^{\text{\rm s}} (  \widehat{\mathbf{Z}}[G]  ^{\times
} )  $ denotes the restriction of ${\rm Det} (  \widehat{\mathbf{Z}}\left[
G\right]  ^{\times} )  $ to $R_{G}^{\text{\rm s}}$. Note that this hermitian
classgroup ${\rm H}^{\text{\rm s}} (  \mathbf{Z}\left[  G\right]   )  $ is
slightly different from the hermitian classgroup ${\rm HCl} (  \mathbf{Z}[G]   )  $ 
used in [CPT1] and [F1]. There is a natural map between
these two classgroups. For details see the Appendix.
\end{definition}

\bigskip

We recall from Definition 3.2 in [CPT2] that the group of Arakelov classes is
defined to be
\begin{equation}
\text{\rm A} (  \mathbf{Z}[G]  )  =\frac{\text{\rm Hom}
_{\Omega} (  R_{G},J_{f} )  \times\text{\rm Hom} (  R_{G}%
,\mathbf{R}_{>0} )  }{\operatorname{\rm Im}(\Delta)\cdot  (  \text{\rm Det} (
\widehat{\mathbf{Z}}[G]  ^{\times} )  \times 1 )  }
\end{equation}
and that the group of symplectic Arakelov classes (see Definition 4.1 in
[CPT2]) is defined to be
\begin{equation}
\text{\rm A}^{\text{\rm s}} (  \mathbf{Z}[G]  )  =\frac
{\text{\rm Hom}_{\Omega} (  R_{G}^{\text{\rm s}},J_{f} )  \times
\text{\rm Hom} (  R_{G}^{\text{\rm s}},\mathbf{R}_{>0} )  }{\operatorname{\rm Im}
(\Delta^{\text{\rm s}})\cdot (  \text{\rm Det}^{\text{\rm s}} (  \widehat{\mathbf{Z}%
}[G] ^{\times} )  \times 1 )  }.
\end{equation}

\begin{remark}{\rm 
Firstly, from the above descriptions, we see that ${\rm A}^{\text{\rm s}} (
\mathbf{Z}[G]  )  $ is naturally a subgroup of
${\rm H}^{\text{\rm s}} (  \mathbf{Z}[G] )  $. Secondly, from
Lemma 2.1 on page 60 of [F2], we note that, since all symplectic characters
are real-valued, there is a natural isomorphism induced by the inclusion
$\overline{\mathbf{Q}}\subset\mathbf{C}$
\[
\text{\rm Hom}_{\Omega} (  R_{G}^{\text{\rm s}},J_{\infty} )  \cong
\text{\rm Hom} (  R_{G}^{\text{\rm s}},\mathbf{R}^{\times} )  .
\]}
\end{remark}

\textbf{\bigskip}

\subsubsection{Rational classes and signature classes.}

Let $-1_{\infty}$ denote the idele which is 1 at all finite primes and
which is $-1$ at all infinite primes.  We then consider the two subgroups of
\begin{align*}
\text{Hom}_{\Omega} (  R_{G}^{\text{s}},J )   & =\text{Hom}_{\Omega
} (  R_{G}^{\text{s}},J_{f} )  \times\text{Hom}_{\Omega} (
R_{G}^{\text{s}},J_{\infty} ) \\
& \cong\text{Hom}_{\Omega} (  R_{G}^{\text{s}},J_{f} )
\times\text{Hom} (  R_{G}^{\text{s}},\mathbf{R}^{\times} )
\end{align*}
given by
\[
\text{R} (  \mathbf{Z}[G]   )  =\text{Hom}_{\Omega
} (  R_{G}^{\text{s}},\mathbf{Q}^{\times} )  \times1,
\]%
\[
\text{S}_{\infty} (  \mathbf{Z}[G]  )  =1\times
\text{Hom} (  R_{G}^{\text{s}},\pm1 )  =\text{Hom} (
R_{G}^{\text{s}},\pm1_{\infty} )  .
\]

\begin{theorem}
The natural map from ${\rm Hom}_{\Omega} (  R_{G}^{\text{\rm s}},J )  $ to
${\rm H}^{\text{\rm s}} (  \mathbf{Z}[G]   )  $ induces an
injection on ${\rm R}(  \mathbf{Z}[G]  )  \times
{\rm S}_{\infty}\left(  \mathbf{Z}[G]  \right)$; thus, in the
sequel, we shall view ${\rm R} (  \mathbf{Z}[G]   )  \times
{\rm S}_{\infty} (  \mathbf{Z}[G]  )  $ as a subgroup of
${\rm H}^{\text{s}} (  \mathbf{Z}[G] )  .$
\end{theorem}

\noindent {\sl Proof.} Let $r\times s\in{\rm R} (  \mathbf{Z}[G]  )
\times{\rm S}_{\infty} (  \mathbf{Z}[G]  )  $. We must show
that if $r\times s\in\operatorname{Im}(\Delta^{\text{s}})\cdot (  \text{Det}
^{\text{s}} (  \widehat{\mathbf{Z}}[G] ^{\times} )
\times1 )  $, then $r=1=s$. Now by the Hasse-Schilling theorem we see
immediately that $s$ is positive and hence 1.  We therefore deduce that
$r\in{\rm R} (  \mathbf{Z}[G] )  \cap{\rm Det}^{\text{s}
} (  \widehat{\mathbf{Z}}[G]  ^{\times} )  $ which is
known to be trivial by Proposition 6.1 in [CNT] (see also [F1] Theorem 17, p. 190). \hfill $\square $
\medskip

The counterpart for Arakelov classes is the following result, which is shown
in 4.D of [CPT2]:

\begin{theorem}
The natural map from ${\rm Hom}_{\Omega} (  R_{G}^{\text{\rm s}},J_{f} )
\times{\rm Hom} (  R_{G}^{\text{\rm s}},\mathbf{R}_{>0} )$ to
${\rm A}^{\text{\rm s}} (  \mathbf{Z}[G]   )  $ induces an
injection on ${\rm R}(  \mathbf{Z}[G]  )  $; thus in the
sequel we may view ${\rm R}(  \mathbf{Z}[G]   )$ as a
subgroup of $A^{\text{\rm s}} (  \mathbf{Z}[G]  )  $.
\end{theorem}

\bigskip Viewing ${\rm A}^{\text{s}} (  \mathbf{Z} [G]   )  $
as a subgroup of ${\rm H}^{\text{s}} (  \mathbf{Z} [G]   )$, we obtain the natural decomposition
\begin{equation}
\text{H}^{\text{s}} (  \mathbf{Z} [G]   )  =\text{A}
^{\text{s}}(  \mathbf{Z} [G]   )  \times\text{S}
_{\infty}(  \mathbf{Z} [G]   )  .
\end{equation}
\medskip

\subsection{\bigskip Formation of Euler characteristics.}

\subsubsection{Definitions.}

\textbf{Symplectic hermitian case}. From now on we fix a set of symplectic
$\mathbf{C}\left[  G\right]  $-representations $W_{m}$ whose characters
$\theta_{m}$ form a $\mathbf{Z}$-basis of $R_{G}^{\text{s}}$. There is of
course a natural $\mathbf{Z}$-basis for $R_{G}^{\text{s}}$ given by the
irreducible symplectic characters and the sums of the irreducible
non-symplectic characters and their contragredients; in the sequel we shall
assume our basis to be of this form. We then fix a non-degenerate
$G$-invariant alternating form $\kappa_{m}$ on $W_{m}$ and we let $ \{
w_{mn} \}  $ denote a hyperbolic basis of $W_{m}$ with respect to
$\kappa_{m}.$

Suppose now that we are given a perfect $\mathbf{Z} [G]   $-complex
$P^{\bullet}$ with $G$-invariant non-degenerate real-valued symmetric forms
$\sigma^{\text{ev}}$ resp. $\sigma^{\text{odd}}$ on ${\rm H}^{\text{ev}} (
P_{\mathbf{Q}}^{\bullet} )  $ resp. ${\rm H}^{\text{odd}} (  P_{\mathbf{Q}
}^{\bullet} )  $. For each prime $p$ of $\mathbf{Z}$ let $ \{
a_{p}^{ij} \}  _{j}$ denote a $\mathbf{Z}_{p} [G]  $ basis
for $\mathbf{Z}_{p}\otimes_{\mathbf{Z}}P^{i}$; similarly we choose a
$\mathbf{Q} [G]  $ basis $ \{  a_{0}^{ij} \}  _{j}$ for
$P_{\mathbf{Q}}^{i}=\mathbf{Q}\otimes_{\mathbf{Z}}P^{i};$ then\ for each prime
$p$ let $\lambda_{p}^{i}$ be the element of $GL (  \mathbf{Q}_{p} [G]   )  $ such that $\lambda_{p}^{i}a_{p}^{ij}=a_{0}^{ij}.$

The following lemma is now clear:

\begin{lemma}
For any free $\mathbf{C} [G]   $-module $U$ with basis $\left\{
u_{i}\right\}$, the map
\[
r_{G}:U\otimes_{\mathbf{C}}W_{m}\rightarrow\left(  U\otimes_{\mathbf{C}}%
W_{m}\right)  ^{G}\text{ }%
\]
defined by $r_{G}\left(  u\otimes w\right)  =%
{\textstyle\sum_{g}}
gu\otimes gw$ is a surjection and $\left\{  r_{G}\left(  u_{i}\otimes
w_{mn}\right)  \right\}  _{i,n}$ is a basis of $\left(  U\otimes_{\mathbf{C}%
}W_{m}\right)  ^{G}$.\hfill$\square$
\end{lemma}

\medskip For each pair $i,m$ we put
\[
b_{jn}^{im}=r_{G} (  a_{0}^{ij}\otimes w_{mn} )  .
\]
Then by the above lemma, since $P_{\mathbf{Q}}^{i}$ is 
$\mathbf{Q} [G]   $-free, $ \{  b_{jn}^{im} \}  _{jn}$ is a $\mathbf{C}
$-basis of $ (  P^{i}\otimes_{\mathbf{Q}}W_{m} )  ^{G}$. As in (2.3)
we shall write $\xi_{m}$ for the canonical isomorphism
\[
\det (   (  P_{\mathbf{Q}}^{\bullet}\otimes_{\mathbf{Q}}W_{m} )
^{G} )  \cong\det (   (  \text{H}^{\bullet} (  P_{\mathbf{Q}
}^{\bullet} )  \otimes_{\mathbf{Q}}W_{m} )  ^{G} )  .
\]
Since all the terms in the complexes $P_{\mathbf{Q}}^{\bullet}$ and
H$^{\bullet} (  P_{\mathbf{Q}}^{\bullet} )  $ are $\mathbf{Q} [G]  $-modules, 
because the representation $W_{m}$ is symplectic, it
follows that all the terms in the complexes $ (  P_{\mathbf{Q}}^{\bullet
}\otimes_{\mathbf{Q}}W_{m} )  ^{G}$ and $ (  \text{H}^{\bullet
} (  P_{\mathbf{Q}}^{\bullet} )  \otimes_{\mathbf{Q}}W_{m} )
^{G}$ are even dimensional. Indeed, let $M$ be an $\mathbf{R} [G] 
$-module: if $W_{m}$ is an irreducible symplectic $\mathbf{R} [G] 
$-module, then $W_{m}$ has real Schur index 2, and so $\dim (  M\otimes_{\R}
W_{m} )  ^{G}$ is even; on the other hand if $W_{m}$ can be written as
$V+V^{\ast}$ for some $\mathbf{C} [G]  $-module $V$ with $V^{\ast}
$ denoting the contragredient of $V$, then $\dim (  M\otimes_{\R} V )
^{G}=\dim (  M\otimes_{\R} V^{\ast} )  ^{G}$ and so again $\dim (
M\otimes_{\R} W_{m} )  ^{G}$ is even. In particular we note that this means
that by Remark 7(b) we may treat the natural isomorphism $\upsilon
_{H_{m}^{\bullet}}$
\begin{equation}
\det (   (  \text{H}^{\bullet} (  P_{\mathbf{Q}}^{\bullet} )
\otimes_{\mathbf{Q}}W_{m} )  ^{G} )  \rightarrow\det (   (
\text{H}^{\text{ev}} (  P_{\mathbf{Q}}^{\bullet} )  \otimes
_{\mathbf{Q}}W_{m} )  ^{G} )  \otimes\det (  (
\text{H}^{\text{odd}} (  P_{\mathbf{Q}}^{\bullet} )  \otimes
_{\mathbf{Q}}W_{m} )  ^{G} )  ^{-1}
\end{equation}
as an identification with no sign changes; similarly, by Remark 7(a), given
another perfect $\mathbf{Z} [G]   $-complex $Q^{\bullet}$, we can
and shall also identify
\begin{equation}
\det (   (   (  P_{\mathbf{Q}}^{\bullet}\oplus Q_{\mathbf{Q}%
}^{\bullet} )  \otimes_{\mathbf{Q}}W_{m} )  ^{G} )
=\det (   (  P_{\mathbf{Q}}^{\bullet}\otimes_{\mathbf{Q}}W_{m} )
^{G} )  \otimes\det (   (  Q_{\mathbf{Q}}^{\bullet}
\otimes_{\mathbf{Q}}W_{m} )  ^{G} )  .
\end{equation}

\begin{definition}
We define $\chi_{\text{\rm H}}^{\text{\rm s}} (  P^{\bullet},\sigma)  \in
{\rm H}^{\text{\rm s}} (  \mathbf{Z} [G]   )  $ to be the class
represented by the character map which sends the character $\theta_{m}$ to
\begin{equation}
{\textstyle\prod_{p<\infty}}
\text{\rm Det} (  \lambda_{p}^{i} )   (  \theta_{m} )  ^{ (
-1 )  ^{i}}\times\text{\rm Pf}_{ (  \sigma\otimes\kappa_{m} )  ^{G}
} \left(  \xi_{m} (  \otimes_{i} (  \wedge_{jn}b_{jn}^{im} )
^{ (  -1 )  ^{i}} )   \right)
\end{equation}
where the terms on the right are taken in lexicographic order. Thus in
particular\ for fixed $i,m$ writing $b_{jn}$ for $b_{jn}^{im}$, then
$\wedge_{jn}b_{jn}$ is taken to mean the exterior product
\[
b_{11}\wedge b_{12}\cdots\wedge b_{1n}\wedge b_{21}\cdots\wedge b_{nm}.
\]
\end{definition}

We now wish to show that this class is independent of all choices.

It is clear that if we change basis from the given $\mathbf{Z}_{p} [G]   $-basis for 
 $\mathbf{Z}_{p}\otimes P^{i}$, $\{
a_{p}^{ij} \}  _{j}$, then we only change the representing character
function by an element in 
${\rm Det}^{\text{s}} (  \mathbf{Z}_{p} [G]  ^{\times} )  \times 1$. Similarly, if we change the given
$\mathbf{Q} [G]   $-basis for $\mathbf{Q}\otimes P^{i}$, $\{
a_{0}^{ij} \}  _{j}$, then we only change the representing character
function by an element in $\operatorname{Im}(\Delta^{s})$.\ 

Next we consider the possible dependence on the alternating forms $\kappa_{m}$
\ and the chosen hyperbolic basis $ \{  w_{mn} \}  $. Let $\eta_{m}$
be a further non-degenerate $G$-invariant alternating form on $W_{m}$, let
$ \{  w_{mn}^{\prime} \}  $ denote a hyperbolic basis of $W_{m}$ with
respect to $\eta_{m}$ and put
\[
b_{jn}^{\prime im}=r_{G} (  a_{0}^{ij}\otimes w_{mn}^{\prime} )\ .
\]
In order to show that the value in (8) does not change, we must show that
\[
\text{Pf}_{ (  \sigma\otimes\kappa_{m} )  ^{G}}\left(  \xi_{m} (
\otimes_{i} (  \wedge_{jn}b_{jn}^{im} )  ^{\left(  -1\right)  ^{i}
})  \right)
\]%
\[
=\text{Pf}_{\left(  \sigma\otimes\eta_{m}\right)  ^{G}}\left(  \xi_{m}(
\otimes_{i}(  \wedge_{jn}b_{jn}^{\prime im})  ^{\left(  -1\right)
^{i}})  \right)  .
\]
This will follow at once from Proposition 23 below and the fact that, by the
corollary to Proposition 4.1 and 4.2 on page 35-37 of [F1], the right hand
expression (12) in Proposition 23 is independent of the particular alternating
form $\kappa_{m}$ used.\hfill $\square$

\bigskip For future reference we record the following two results which follow
at once from the above definition and Lemma 3.

\begin{lemma}
(a) Suppose the complex $P^{\bullet}$ is acyclic, and let $0$ denote the
trivial form on the trivial vector space ${\rm H}^{\bullet} (  P^{\bullet
} )  = \{  0 \}$; then $\chi_{\text{\rm H}}^{\text{\rm s}} (
P^{\bullet},0 )  $ is the trivial class.

(b) Given two perfect $\mathbf{Z} [G]  $-complexes $P_{i}%
^{\bullet}$ for  $i=1,2$ with non-degenerate $G$-invariant forms $\sigma
_{i}^{\text{ev}},\sigma_{i}^{\text{odd}}$ on ${\rm H}^{\text{\rm ev}} (
P_{i,\mathbf{Q}}^{\bullet} )$, ${\rm H}^{\text{\rm odd}} (  P_{i,\mathbf{Q}%
}^{\bullet} )  $; we view $\sigma_{1}^{\text{\rm ev}}\oplus\sigma
_{2}^{\text{\rm ev}}$ as a form on the even part of the cohomology of
$P_{1}^{\bullet}\oplus P_{2}^{\bullet}$ via the identification ${\rm H}^{\bullet
} (  P_{1}^{\bullet}\oplus P_{2}^{\bullet} )  ={\rm H}^{\bullet} (
P_{1}^{\bullet} )  \oplus{\rm H}^{\bullet} (  P_{2}^{\bullet} )  $;
we then have the equality of classes
\[
\chi_{\text{\rm H}}^{\text{\rm s}}\left(  P_{1}^{\bullet}\oplus P_{2}^{\bullet}%
,\sigma_{1}\oplus\sigma_{2}\right)  =\chi_{\text{\rm H}}^{\text{\rm s}}\left(
P_{1}^{\bullet},\sigma_{1}\right)  \chi_{\text{\rm H}}^{\text{\rm s}}\left(
P_{2}^{\bullet},\sigma_{2}\right)  .
\]
\end{lemma}

\textbf{Arakelov case}. Here we briefly recall the construction of the
Arakelov Euler characteristic given in [CPT2]. Let $\left\{  V_{r}\right\}  $
denote the distinct simple two sided ideals of the complex group algebra
$\mathbf{C} [G]  ,$ and let $\nu_{r}$ denote the hermitian form on
$V_{r}$ given by the restriction of the standard non-degenerate $G$-invariant
hermitian form $ \nu_{\mathbf{C}}:\mathbf{C} [  G ]
\times\mathbf{C} [  G ]  \rightarrow\mathbf{C}$%
\[
\nu_{\mathbf{C}} (
{\textstyle\sum_{g\in G}}
l_{g}g,%
{\textstyle\sum_{h\in G}}
m_{h}h )  =\left|  G\right|
{\textstyle\sum_{g\in G}}
l_{g}\overline{m_{g}}%
\]
and we let $\left\{  v_{rs}\right\}  $ denote an orthonormal basis of $V_{r}$
with respect to $\nu_{r}$.

We next suppose that we are given a perfect $\mathbf{Z} [G] 
$-complex $P^{\bullet}$\ with \linebreak metrics $h=\left\{  h_{r}\right\}  $
on the equivariant determinant of cohomology i.e. each $h_{r}$ is a metric on
the line $\det(  (  H^{\bullet}(  P^{\bullet})
\otimes_{\mathbf{Q}}V_{r})  ^{G})  $. We again let $\{
a_{p}^{ij}\}  _{j}$ denote a $\mathbf{Z}_{p} [G]  $-basis
for $\mathbf{Z}_{p}\otimes P^{i}$\ and let $\{  a_{0}^{ij}\}
_{j} $ denote a $\mathbf{Q} [G]   $-basis for $\mathbf{Q}\otimes
P^{i}$;  as previously, we let $\lambda_{p}^{i}$ be the element of $GL (
\mathbf{Q}_{p} [  G ]  )  $ such that $\lambda_{p}^{i}
a_{p}^{ij}=a_{0}^{ij}$. Then for each pair $i,r$ we put
\[
c_{js}^{ir}=r_{G} (  a_{0}^{ij}\otimes v_{rs} )
\]
and again by Lemma 15 we know that $ \{  c_{js}^{ir} \}  $ is a
$\mathbf{C}$-basis of $ (  P^{i}\otimes_{\mathbf{Q}}V_{r} )  ^{G}$.

\begin{definition}
The equivariant Arakelov class $\chi_{\text{\rm A}} (  P^{\bullet},h )
\in{\rm A}(  \mathbf{Z} [  G]  )  $  is defined to be the
class represented by the following homomorphism on characters: if $V_{r}$ has
character $\chi_{r}$, then the complex conjugate $\overline{\chi}_{r}$ is sent
to the value
\[%
{\textstyle\prod_{p<\infty}}
\text{\rm Det} (  \lambda_{p}^{i} )   (  \overline{\chi}_{r} )
^{ (  -1 )  ^{i}}\times h_{r} (  \xi_{r} (\otimes_{i} (  \wedge_{js}c_{js}^{ir} )  ^{ (  -1 )  ^{i}} )
 )  ^{1/\chi_{r} (  1 )  }.
\]
From 3.3 in [CPT2] we know that the class given by this character map is again independent of
choices. The symplectic Arakelov class $\chi_{\text{\rm A}}^{\text{\rm s}} (
P^{\bullet},h ) \in {\rm A}^{\text{\rm s}} (  \mathbf{Z}[  G ]
 )  $ is then given by restricting the above character map to symplectic characters.
\end{definition}

\bigskip

\subsubsection{The hermitian metrics associated to a symmetric bilinear form.}

With the notation of the previous subsection, we suppose we are given
non-degenerate $G$-invariant real-valued symmetric bilinear forms
$\sigma^{\text{ev}}$, $\sigma^{\text{odd}}$ on ${\rm H}^{\text{ev}} (
P_{\mathbf{Q}}^{\bullet} )$, ${\rm H}^{\text{odd}} (
P_{\mathbf{Q}}^{\bullet} )  $. We now briefly recall how this data
naturally determines a system of metrics on the equivariant determinant of
cohomology of $P^{\bullet}$. We observe that $ (  \sigma^{\text{ev}%
}\otimes\nu_{r} )  ^{G}$ is a non-degenerate hermitian form on the vector space $ (
\text{H}^{\text{ev}} (  P_{\mathbf{Q}}^{\bullet} )  \otimes
_{\mathbf{Q}}V_{r} )  ^{G}$; the determinant of this form affords a
hermitian form $\det (   (  \sigma^{\text{ev}}\otimes\nu_{r} )
^{G} )  $ on the complex line $\det (   (  \text{H}^{\text{ev}%
} (  P_{\mathbf{Q}}^{\bullet} )  \otimes_{\mathbf{Q}}V_{r} )
^{G} )  $ which may be either positive or negative definite; multiplying
by $-1$ if the form is negative definite, in all cases we then obtain a
positive definite form which we denote by $\left|  \det (   (
\sigma^{\text{ev}}\otimes\nu_{r} )  ^{G} )  \right|  $. The positive
definite form $\left|  \det (   (  \sigma^{\text{odd}}\otimes\nu
_{r} )  ^{G} )  \right|  $ is defined similarly. We then write
$h_{r}$ for the metric on the complex line $\det (   (  \text{H}^{\bullet} (
P_{\mathbf{Q}}^{\bullet} )  \otimes_{\mathbf{Q}}V_{r} )
^{G} )  $ corresponding via $\upsilon_{\text{H}^{\bullet} (
P^{\bullet} )  } $ to the positive definite form
\[
\left|  \det (   (  \sigma^{\text{ev}}\otimes\nu_{r} )
^{G} )  \right|  \otimes\left|  \det ( (  \sigma^{\text{odd}
}\otimes\nu_{r} )  ^{G} )  \right|  ^{-1}.
\]
Then $h= \{  h_{r} \}  $ is then the required system of metrics on
the equivariant determinant of cohomology of $P^{\bullet}$.

\subsubsection{\bigskip Independence under quasi-isomorphism.}

We first recall the following result from Theorem 3.9 in [CPT2]: suppose
$P_{1}^{\bullet}$, resp. $P_{2}^{\bullet}$ is a perfect $\mathbf{Z} [
G ]  $-complex which supports metrics $h^{1}= \{  h_{r}^{1} \}
_{r}$ resp. $h^{2}= \{  h_{r}^{2} \}  _{r}$ on its equivariant
determinant of cohomology. Suppose further that there is a quasi-isomorphism
$\phi:P_{1}^{\bullet}\dashrightarrow P_{2}^{\bullet}$ in the derived category
of bounded complexes of finitely generated $\mathbf{Z} [  G ]
$-modules, which has the property that $\phi^{\ast}h^{1}=h^{2}$. Then we
know that the formation of Arakelov classes is natural with respect to
quasi-isomorphisms in the sense that
\[
\chi_{\text{A}} (  P_{1}^{\bullet},h^{1} )  =\chi_{\text{A}} (
P_{2}^{\bullet},h^{2} )  .
\]

We now establish the corresponding result for hermitian classes. Prior to
stating the result, we first need some more notation: for each $m,$ for
brevity we let $\det (  P_{i}^{\bullet} )  _{m}$ denote $\det (
 (  P_{i}^{\bullet}\otimes_{\mathbf{Q}}W_{m} )  ^{G} )$; we
write $\det(\phi)_{m}$ for the isomorphism
\[
\det(  P_{1}^{\bullet})  _{m}\cong\det(  P_{2}^{\bullet
})  _{m}
\]
induced by $\phi$, and we let $\xi_{m}^{i}$ denote the canonical isomorphism
\[
\det(  P_{i}^{\bullet})  _{m}\cong\det(  \text{H}^{\bullet
}(  P_{i}^{\bullet})  )  _{m}.
\]
Then we have the following result:

\begin{theorem}
Suppose $P_{1}^{\bullet}$, $P_{2}^{\bullet}$ are perfect 
$\mathbf{Z} [G]  $-complexes which support non-degenerate $G$-invariant
real-valued symmetric forms $\sigma_{1}^{\text{\rm ev}}$, $\sigma_{1}^{\text{\rm odd}
}$, $\sigma_{2}^{\text{\rm ev}}$, $\sigma_{2}^{\text{\rm odd}}$, on their rational
cohomology, and suppose, as previously, that there is a quasi-isomorphism
$\phi:P_{1}^{\bullet}\dashrightarrow P_{2}^{\bullet}$ in the derived category
of bounded complexes of finitely generated $\mathbf{Z}[G]
$-modules, which has the property that
\begin{equation}
\text{\rm Pf}_{(  \sigma_{1,\mathbf{C}}\otimes\kappa_{m})  ^{G}}%
\circ\xi_{m}^{1}=\text{\rm Pf}_{ (  \sigma_{2,\mathbf{C}}\otimes\kappa
_{m} )  ^{G}}\circ\xi_{m}^{2}\circ\det\left(  \phi_{m}\right)  .
\end{equation}
Then there is an equality of hermitian classes:
\[
\chi_{\text{\rm H}}^{\text{\rm s}} (  P_{1}^{\bullet},\sigma )
=\chi_{\text{\rm H}}^{\text{\rm s}} (  P_{2}^{\bullet},\sigma_{2} )  .
\]
\end{theorem}

\noindent \textit{Proof.} We may write the quasi-isomorphism $\phi:P_{1}^{\bullet
}\dashrightarrow P_{2}^{\bullet}$ as
\[
P_{1}^{\bullet}\overset{\psi_{1}}{\longleftarrow}Q^{\bullet}\overset{\psi_{2}%
}{\rightarrow}P_{2}^{\bullet}%
\]
with $Q^{\bullet}$ a bounded complex of finitely generated 
$\mathbf{Z} [G]   $-modules and with the $\psi_{i}$ quasi-isomorphic chain maps. As a
first step, we observe that we can in fact choose $Q^{\bullet}$ to be a
perfect $\mathbf{Z} [G]  $-complex: this follows from a standard
argument and the reader is referred to Lemma 5.1 in [CPT2] for the details.

Next we observe that, by adding on a sufficiently large acyclic free complex
$A_{1}^{\bullet}$, we can ensure that there is chain map $\alpha$ such that
\[
\psi_{1}^{\prime}=\psi_{1}\oplus\alpha:Q^{\bullet}\oplus A_{1}^{\bullet
}\rightarrow P_{1}^{\bullet}%
\]
is surjective. Then by VI.8.17 in [M] we can split $\psi_{1}^{\prime}$ by a
quasi-isomorphism $\beta_{1}:P_{1}^{\bullet}\rightarrow Q^{\bullet}\oplus
A_{1}^{\bullet}$, and so we obtain a direct decomposition
\[
Q^{\bullet}\oplus A_{1}^{\bullet}=\operatorname{Im}\left(  \beta_{1}\right)
\oplus\ker\left(  \psi_{1}^{\prime}\right)  .
\]
Let $\tau_{1}^{\text{ev}}$, $\tau_{1}^{\text{odd}}$ denote the transport of
$\sigma_{1}^{\text{ev}}$, $\sigma_{1}^{\text{odd}}$ to ${\rm H}^{\text{ev}}\left(
\operatorname{Im}\left(  \beta_{1}\right)  \right)$, ${\rm H}^{\text{odd}}\left(
\operatorname{Im}\left(  \beta_{1}\right)  \right)  $ and of course we endow
the (zero) cohomology of the acyclic complex $\ker\left(  \psi_{1}^{\prime
}\right)  $ with the zero form. We then repeat these constructions for
$P_{2}^{\bullet}$. Putting this together and using Lemma 17 we get
\begin{align*}
\chi_{\text{H}}^{\text{s}}\left(  A_{2}^{\bullet}\oplus Q^{\bullet}\oplus
A_{1}^{\bullet},0\oplus\tau_{1}\oplus0\right)   & =\chi_{\text{H}}^{\text{s}%
}\left(  A_{2}^{\bullet},0\right)  \chi_{\text{H}}^{\text{s}}\left(
Q^{\bullet}\oplus A_{1}^{\bullet},\tau_{1}\oplus0\right) \\
& =\chi_{\text{H}}^{\text{s}}\left(  Q^{\bullet}\oplus A_{1}^{\bullet}%
,\tau_{1}\oplus0\right)
\end{align*}
while
\begin{align*}
\chi_{\text{H}}^{\text{s}}\left(  Q^{\bullet}\oplus A_{1}^{\bullet},\tau
_{1}\oplus0\right)   & =\chi_{\text{H}}^{\text{s}}\left(  \operatorname{Im}%
\left(  \beta_{1}\right)  \oplus\ker\left(  \psi_{1}\right)  ,\tau_{1}%
\oplus0\right) \\
& =\chi_{\text{H}}^{\text{s}}\left(  \operatorname{Im}\left(  \beta
_{1}\right)  ,\tau_{1}\right)  \chi_{\text{H}}^{\text{s}}\left(  \ker\left(
\psi_{1}\right)  ,0\right) \\
& =\chi_{\text{H}}^{\text{s}}\left(  P_{1}^{\bullet},\sigma_{1}\right)
\chi_{\text{H}}^{\text{s}}\left(  \ker\left(  \psi_{1}\right)  ,0\right)
=\chi_{\text{H}}^{\text{s}}\left(  P_{1}^{\bullet},\sigma_{1}\right)  .
\end{align*}
Thus we see that for $i=1,2$%
\[
\chi_{\text{H}}^{\text{s}}\left(  A_{2}^{\bullet}\oplus Q^{\bullet}\oplus
A_{1}^{\bullet},0\oplus\tau_{i}\oplus0\right)  =\chi_{\text{H}}^{\text{s}%
}\left(  P_{i}^{\bullet},\sigma_{i}\right)
\]
and so to prove the theorem it will suffice to show that
\[
\chi_{\text{H}}^{\text{s}}\left(  A_{2}^{\bullet}\oplus Q^{\bullet}\oplus
A_{1}^{\bullet},0\oplus\tau_{1}\oplus0\right)  =\chi_{\text{H}}^{\text{s}%
}\left(  A_{2}^{\bullet}\oplus Q^{\bullet}\oplus A_{1}^{\bullet},0\oplus
\tau_{2}\oplus0\right)  .
\]
(Note here that the $0\oplus\tau_{i}\oplus0$ refer to different direct sum
decompositions for $i=1,2$.)\ In order to show this it will suffice to show
that for each $m$
\begin{equation}
\text{Pf}_{\left(  \left(  0\oplus\tau_{1}\oplus0\right)  \otimes\kappa
_{m}\right)  ^{G}}=\text{Pf}_{\left(  \left(  0\oplus\tau_{2}\oplus0\right)
\otimes\kappa_{m}\right)  ^{G}}.
\end{equation}
To establish this equality we consider the commutative diagram
\[%
\begin{array}
[c]{ccc}%
\det\left(  Q^{\bullet}\right)  _{m} & \overset{\det\psi_{m}}{\longrightarrow} &
\det\left(  P_{1}^{\bullet}\right)  _{m}\\
\downarrow\xi_{Q,m} &  & \downarrow\xi_{m}^{1}\\
\det\left(  \text{H}^{\bullet}\left(  Q^{\bullet}\right)  \right)  _{m} &
\overset{\det\text{H}\left(  \psi\right)  _{m}}{\longrightarrow} & \det\left(
\text{H}^{\bullet}\left(  P_{1}^{\bullet}\right)  \right)  _{m}.
\end{array}
\]
Writing $\det\left(  \phi\right)  _{m}=\det\left(  \psi_{2}\right)  _{m}%
\det\left(  \psi_{1}\right)  _{m}^{-1}$, from (9) we get
\[
\text{Pf}_{\left(  \sigma_{1,\mathbf{C}}\otimes\kappa_{m}\right)  ^{G}}%
\circ\xi_{m}^{1}\circ\det\left(  \psi_{1}\right)  _{m}=\text{Pf}_{\left(
\sigma_{2,\mathbf{C}}\otimes\kappa_{m}\right)  ^{G}}\circ\xi_{m}^{2}\circ
\det\left(  \psi_{2}\right)  _{m}%
\]
and so from the commutative diagram we deduce that
\begin{equation}
\text{Pf}_{\left(  \sigma_{1,\mathbf{C}}\otimes\kappa_{m}\right)  ^{G}}%
\circ\det\left(  \text{H}\left(  \psi_{1}\right)  \right)  _{m}=\text{Pf}%
_{\left(  \sigma_{2,\mathbf{C}}\otimes\kappa_{m}\right)  ^{G}}\circ\det\left(
\text{H}\left(  \psi_{2}\right)  \right)  _{m}.
\end{equation}
The equality (10) then follows since, by construction, under the isomorphism
\[
\text{H}^{\bullet}\left(  P_{1}^{\bullet}\right)  \overset{\text{H}\left(
\psi_{1}\right)  }{\leftarrow}\text{H}^{\bullet}\left(  Q^{\bullet}\right)
=\text{H}^{\bullet}\left(  A_{2}^{\bullet}\oplus Q^{\bullet}\oplus
A_{1}^{\bullet}\right)  =\text{H}^{\bullet}\left(  A_{2}^{\bullet}%
\oplus\operatorname{Im}\left(  \beta_{1}\right)  \oplus\ker\left(  \psi
_{1}\right)  \oplus A_{1}^{\bullet}\right)
\]
$\sigma_{1}$ transports to $0\oplus\tau_{1}\oplus0;$ thus\ identifying
H$^{\bullet}\left(  Q^{\bullet}\right)  $ with H$^{\bullet}\left(
A_{2}^{\bullet}\oplus Q^{\bullet}\oplus A_{1}^{\bullet}\right)  $ we read (11)
as (10), as required.\hfill $\square$

\subsubsection{Evaluation of Pfaffians.}

The results we require in this subsection come from Appendix C in [CPT1].
Throughout this sub-section we again  assume that all real and complex vector
spaces are finite dimensional. We begin by considering the hermitian form
associated to a $G$-invariant symmetric form.

\bigskip

Suppose now that $U$ is a free $\mathbf{R} [  G ]  $-module with
basis $ \{  u_{i}:1\leq i\leq q \}  $ which again supports a
real-valued non-degenerate $G$-invariant symmetric form $\sigma$. We then
write
\[
\widetilde{\sigma}:U\times U\rightarrow\mathbf{R} [G] 
\]
for the associated group ring valued hermitian form (ref. page
25 in [F1]); thus for $u,u^{\prime}\in U$%
\[
\widetilde{\sigma}\left(  u,u^{\prime}\right)  =%
{\textstyle\sum_{g\in G}}
\sigma\left(  gu,u^{\prime}\right)  g^{-1}.
\]
Let $r\mapsto$ $\overline{r}$ denote the $\mathbf{R}$-linear involution on
$\mathbf{R} [  G ]  $ induced by group inversion. Note that
\[
\overline{\widetilde{\sigma}\left(  u_{i},u_{j}\right)  }=\widetilde{\sigma
}\left(  u_{j},u_{i}\right)  .
\]

\begin{proposition}
Let $V$ be a $\mathbf{C} [  G ]  $-module which supports a
$G$-invariant non-degenerate form $h$; then $q\times q$ matrices with entries
in $\mathbf{C} [  G ]  $ act on the direct sum of $q $ copies of V,
and the matrix $T=\left(  \widetilde{\sigma}\left(  u_{i},u_{j}\right)
\right)  _{i,j}$ is self adjoint with respect to $h^{\left(  q\right)  }$, the
orthogonal direct sum of $h$ on $q$-copies of $V$.
\end{proposition}

\noindent{\sl Proof.} Indeed, writing $v_{ (  i )  }$ for the vector in
$\oplus_{i=1}^{q}V$ which is $v$ in the $i$-th position and zero elsewhere, we
have
\[
h^{ (  q )  } (  v_{ (  i )  },Tv_{ (  j )
}^{\prime} )  =h^{ (  q )  } (  v_{ (  i )  },
 {\textstyle\sum_{k}}
 \widetilde{\sigma} (  u_{k},u_{j} )  v_{ (  k )  }^{\prime
} )=
\]
\[
=h^{(  q )  } (  v_{\left(  i\right)  },\widetilde{\sigma
}\left(  u_{i},u_{j}\right)  v_{\left(  i\right)  }^{\prime} )
=h^{ (  q )  } (  \overline{\widetilde{\sigma}\left(  u_{i}%
,u_{j}\right)  }v_{\left(  i\right)  },v_{\left(  i\right)  }^{\prime} )=
\]
\[
=h^{ (  q )  } (  \widetilde{\sigma}\left(  u_{j},u_{i}\right)
v_{\left(  i\right)  },v_{\left(  i\right)  }^{\prime} )  =h (
\widetilde{\sigma}\left(  u_{j},u_{i}\right)  v,v^{\prime} )  ;
\]
while similarly
\[
h^{(  q)  }(  Tv_{\left(  i\right)  },v_{\left(  j\right)
}^{\prime})  =h^{(  q)  }(
{\textstyle\sum_{l}}
\widetilde{\sigma}(  u_{l},u_{i} )  v_{\left(  l\right)
},v_{\left(  j\right)  }^{\prime})=
\]%
\[
=h^{(  q)  }(  \widetilde{\sigma}\left(  u_{j},u_{i}\right)
v_{\left(  j\right)  },v_{\left(  j\right)  }^{\prime})  =h(
\widetilde{\sigma}(  u_{j},u_{i})  v,v^{\prime})
.\;\;\;\square
\]

\begin{lemma}
 Suppose, as previously, that $W_{m}$ is a symplectic $\mathbf{C}%
 [  G ]  $-module with non-degenerate $ G$-invariant alternating
form $\kappa_{m}$. The map $r_{G}:W_{m}\rightarrow (  \mathbf{R} [
G ]  \otimes_{\mathbf{R}}W_{m} )  ^{G}$ (see Lemma 15) given by
\[
r_{G}(  w )  = 
{\textstyle\sum_{g\in G}}
g\otimes gw
\]
has the property that $|G|^{-1} r_{G}$ is a $G$-isometry
when $ (  \mathbf{R} [  G ]  \otimes_{\mathbf{R}}W_{m} )
^{G} $ is endowed with the form $ (  \nu_{\mathbf{R}}\otimes\kappa
_{m} )  ^{G}$. (Recall that $\nu_{\mathbf{C}}$ was defined after Lemma
17; $\nu_{\mathbf{R}}$ is the restriction of $\nu_{\mathbf{C}}$ to the real
group algebra $\mathbf{R} [  G ]$.)
\end{lemma}

\noindent {\sl Proof.} First note that $|G|^{-1}r_{G}$ has inverse $%
{\textstyle\sum}
a_{g}g\otimes w\mapsto%
{\textstyle\sum}
a_{g}g^{-1}w$; we then observe that
\[
\left(  \nu\otimes\kappa_{m}\right)  ^{G} (  |G|^{-1}
r_{G}\left(  w\right)  ,|G|^{-1}r_{G} (  w^{\prime} )   )  =|G|^{-2}(\nu\otimes\kappa
_{m}) (
{\textstyle\sum_{g}}
g\otimes gw,%
{\textstyle\sum_{h}}
h\otimes hw )
\]%
\[
=|G|^{-1}
{\textstyle\sum_{g}}
\kappa_{m}\left(  gw,gw^{\prime}\right)  =\kappa_{m}\left(  w,w^{\prime
}\right)  .\;\;\;\square
\]

\begin{remark}
{\rm Note for future reference that ${\rm Im} (  r_{G} )  $ is a left
$\mathbf{R} [  G ]  $-module by transport of structure; to be more
precise, for $h\in G$ we define $h\cdot r_{G}\left(  w\right)  $ to be
\[
r_{G}\left(  hw\right)  =%
{\textstyle\sum_{g\in G}}
g\otimes ghw=r_{G}\left(  w\right)  \left(  h^{-1}\otimes1\right)  .
\]}
\end{remark}

For a symplectic $\mathbf{C}\left[  G\right]  $-module $W_{m}$ we write
$T_{W_{m}}^{\left(  q\right)  }$ for the representation of $G$ afforded by the
direct sum of $q$-copies of $W_{m}$; thus $T_{W_{m}}^{\left(  q\right)  }$ 
provides an action of $q\times q$ matrices with entries in $\mathbf{C}
[  G]  $ on the direct sum of $q$ copies of $W_{m}$.

We shall henceforth identify $U$ with $\oplus_{i=1}^{q}\mathbf{R}[
G] $ and let $\{u_{i}\}$ be the $\mathbf{R}[  G]  $-basis
of $U$ given by the canonical basis $\{1_{i}\}$ of $\oplus_{i=1}
^{q}\mathbf{R}[  G]  $. We define a form $\nu$ on $U$ by the rule
\[
\nu\left(  \lambda u_{i},\mu u_{j}\right)  =\delta_{ij}\nu_{\mathbf{R}}\left(
\lambda,\mu\right)  .
\]
For $w,w^{\prime}\in W_{m}$ and given $i,j$ we consider $x=\left|  G\right|
^{-1}r_{G} (  w_{\left(  i\right)  } )$, $y=\left|  G\right|
^{-1}r_{G} (  w_{\left(  j\right)  }^{\prime} )  $; note that by
Lemma 21, if $w$ ranges through a hyperbolic basis of $W_{m}$ with respect to
$\kappa_{m}$ and if $ i$ ranges from $1$ to $q$, then $x$ ranges through a
hyperbolic basis of $ (  U\otimes W_{m} )  ^{G}$ with respect to
$ (  \nu\otimes\kappa_{m} )^{G}$.

The following result will be fundamental in enabling us to calculate with Pfaffians.

\begin{proposition}
With the above notation  $\left|  G\right|  ^{-1}r_{G}$ defines an isometry
\[
\oplus_{1}^{q}W_{m}\cong\left(  U\otimes W_{m}\right)  ^{G}%
\]
where $\oplus_{1}^{q}W_{m}$ is endowed with the form $\kappa_{m}^{\left(
q\right)  }$ and $\left(  U\otimes W_{m}\right)  ^{G}$ is endowed with the
form $\left(  \nu\otimes\kappa_{m}\right)  ^{G}$. By Proposition 20,
$T_{W_{m}}^{ (  q )  } (  \widetilde{\sigma} (  u_{i}
,u_{j} )  )  $ is a self-adjoint automorphism with respect to
$\kappa_{m}^{ (  q )  }$. There is an equality of Pfaffians
\begin{equation}
\text{\rm Pf}_{\left(  \sigma\otimes\kappa_{m}\right)  ^{G}}\left(  \wedge
_{in}\left|  G\right|  ^{-1}r_{G}\left(  u_{i}\otimes w_{mn}\right)  \right)
=\mathbf{pf}_{\kappa_{m}^{\left(  q\right)  }}\left(  \left|  G\right|
^{-1}T_{W_{m}}^{\left(  q\right)  }\left(  \widetilde{\sigma}\left(
u_{i},u_{j}\right)  \right)  \right)
\end{equation}
where, as previously, $ \{  w_{mn} \}  $ denotes a hyperbolic basis
of $W_{m}$ with respect to $\kappa_{m}$ and the wedge product is of course
taken in lexicographic order.
\end{proposition}

\noindent \textit{Proof}. To prove the result, we claim that it will suffice to
establish the equality
\begin{equation}
 (  \sigma\otimes\kappa_{m} )   (  \left|  G\right|  x,\left|
G\right|  y )  =\kappa_{m}^{ (  q )  } (  w_{ (
i )  },\left|  G\right|  Tw_{ (  j )  }^{\prime} )  ,
\end{equation}
where we put $T=T_{W_{m}}^{(q)  } (  \widetilde{\sigma
} (  u_{i},u_{j} )   )  $. Indeed, assuming (13), we see that
\[
\left(  \sigma\otimes\kappa_{m}\right)  \left(  x,y\right)  =\kappa
_{m}^{(q) } (  w_{ (  i )  },\left|  G\right|
^{-1}Tw_{ (  j )  }^{\prime} )
\]
and so by Lemma 21
\begin{align*}
\kappa_{m}^{(q)} (  w_{ (i)  },\left|
G\right|  ^{-1}Tw_{(j)}^{\prime} )   & =(\nu\otimes
\kappa_{m}) (  \left|  G\right|  ^{-1}r_{G} (  w_{\left(  i\right)
} )  ,\left|  G\right|  ^{-1}r_{G} (  \left|  G\right|
^{-1}Tw_{ (  j) }^{\prime} )   ) \\
& =(\nu\otimes\kappa_{m}) (  x,\left|  G\right|  ^{-1}Ty )  .
\end{align*}
Hence by Proposition 6 it follows that
\begin{align*}
\text{Pf}_{\left(  \sigma\otimes\kappa_{m}\right)  ^{G}}\left(  \wedge
_{in}\left|  G\right|  ^{-1}r_{G}\left(  u_{i}\otimes w_{mn}\right)  \right)
& =\mathbf{pf}_{\kappa_{m}^{ (  q )  }} (  \left|  G\right|
^{-1}T )  \text{Pf}_{  (  \nu\otimes\kappa_{m}^{ (  q )
} )  ^{G}} \left(  \wedge_{in} (  \left|  G\right|  ^{-1}r_{G} (
w_{mn, (  i )}))\right) \\
& =\mathbf{pf}_{\kappa_{m}^{(q)}} (  |G|^{-1}T )
\end{align*}
with the last equality holding since the $ \{  \left|  G\right|^{-1}r_{G} (  w_{mn, (  i )  } )  \}  $ is a
hyperbolic basis of the vector space $ (U\otimes W_{m} )^{G}$ endowed with the
form $(\nu\otimes\kappa_{m}^{(q)})^{G}$.

To show (13), we consider the left-hand side, which can be written as
\[%
{\textstyle\sum_{g,h}}
\sigma\left(  gu_{i},hu_{j}\right)  \kappa_{m}\left(  gw,hw^{\prime}\right)
\]
which by $G$-invariance is
\[
\left|  G\right|
{\textstyle\sum_{h}}
\sigma\left(  u_{i},hu_{j}\right)  \kappa_{m}\left(  w,hw^{\prime}\right)
=\left|  G\right|  \kappa_{m}\left(  w,%
{\textstyle\sum_{h}}
\sigma\left(  u_{i},hu_{j}\right)  hw^{\prime}\right)  .
\]
The result then follows since
\[
\kappa_{m}^{\left(  q\right)  } (  w_{\left(  i\right)  },\left|
G\right|  Tw_{\left(  j\right)  }^{\prime} )  =\kappa_{m}^{\left(
q\right)  } (  w_{\left(  i\right)  },\left|  G\right|
{\textstyle\sum_{k}}
T_{kj}w_{\left(  k\right)  }^{\prime} )=
\]%
\[
=\kappa_{m}^{\left(  q\right)  } (  w_{\left(  i\right)  },\left|
G\right|  T_{ij}w_{\left(  i\right)  }^{\prime} )  =\kappa_{m} (
w,\left|  G\right|  T_{ij}w^{\prime} )=
\]%
\[
=\left|  G\right|  \kappa_{m} (  w,\widetilde{\sigma}\left(  u_{i}%
,u_{j}\right)  w^{\prime} )=\left|  G\right|  \kappa_{m}(w, {\textstyle\sum_{h}}
\sigma\left(  u_{i},hu_{j}\right)  hw^{\prime} )  .\;\;\;\square 
\]
\medskip

For future reference we also record the corresponding result for hermitian forms:

\bigskip

\begin{proposition}
Suppose that $V$ is a left ideal of $\mathbf{C}[G]$ endowed
with the non-degenerate  $G$-invariant hermitian form $\nu_{V}$ given by the
restriction of $\nu_{\mathbf{C}}$. As previously, 
we let $T_{V}^{(q)}$ denote the $G$-representation afforded by the direct sum of
$q$ copies of $V$. Let $h_{V}$ denote the hermitian form on the complex
vector space $ (  U\otimes_{\mathbf{R}}V )  ^{G}$ given by
restricting $\sigma\otimes\nu_{V}$. Let $ \{  v_{Vs} \}  $ be an
orthonormal basis of $V$\ with respect to $\nu_{V}.$ Then $T_{V}^{(q)}
( \widetilde{\sigma} (  u_{i},u_{j})  ) $ is
a self-adjoint automorphism with respect to $h_{V}$ and
\[
h_{V} (  \wedge_{is}r_{G}\left(  u_{i}\otimes v_{Vs}\right)   )
=\left|  \det (  \left|  G\right|  T_{V}^{\left(  q\right)  }\left(
\widetilde{\sigma}\left(  u_{i},u_{j}\right)  \right)   )  \right|
^{1/2}.
\]
\end{proposition}

\noindent{\sl Proof.} First we note that $\det\left(  \sigma\otimes\nu_{V}\right)  \left(
\wedge_{is}r_{G}\left(  u_{i}\otimes v_{Vs}\right)  \right)$ is equal to
the determinant
\[
\det\left(   (  \sigma\otimes\nu_{V}\right)  \left(
{\textstyle\sum_{g\in G}}
gu_{i}\otimes gv_{Vs},%
{\textstyle\sum_{h\in G}}
hu_{i^{\prime}}\otimes hv_{Vs^{\prime}} )  \right)
\]%
\[
=\det\left(
{\textstyle\sum_{g,h\in G}}
\sigma\left(  gu_{i},hu_{i^{\prime}}\right)  \nu_{V}\left(  gv_{Vs}%
,hv_{Vs^{\prime}}\right)  \right),
\]
which by setting $k=g^{-1}h$ and using the $G$-invariance of $\sigma$ and $\nu_{V}$ is
\[
=\det\left(  \left|  G\right|
{\textstyle\sum_{k\in G}}
\sigma\left(  u_{i},ku_{i^{\prime}}\right)  \nu_{V}\left(  v_{Vs}%
,kv_{Vs^{\prime}}\right)  \right)  .
\]
The matrix $\left(  \nu_{V}\left(  v_{Vs},kv_{Vs^{\prime}}\right)  \right)
_{s,s^{\prime}}$ is exactly the matrix of $T_{V}^{(q)}(
k)$ relative to the orthonormal basis $ \{  v_{Vs} \}  $ of
$V$ with respect to $\nu_{V}$. It therefore follows that the above determinant
is equal to
\[
\det (  \left|  G\right|  T_{V}^{(q)} (  \widetilde
{\sigma} (  u_{i},u_{i^{\prime}} )  )   )
\]
since $\widetilde{\sigma}\left(  u_{i},u_{i^{\prime}}\right)  =%
{\textstyle\sum_{k\in G}}
\sigma\left(  u_{i},ku_{i^{\prime}}\right)  k.\;\;\;\;\;\square$

\bigskip

\subsubsection{Comparison of hermitian and Arakelov classes.}

Recall that $\sigma^{\text{ev}}$ and $\sigma^{\text{odd}}$ denote
$G$-invariant real-valued symmetric forms on ${\rm H}^{\text{ev}}(  P^{\bullet
})  $ and ${\rm H}^{\text{odd}}(P^{\bullet})$ which, as per
3.2.2, induce metrics $h_{\sigma}= \{  h_{r} \}  $ on the equivariant
determinant of cohomology of $P^{\bullet}$. In this subsection we shall
compare the invariants $\chi_{\text{H}}^{\text{s}} (  P^{\bullet}
,\sigma )  $ and $\chi_{\text{A}}^{\text{s}} (  P^{\bullet
},h_{\sigma} )$. First we need the following algebraic result:

\begin{proposition}
Given an $\mathbf{R}[G]$-module $M$ and a
non-degenerate $G$-invariant symmetric form $\sigma$ on $M$, there exists a
$G$-decomposition $M=M^{+}\oplus M^{-}$ where $\sigma$ is positive definite on
$M^{+}$ and negative definite on $M^{-}$. This decomposition is not
necessarily unique, but the characters of the action of $G$ on $M^{+}$ and
$M^{-}$ are independent of choices.
\end{proposition}

\noindent{\sl Proof.} For full details see page 578 in [AS]; we briefly sketch a proof for
the reader's convenience. First we choose a $G$-invariant positive definite
symmetric form $\tau$ on $M$; there is then a unique automorphism $A$ of $M$
such that for all $x,y\in M$
\[
\sigma(  x,y)  =\tau(  x,Ay).
\]
As both $\sigma$ and $\tau$ are symmetric, $A$ is self adjoint with respect to
$\tau$; furthermore, since both $\sigma$ and $\tau$ are $G$-invariant, $A$
commutes with the action of $G$; thus the different eigenspaces of $A$ are
preserved by $G$; then, by considering the sums of eigenspaces for positive
and negative eigenvalues, we obtain the required decomposition $M=M^{+}\oplus
M^{-}$.

Clearly the above decomposition depends on the choice of $\tau$. To see that
the characters of $M^{+}$ and $M^{-}$ are independent of the choice of $\tau$,
we note that: the space of positive definite $G$-invariant forms on $M$
is connected; the maps $\tau\mapsto {\rm char}(M^{\pm})$ are
continuous; ${\rm char}(  M^{\pm})$ takes values in the discrete group
$R_{G}$.\hfill $\square$
\medskip

A particularly simple, but nonetheless useful, instance of the above
decomposition occurs when $ (  M,\sigma )  $ is hyperbolic. To state
this result we first need some notation. Recall that for an $\mathbf{R}[
G]$-module $V$ the hyperbolic space is ${\rm Hyp} (  V )  =V\oplus
V^{D}$ endowed with the form $h$ such that
\[
h\left(  v\oplus f,v^{\prime}\oplus f^{\prime}\right)  =f\left(  v^{\prime
}\right)  +f^{\prime}\left(  v\right)  \;\;\text{\ for }v,v^{\prime}\in
V,\;f,f^{\prime}\in V^{D}.
\]

\begin{lemma}
There are $\mathbf{R} [  G ]  $-isomorphisms ${\rm Hyp}(V)
^{+}\cong V\cong{\rm Hyp} (  V )  ^{-}.$
\end{lemma}

\noindent {\sl Proof.} First note that since $V$ is defined over $\mathbf{R,}$ we know that
$V\cong V^{D}$ as $\mathbf{R}[G ]  $-modules. As ${\rm Hyp} (
V\oplus W )  \cong{\rm Hyp}(V )  \oplus{\rm Hyp} (  W )
$, we see that it will suffice to prove the lemma when $V$ is irreducible over
$\mathbf{R}$. The result then follows immediately, since ${\rm Hyp}\left(  V\right)
^{+}$, ${\rm Hyp}(  V )  ^{-}$ are both $\mathbf{R} [  G ]
$-submodules of ${\rm Hyp} (  V )  $ and since hyperbolic spaces have zero
signature.\hfill $\square $
\medskip

In the final section we shall need the following result on hyperbolic
summands of quadratic modules:

\begin{lemma}
Let $K$ be an arbitrary field of characteristic zero, and let $\sigma$ be a
non-degenerate $K$-valued $G$-invariant symmetric form on a finite dimensional
$K [  G ]  $-module $V$. Suppose that $W$ is an isotropic $K [
G ]  $-submodule of $V$ and let $W^{\bot}$ denote the space of vectors
orthogonal to $W$. Then there is an orthogonal decomposition of $K [
G ]  $-modules
\[
V\cong {\rm Hyp} (  W )  \oplus\frac{W^{\bot}}{W}.
\]
Suppose further that $ (  V,\sigma )  $ is a filtered quadratic
$K [  G ]  $-space in the following sense: we are given an increasing
filtration $ \{  F_{i} \}  $ of $K[G]$-submodules with
$F_{-N}= (  0 )  $ and $F_{N}=V$ for $N>>0,$ and with $F_{i}^{\perp
}=F_{-i-1};$ thus for all $i$, $\sigma$ induces isomorphisms
\[
{\rm Gr}_{-i}\cong {\rm Gr}_{i}^{D}%
\]
where ${\rm Gr}_{i}$ denotes the $i$-th graded piece $F_{i}/F_{i-1}$. Then there is
a (non-canonical) $K [  G ]  $-decomposition of quadratic modules
\[
V\cong\oplus_{i<0}\text{\rm Hyp} (  {\rm Gr}_{i} )  \oplus {\rm Gr}_{0}.
\]
\end{lemma}

\noindent{\sl Proof.} To prove the first part for simplicity we may suppose without loss of
generality that $W$ is irreducible. First choose an arbitrary decomposition of
$K [  G ]  $-modules $W^{\bot}=W\oplus U$; this is trivially an
orthogonal decomposition. We then choose an arbitrary further decomposition of
$K\left[  G\right]  $-modules $V=W^{\bot}\oplus W^{\prime}$.   Then the
form $\sigma$ induces a map
\begin{equation}
W^{\prime}\overset{\sigma^{\prime}}{\rightarrow}U^{D}\oplus W^{D}%
\overset{\text{proj}}{\rightarrow}W^{D}%
\end{equation}
and the composite is an isomorphism. We may then alter the initial
decomposition $V=W^{\bot}\oplus W^{\prime}$ by a homomorphism from $W^{\prime
}$ to $U$ to guarantee that the composition of $\sigma^{\prime}$  with
projection to $U^{D}$ is zero, as required.

The second part of the lemma then follows at once from the first part.
\hfill $\square$

\begin{proposition}
 Let $U$ be a free $\mathbf{R} [  G ]  $-module with basis,
$ \{  u_{i} \} $ $i=1,\ldots,q$, and suppose that $U$ supports a
non-degenerate real-valued $G$-invariant form $\sigma$. Choose a
decomposition $U=U^{+}\oplus U^{-}$, as in Proposition 25, and for each $m$
define $n_{m}^{\pm}\left(  \sigma\right)  =\dim (  U^{\pm}\otimes
W_{m} )  ^{G}$. Then
\[
\text{\rm sign}\left(  \mathbf{pf}_{\kappa_{m}^{\left(  q\right)  }} (
T_{W_{m}}^{\left(  q\right)  }\left(  \widetilde{\sigma}\left(  u_{i}%
,u_{j}\right)  \right)   )  \right)  = (  \sqrt{-1} )
^{n_{m}^{-} (  \sigma )  }.
\]
Note that the integers $n_{m}^{\pm}\left(  \sigma\right)  $ are all even,
since they are the multiplicities of symplectic representations in real representations.
\end{proposition}

\begin{remark}
{\rm The authors are grateful to Boas Erez for stressing the importance of
relating the sign of Fr\"{o}hlich's Pfaffian to signature invariants.}
\end{remark}

\noindent{\sl Proof.} As previously, let $\nu_{\mathbf{R}}$ denote the standard $G$-invariant
form on $\mathbf{R} [  G ]  $ (see Lemma 21); we again define the
form $\nu$ on $U$ by the rule
\[
\nu (  \lambda u_{i},\mu u_{j} )  =\delta_{ij}\nu_{\mathbf{R} [
G ]  }\left(  \lambda,\mu\right)  .
\]
As in Proposition 25, there is a unique $\mathbf{R} [  G ]
$-automorphism $A$ of $U$, which is self-adjoint with respect to $\nu$, such
that for all $x$,$y\in U$
\[
\sigma (  x,y )  =\nu (  x,Ay )  .
\]
Therefore the decomposition $U=U^{+}\oplus U^{-}$ induces a decomposition
\[
\left(  U\otimes W_{m}\right)  ^{G}=\left(  U^{+}\otimes W_{m}\right)
^{G}\oplus\left(  U^{-}\otimes W_{m}\right)  ^{G}%
\]
and $A\otimes1$ is diagonalisable on the subspaces $\left(  U^{\pm}\otimes
W_{m}\right)  ^{G}$ with positive resp. negative eigenvalues on $\left(
U^{+}\otimes W_{m}\right)  ^{G}$ resp. $\left(  U^{-}\otimes W_{m}\right)
^{G}$. By construction $ (  A\otimes 1 )  ^{G}$ induces an isometry
of alternating forms
\[
\left(  \sigma\otimes\kappa_{m}\right)  ^{G}\cong\left(  \nu\otimes\kappa
_{m}\right)  ^{G}%
\]
which we denote $ (  A\otimes 1 )  _{m}^{G}$. Therefore by Proposition
5 we see that
\[
\text{Pf}_{\left(  \sigma\otimes\kappa_{m}\right)  ^{G}}=\mathbf{pf}_{\left(
\nu\otimes\kappa_{m}\right)  ^{G}} (   (  A\otimes1 )  _{m}%
^{G} )\cdot  \text{Pf}_{\left(  \nu\otimes\kappa_{m}\right)  ^{G}}.
\]
We then evaluate both sides on $\wedge_{in}\left|  G\right|  ^{-1}r_{G}\left(
u_{i}\otimes w_{mn}\right)  $ and, noting that by Lemma 21 $ \{|G|^{-1}r_{G}\left(  u_{i}\otimes w_{mn}\right)   \}  $ is
a hyperbolic basis with respect to $\left(  \nu\otimes\kappa_{m}\right)  ^{G}%
$, we see that
\[
\text{Pf}_{\left(  \sigma\otimes\kappa_{m}\right)  ^{G}}\left(  \wedge
_{in}\left|  G\right|  ^{-1}r_{G}\left(  u_{i}\otimes w_{mn}\right)  \right)
=\mathbf{pf}_{\left(  \nu\otimes\kappa_{m}\right)  ^{G}} (  \left(
A\otimes1\right)  _{m}^{G} )  .
\]
On the other hand by Proposition 23
\[
\text{Pf}_{\left(  \sigma\otimes\kappa_{m}\right)  ^{G}}\left(  \wedge
_{in}\left|  G\right|  ^{-1}r_{G}\left(  u_{i}\otimes w_{mn}\right)  \right)
=\mathbf{pf}_{\kappa_{m}^{\left(  q\right)  }}(  \left|  G\right|
^{-1}T_{W_{m}}^{\left(  q\right)  }\left(  \widetilde{\sigma}\left(
u_{i},u_{j}\right)  \right)  )  .
\]
The result then follows by repeated use of the fact that (see page 40 in
[F1])
\begin{equation*}
\mathbf{pf}\left(
\begin{array}
[c]{cc}%
d & 0\\
0 & d
\end{array}
\right)  =d.  \ \ \ \ \square  
\end{equation*} 

Recall that ${\rm S}_{\infty} (  \mathbf{Z} [  G ]  )  $ is the
signature group defined in 3.1.2. We complete this section by stating the
following result, which we shall prove in (5.3) of the Appendix.

\begin{theorem}
The class $\chi_{\text{\rm H}}^{\text{\rm s}} (  P^{\bullet},\sigma )
\chi_{\text{\rm A}}^{\text{\rm s}} (  P^{\bullet},h_{\sigma} )  ^{-1}$ lies
in ${\rm S}_{\infty} (  \mathbf{Z} [  G ]  )  $ and is
represented by the character function which maps $\theta_{m}$ to $i^{n_{m}%
^{-} (  \sigma )  }$ where
\[
n_{m}^{-} (  \sigma )  =n_{m}^{-} (  \sigma^{\rm ev} )
-n_{m}^{-} (  \sigma^{\rm odd} )
\]
denotes the virtual dimension of a maximal negative definite subspace of the
 $W_{m}$-isotypic component of the cohomology of $P^{\bullet}$.
\end{theorem}

\section{ DE RHAM DISCRIMINANTS.}

Throughout this section we again adopt the notation given in the Introduction.
Thus the scheme $\mathcal{X}$ is projective and flat over ${\rm Spec} (
\mathbf{Z} )  $ of relative dimension $d$ and  $\pi:\mathcal{X\rightarrow Y}$ is a $G$-cover which
satisfies hypotheses (T1) and (T2) given in the Introduction.  In most of this
section we again let $X=\mathcal{X}\times_{\text{Spec}\left(  \mathbf{Z}%
\right)  }{\rm Spec}\left(  \mathbf{Q}\right)  $ resp. $Y=\mathcal{Y}%
\times_{\text{Spec}\left(  \mathbf{Z}\right)  }{\rm Spec}\left(  \mathbf{Q}%
\right)  $ denote the generic fibre of $\mathcal{X}$  resp. $\mathcal{Y}$.

\subsection{De Rham pairings.}

As in III Sect. 7 of [H] for $0\leq i,j\leq d$ we have $G$-equivariant
duality pairings 
\begin{equation*}
\sigma _{ij}:\text{\textrm{H}}^{i} ( X,\Omega _{X}^{j} ) \times 
\text{\textrm{H}}^{d-i} ( X,\Omega _{X}^{d-j} ) \overset{\cup }{%
\rightarrow }\mathrm{H}^{d} ( X,\Omega _{X}^{d} ) \overset{\left|
G\right| ^{-1}{\rm Tr}}{\rightarrow }\mathbf{Q}
\end{equation*}
where $\rm Tr$ is the trace map described after (17) below; note that here we
divide the pairings used in [CPT1] by the group order; we used this
normalisation in [CPT2], and the reason for choosing this normalisation will
be explained after Lemma 32 below.  For arbitrary coherent $X$-sheaves $%
\mathcal{F}$, $\mathcal{G}$  and for $x\in \mathrm{H}^{i}( X,\mathcal{%
F})$, $y\in \mathrm{H}^{j} ( X,\mathcal{G} )$, we know that $%
x\cup y= ( -1 ) ^{ij}y\cup x$ after identifying $\mathcal{F}\otimes 
\mathcal{G}$ and $\mathcal{G}\otimes \mathcal{F}$ by the ``flip''
isomorphism. Taking $\mathcal{F}=\Omega _{X}^{a}$, $\mathcal{G}=\Omega
_{X}^{b}$ we see that $x\cup y= ( -1 ) ^{ij+ab}y\cup x$ in $%
\mathrm{H}^{i+j} ( X,\Omega _{X}^{a+b} ) $; it therefore follows
that 
\begin{equation}
\sigma _{i,j} ( x,y ) = ( -1) ^{ ( d+1 )  (
i+j ) }\sigma _{d-i,d-j} ( y,x ) .
\end{equation}
We then symmetrise these pairings by the construction given in Sect. 3 of
[CPT1]: namely, we define  the twisted pairing $\sigma _{i,j}^{\prime }$ 
\begin{equation}
\sigma _{d-i,d-j}^{\prime }( y,x) =\sigma _{i,j}( x,y)
=( -1) ^{( d+1) ( i+j) }\sigma
_{d-i,d-j}( y,x) .
\end{equation}
We then define pairings $\sigma^{t}$ on the hypercohomology of $R\Gamma (
X, L\wedge^{\bullet}\Omega^1_{X/\mathbf{Q}} )   [  d ]$ as follows:
\[
\text{for\ \ }t<0\;\;\text{we put\ \ }\sigma^{t}=\oplus_{i+j=t+d}%
\;\sigma_{i,j}%
\]%
\[
\text{for\ \ }t>0\;\;\text{we put\ \ }\sigma^{t}=\oplus_{i+j=t+d}%
\;\sigma_{i,j}^{\prime}%
\]
and for $t=0$ we set
\[
\sigma^{0}=\oplus_{i<d/2}\;\sigma_{i,d-i}\oplus\sigma_{d/2,d/2}\oplus
_{i>d/2}\;\sigma_{i,d-i}^{\prime}.
\]
Here it is to be understood that the term $\sigma_{d/2,d/2}$ occurs only when
$d$ is even. We note that in all cases $\sigma^{t}\;$is symmetric by (15) and
(16), and we then define
\[
\sigma^{\text{ev}}=\oplus_{t\text{ even}}\sigma^{t},\;\;\sigma^{\text{odd}%
}=\oplus_{t\text{ odd}}\sigma^{t}.
\]
Note that in all cases $\sigma^{\text{odd}}$ is a hyperbolic pairing; and
$\sigma^{\text{ev}}$ is hyperbolic whenever $d$ is odd. To be more precise we have:

\begin{proposition}
There is a $\mathbf{Q}[  G]  $-isometry
\[
\left(  \text{\rm H}^{i} (  X,\Omega_{X}^{j} )  \oplus\text{\rm H}%
^{d-i} (  X,\Omega_{X}^{d-j} )  ,\sigma_{ij}\oplus\sigma
_{d-i,d-j}^{\prime}\right)  =\text{\rm Hyp}\left(  \text{\rm H}^{i} (
X,\Omega_{X}^{j} )  \right)
\]
unless $d$ is even and $i=j=d/2$.
\end{proposition}

Since the signature of a hyperbolic form is always zero we know

\begin{lemma}
\bigskip For any symplectic representation $W_{m}$ of $G$
\[
n_{m}^{+} (  \sigma )  -n_{m}^{-} (  \sigma )  =n_{m}
^{+} (  \sigma_{d/2,d/2} )  -n_{m}^{-} (  \sigma_{d/2,d/2})
\]
where the right hand side is to be interpreted as zero if $d$ is odd.
\end{lemma}

We conclude this subsection by considering the complex hermitian form
associated to the de Rham pairing $\sigma$ by 3.2.2. We begin by considering the 
$L^{2}$-norms on cohomology.

Given a K\"{a}hler metric $h_{Y}$ on the complex tangent space of an
arithmetic variety $\mathcal{Y}$ which is invariant under complex conjugation,
we denote by $h_{X}=h^{TX}$ the K\"{a}hler metric on $X(\C)$ given by the
pullback of $h_{Y}$; this is then also invariant under complex conjugation.
Define $h_{X}^{D}$ to be the metric on the complex cotangent space of $X(\C)$
which is dual to $h_{X}$.

Let $d_{X}$ denote the volume form given by the $d$-th exterior power of the
$\left(  1,1\right)  $-form associated to $h_{X}^{D}$. Define the $L^{2}%
$-metric on $\Omega_{X(\C)}^{p,0}\otimes\Omega_{X(\C)}^{0,q}$ by
\[
 \langle s,t \rangle _{X}=\frac{1}{d!}\int_{X(\C)}\left|  G\right|
^{-1}\wedge^{p+q}h_{X}^{D} (  s (  x )  ,t (  x )
 )  \left(  \frac{i}{2\pi} \right)  ^{d}d_{X}%
\]
where $\wedge^{p+q}h_{X}^{D} (  -,- )  $ denotes the inner product on
$p+q$ forms given by the $p+q$-th exterior product of $h_{X}^{D}$ (see for
instance page 131 in [So]). The reason for the normalisation factor $\left(
i/2\pi\right)  ^{d}$ on the volume form will become apparent below: basically
it will ensure that the $L^{2}$-norm is compatible with Serre duality. The
reason for normalising by the factor $\left|  G\right|  ^{-1}$ is that, since
$X\rightarrow Y$ is etale, our metrics are then natural with respect to
pullback in the sense that for $p$-forms $s^{\prime},t^{\prime}$ on $Y$, we
then have $ \langle \pi^{\ast}s^{\prime},\pi^{\ast}t^{\prime} \rangle
_{X}= \langle s^{\prime},t^{\prime} \rangle _{Y}$ where
\[
 \langle s^{\prime},t^{\prime} \rangle _{Y}=\frac{1}{d!}\int
_{Y(\C)}\wedge^{p+q}h_{Y}^{D} (  s^{\prime}\left(  y\right)  ,t^{\prime
}\left(  y\right)   )  \left(  \frac{i}{2\pi}\right)  ^{d}d_{Y}%
\]
and $d_{Y}$ is the volume form given by the $d$-th exterior power of the
$\left(  1,1\right)  $-form associated to $h_{Y}^{D}$. Let $\Delta
^{q}=\overline{\partial}\overline{\partial}^{\ast}+\overline{\partial}^{\ast
}\overline{\partial}$ be the Laplace operator on $\Omega_{X(\C)}^{p,0}%
\otimes\Omega_{X(\C)}^{0,q}$. The Hodge isomorphism
\[
\text{H}^{q} (  X(\C),\Omega_{X(\C)}^{p,0} )  =\ker\left(
\Delta^{q}\right)
\]
then gives an $L^{2}$-metric on H$^{q} (  X(\C),\Omega_{X(\C)}%
^{p,0} )  $. We will denote the resulting $L^{2}$-metric on the
determinant of cohomology of $\Omega_{X(\C)}^{p,0}$ by $\left|
G\right|  ^{-1}\wedge^{\bullet}\left| \ .\ \right|  _{L^{2}}$ in order to
emphasise the appearance of the scaling factor $\left|  G\right|  ^{-1}$. We
then construct the associated Quillen metrics on the equivariant determinant
of cohomology of $\Omega_{X(\C)}^{p,0}$ by multiplying the above $L^{2}$-metrics
by the inverse of the equivariant analytic torsion associated to $\left|
G\right|  ^{-1}\wedge^{p}h_{X}^{D}.$ This construction is described in more
detail in the proof of Proposition 34. For a full discussion of this
construction see Section 6 of [CPT2].

Identifying ${\rm H}^{d} (  X_\C,\Omega_{X_{\mathbf{C}}}^{d} )  $ with the
Dolbeault cohomology group ${\rm H}_{\overline{\partial}}^{d,d} (  X )  $
and then integrating over $X $ affords a surjection ${\rm H}^{d} (
X,\Omega_{X}^{d} )  \otimes\mathbf{C}={\rm H}^{d} 
(  X_\C,\Omega_{X_{\mathbf{C}}}^{d} )\overset{\int_{X}%
}{\rightarrow}\mathbf{C}$. From the above discussion we know that the
following diagram commutes
\begin{equation}%
\begin{array}
[c]{ccc}%
\text{H}^{d} (  Y,\Omega_{Y}^{d} )  & \overset{\int_{Y}}{\rightarrow}%
& \mathbf{C}\\
\pi^{\ast}\downarrow &  &\ \  \downarrow=\\
\text{H}^{d} (  X,\Omega_{X}^{d} )  & \overset{\left|  G\right|
^{-1}\int_{X}}{\rightarrow} & \mathbf{C}%
\end{array}
\end{equation}
We then define the trace map
\[
\text{H}^{d} (  X,\Omega_{X}^{d} )  \overset{\left|  G\right|
^{-1}{\rm Tr}}{\longrightarrow}\mathbf{Q}%
\]
to be induced by the map $\displaystyle{\frac{i^{d}}{\left(  2\pi\right)  ^{d}d!\left|
G\right|  }}\int_{X}$ .  Recall that the following diagram commutes up to sign
(see page 102 in [GH])
\[
\begin{array}
[c]{ccccc}
\text{H}^{i} (  X,\Omega_{X}^{j} )  & \times & \text{H}^{d-i} (
X,\Omega_{X}^{d-j} )  & \overset{\cup}{\rightarrow} & \text{H}%
^{d} (  X,\Omega^{d}_X ) \\
\downarrow &  & \downarrow &  & \downarrow\\
{\rm H}_{\overline{\partial}}^{i,j} (  X )  & \times &  {\rm H}_{\overline
{\partial}}^{d-i,d-j} (  X )  & \overset{\wedge}{\rightarrow} &
{\rm H}_{\overline{\partial}}^{d,d} (  X )
\end{array}
\]
where the upper horizontal map is cup product, the lower horizontal  map is
the exterior product of differential forms, and the vertical arrows are
Dolbeault isomorphisms; 
in fact we shall compute this sign in 4.3.2. It now follows that the metrics $h_{\sigma}= \{
h_{r} \}  $ on the equivariant determinant of cohomology induced by
$\sigma$ coincide with the metrics on the equivariant determinant of
cohomology induced by $|G|^{-1} |\ .\ |_{L^{2}}$.
(See  1.4 in [GSZ], and especially Theorem 7.8 in [CPT2] for a full
account of the duality and metrics.) In summary we have now shown

\begin{lemma}
The metrics $h_{\sigma}= \{  h_{r} \}  $, associated to the de Rham
pairings $\sigma_{i,j}$ by 3.2.2, coincide with the $L^{2}$-metrics, and so we have the
equality of Arakelov Euler characteristics
\begin{equation}
\chi_{\text{\rm A}} (  R\Gamma (  \mathcal{X},L\wedge^{\bullet}%
\Omega_{{\mathcal X}/\Z}^{1} )  [  d ]  ,  h_{\sigma} )
=\chi_{\text{\rm A}} (  R\Gamma (  \mathcal{X},L\wedge^{\bullet}%
\Omega_{{\mathcal X}/\Z}^{1} )  [  d ]  ,   |  G |
^{-1} |  . |_{L^{2}} )  .
\end{equation}
\end{lemma}

Next we use a result of Ray-Singer (see Theorem 3.1 of [RS])  to show that:

\begin{proposition}
The equivariant analytic torsion of the total de Rham complex vanishes (see
below); thus we can write
\begin{equation}
\chi_{\text{\rm A}} (  R\Gamma (  \mathcal{X},L\wedge^{\bullet}%
\Omega_{{\mathcal X}/\Z}^{1} )  [  d ]  ,   |  G |
^{-1}\left|  .\right|  _{L^{2}} )  =\chi_{\text{\rm A}} (  R\Gamma (
\mathcal{X},L\wedge^{\bullet}\Omega_{{\mathcal X}/\Z}^{1} )   [
d ]  ,  |  G |  ^{-1}\wedge^{\bullet}h_{X,Q}^{D} )
\end{equation}
where $ |  G |  ^{-1}\wedge^{\bullet}h_{X,Q}^{D}$ denotes the
equivariant Quillen metrics\ (see [B]) on the equivariant determinant of
cohomology induced by the $ |  G |  ^{-1}\wedge^{p}h_{X}^{D}$.
\end{proposition}

\noindent{\sl Proof.} For an irreducible character $\phi$ of $G$, we let $T_{\phi} (
\Omega_{X(\C)}^{p}, |  G |  ^{-1}\wedge^{p}h_{X}^{D} )  $
denote the analytic torsion associated to the hermitian sheaf $ (
\Omega_{X/(\C)}^{p}, |  G |  ^{-1}\wedge^{p}h_{X}^{D} )
$. Thus, by definition,
\[
(\left|  G\right|  ^{-1}\wedge^{\bullet}h_{X}^{D})_{Q,\phi}=T_{\phi} (
\Omega_{X/\mathbf{Q}}^{p},\left|  G\right|  ^{-1}\wedge^{p}h_{X}^{D} )
^{-1} (  \left|  G\right|  ^{-1}\left|  .\right|  _{L^{2}} )  _{\phi
}.
\]
From Theorem 3.1 in [RS] we know that
\[%
{\textstyle\prod_{p=0}^{d}}
T_{\phi} (  \Omega_{X(\C)}^{p},\wedge^{p}h_{X}^{D} )
^{\left(  -1\right)  ^{p}}=1.
\]
Moreover, it is standard (see for instance page 153 in [R]) that, if we scale
the metrics $\wedge^{p}h_{X}^{D}$ to $c^{2}\wedge^{p}h_{X}^{D}$ for a positive
real number $c$, then the total analytic torsion changes by a factor $c$ to
the power
\[%
{\textstyle\sum_{p,q}}
\left(  -1\right)  ^{p+q}q\zeta_{p,q}\left(  0,\phi\right)
\]
where $\zeta_{p,q}\left(  s,\phi\right)  $ denotes the $\zeta$-function for
$\phi$ associated to $\Omega_{X(\C)}^{p,q}$ with the metric induced by
$\wedge^{p+q}h_{X}^{D}$. However, from equation (3.2) in the proof of Theorem
3.1 in [RS] we know that for each $q$%
\[%
{\textstyle\sum_{p}}
\left(  -1\right)  ^{p}\zeta_{p,q}\left(  s,\phi\right)  =0
\]
and so we deduce that
\[%
{\textstyle\sum_{p,q}}
\left(  -1\right)  ^{p+q}q\zeta_{p,q}\left(  0,\phi\right)  =0.\;\;\;\square
\]

\bigskip

Recall that the principal goal of this paper is to describe the hermitian
Euler characteristic
\[
\chi_{\text{H}}^{\text{s}} (  R\Gamma({\mathcal X},  L\wedge^{\bullet}\Omega
_{\mathcal{X}/\mathbf{Z}},\sigma_X )  \in\text{H}^{\text{s}} (
\mathbf{Z} [  G ]   )  ;
\]
while from (5) of 3.1.2 we have the decomposition
\[
\text{H}^{\text{s}} (  \mathbf{Z} [  G ]   )  =\text{A}%
^{\text{s}} (  \mathbf{Z} [  G ]   )  \oplus\text{S}%
_{\infty} (  \mathbf{Z} [  G ]   )  .
\]
From Theorem 30 we know that the image of $\chi_{\text{H}}^{\text{s}} (
R\Gamma (  \mathcal{X},L\wedge^{\bullet}\Omega_{{\mathcal X}/\Z}^{1} )
 [  d ]  ,\sigma_X )  $ in ${\rm A}^{\text{s}} (  \mathbf{Z} [
G ]   )  $ is the Arakelov class $\chi_{\text{A}}^{\text{s}} (
R\Gamma (  \mathcal{X},L\wedge^{\bullet}\Omega_{{\mathcal X}/\Z}^{1} )
 [  d ]  ,h_{\sigma} )  $ where $h_{\sigma}$ denotes the metrics
on the equivariant  determinant of cohomology afforded by the absolute values
of the equivariant determinants of $\sigma=\sigma_X$. With the above choices we have
seen that $h_{\sigma}$ coincides with the $ |  G |  ^{-1}%
\wedge^{\bullet} |  . |  _{L^{2}}$-norm and so by (18) and (19)
\[
\chi_{\text{A}}^{\text{s}} (  R\Gamma (  \mathcal{X},L\wedge^{\bullet
}\Omega_{{\mathcal X}/\Z}^{1} )   [  d ]  ,h_{\sigma} )
=\chi_{\text{A}}^{\text{s}} (  R\Gamma (  \mathcal{X},L\wedge
^{\bullet}\Omega_{{\mathcal X}/\Z}^{1} )  [  d ]  , \left|
G\right|  ^{-1}\left|  .\right|  _{L^{2}} )=
\]%
\[
=\chi_{\text{A}}^{\text{s}} (  R\Gamma(  \mathcal{X},L\wedge
^{\bullet}\Omega_{{\mathcal X}/\Z}^{1} )   [  d ]  , \left|
G\right|  ^{-1}\wedge^{\bullet}h_{X,Q}^{D} )  .
\]
The crucial point here is that the latter Arakelov class was determined in
Theorem 8.4 in [CPT2]: in particular, we show that on symplectic characters
$\theta$ of degree zero this class characterises the global constant
$\varepsilon (  \mathcal{X},\theta )  $. Thus, to complete our
description of the class $\chi_{\text{H}}^{\text{s}} (  R\Gamma (
\mathcal{X}, L\wedge^{\bullet}\Omega_{{\mathcal X}/\Z}^{1} )   [
d ]  , \sigma )  $, it remains to describe the image $\chi_{2}$ of
the class $\chi_{\text{H}}^{\text{s}} (  R\Gamma (  \mathcal{X}, L\wedge
^{\bullet}\Omega_{{\mathcal X}/\Z}^{1} )  [  d ]  ,\sigma )$ 
in ${\rm S}_{\infty} (  \mathbf{Z} [  G ] )$. From Theorem
30 we know that for $\theta_{m}$ in $R_{G}^{\text{s}}$%
\begin{equation}
\chi_{2} (  \theta_{m} )  =i^{n_{m}^{-} (  \sigma )  }.
\end{equation}
Thus, in order to prove Theorem 1, we are now required to show that
\begin{equation}
i^{\delta (  Y )  \theta_{m} (  1 )  }\varepsilon_{\infty
} (  \mathcal{X},\theta_{m} )  =i^{n_{m}^{-} (  \sigma )  }.
\end{equation}

\subsection{Archimedean $\varepsilon$-constants.}

Here we recall a number of results from Sect. 5 of [CEPT2]. Let $F_{\infty
}: X (  \mathbf{C} )  \rightarrow X (  \mathbf{C} )  $ denote
the involution induced by complex conjugation on $X (  \mathbf{C} )$,
the space of complex points of $\mathcal{X}$; then $F_{\infty}$ acts on
the Betti cohomology ${\rm H}_{B}^{i} (  X (  \mathbf{C} )
,\mathbf{Q} )  $ and, for a complex representation $V$ of $G$ with
contragredient $V^{\ast}$, we write ${\rm H}_{B+}^{i} (  X\otimes_{G}V )  $
resp. ${\rm H}_{B-}^{i} (  X\otimes_{G}V )  $ for the subspace of $ (
\text{H}_{B}^{i} (  X\left(  \mathbf{C}\right)  ,\mathbf{Q} )
\otimes_{\mathbf{Q}}V^{\ast} )  ^{G}$ on which $F_{\infty}$ acts by $+1$
resp. $-1$. (For a discussion of the motives $X\otimes_{G}V$ see Section 2 of
[CEPT2].) We then set
\[
\chi_{\pm}\left(  X\otimes_{G}V\right)  =%
{\textstyle\sum_{i=0}^{2d}}
\left(  -1\right)  ^{i}\dim_{\mathbf{C}} (  \text{H}_{B\pm}^{i} (
X\otimes_{G}V )  )
\]
and we may extend $\chi_{\pm} (  X\otimes_{G}V )  $ to virtual
representations, since it is additive in $V$.

The archimedean constant $\varepsilon_{\infty} (  X\otimes_{G}V )  $
is constructed from the Hodge structure of the motive $X\otimes_{G}V$; again
it is additive in $V$ and thus extends to virtual $V$.

\begin{lemma}
\bigskip Let $ W$ be a virtual symplectic complex representation of $G$.

(a) Both $\chi_{\pm} (  X\otimes_{G}W )  $ are even integers.

(b) If $d$ is odd, then $\varepsilon_{\infty} (  X\otimes_{G}W )  =1.$

(c) If $d$ is even, then writing $\pm$ for the sign of $\left(  -1\right)
^{d/2+1}$ we have
\[
\varepsilon_{\infty} (  X\otimes_{G}W )  =i^{\chi_{\pm} (
X\otimes_{G}W )  }%
\]

and, moreover if $\dim_{\mathbf{C}} (  W )  =0,$ then $\varepsilon
_{\infty} (  X\otimes_{G}W )  =i^{\chi_{+} (  X\otimes
_{G}W )  }=i^{\chi_{-} (  X\otimes_{G}W )  }$.
\end{lemma}

\noindent{\sl Proof.} Part (a) follows from the discussion in 3.2.1 after Lemma 15 which shows
that each $\dim_{\mathbf{C}} (  \text{H}_{B\pm}^{i} (  X\otimes
_{G}W )   )  $ is even; (b) and (c) come from Lemma 5.1.1 in
[CEPT2].\ \ \ $\square$

\bigskip

\subsection{Signature of cohomology.}

Throughout all of this sub-section we shall suppose that the fibral dimension
$d$ is \textit{even}.

\subsubsection{\bigskip Betti cohomology.}

Since $d$ is even, $X (  \mathbf{C} )  $ has real dimension divisible
by 4; hence the cup product $c^{d}$ is a non-degenerate symmetric
$G$-invariant form on ${\rm H}_{B}^{d} (  X (  \mathbf{C} )
,\mathbf{R} )  $ via the map ${\rm H}_{B}^{2d} (  X (  \mathbf{C}%
 )  ,\mathbf{R} )  \rightarrow\mathbf{R}$. By Proposition 25 we
know that ${\rm H}_{B}^{d} (  X (  \mathbf{C} )  , \mathbf{R} )
$ admits a non-canonical decomposition of $G$-modules 
\[
\text{H}_{B}^{d} (  X (  \mathbf{C} )  ,\mathbf{R} )
=\text{H}_{B}^{+}\oplus\text{H}_{B}^{-}%
\]
where ${\rm H}_{B}^{+}$ is a maximal positive definite subspace and ${\rm H}_{B}^{-}$ is
a maximal negative definite subspace of ${\rm H}_{B}^{d} (  X (
\mathbf{C} )  ,\mathbf{R} )  $ with respect to $c^{d}$.

For $t<d$ we let $c^{t}$ denote the symmetrised $G$-invariant form on
 ${\rm H}_{B}^{t} (  X (  \mathbf{C} )  ,\mathbf{R} )
\oplus{\rm H}_{B}^{2d-t} (  X (  \mathbf{C} )  ,\mathbf{R} )  $
induced by the cup product
\[
\text{H}_{B}^{t}  ( X (  \mathbf{C} )  ,\mathbf{R} )
\times\text{H}_{B}^{2d-t} (  X (  \mathbf{C} )  ,\mathbf{R}%
 )  \rightarrow\text{H}_{B}^{2d} (  X (  \mathbf{C} )
,\mathbf{R} )  \rightarrow\mathbf{R}%
\]
as per the construction of $\sigma^{t}$ in 4.1. Note that the symmetrisation
here is the same as that used in 4.1: indeed, for $x\in{\rm H}_{B}^{t} (
X (  \mathbf{C} )  ,\mathbf{R} )$, $y\in {\rm H}_{B}%
^{2d-t} (  X(  \mathbf{C} )  ,\mathbf{R} )$,
 $t<d$
\[
c\left(  y,x\right)  =\left(  -1\right)  ^{t}c\left(  x,y\right)  .\text{ }%
\]
Whereas by (15) for $w\in {\rm H}^{i} (  X,\Omega_{X}^{j} )$, 
$ z\in {\rm H}^{d-i} (  X,\Omega_{X}^{d-j} )$, if we set $t=i+j$, then, as
$d$ is even, we have seen that
\[
\sigma_{d-i,d-j}\left(  z,w\right)  =\left(  -1\right)  ^{\left(  i+j\right)
}\sigma_{i,j}\left(  w,z\right)  =\left(  -1\right)  ^{t}\sigma_{i,j}\left(
w,z\right)  .
\]
Thus for $t<d$, $c^{t}$ is hyperbolic and by Proposition 25 we have a
decomposition of $\mathbf{R} [  G ]  $-modules
\[
\text{H}_{B}^{\text{odd}} (  X (  \mathbf{C} )  ,\mathbf{R}%
 )  =\text{H}_{B}^{\text{odd}+}\oplus\text{H}_{B}^{\text{odd}-}%
\]
into positive and negative subspaces.

Applying Proposition 25 once again we obtain a decomposition
\[
\text{H}_{B}^{\text{ev}} (  X (  \mathbf{C} )  ,\mathbf{R}%
 )  =\text{H}_{B}^{\text{ev}+}\oplus\text{H}_{B}^{\text{ev}-}%
\]
where ${\rm H}_{B}^{+}\subset{\rm H}_{B}^{\text{ev}+}$, ${\rm H}_{B}^{-}\subset
{\rm H}_{B}^{\text{ev}-}.$

Furthermore, by Lemma 26 and by hyperbolicity, we know that as $\mathbf{R}%
$-vector spaces
\begin{align*}
\text{H}_{B}^{\text{ev}+}/\text{H}_{B}^{+}  & \cong\text{H}_{B}^{\text{ev}%
-}/\text{H}_{B}^{-}\text{ }\cong\oplus_{t\text{ even, }t<d}\text{H}_{B}%
^{t} (  X (  \mathbf{C} )  ,\mathbf{R} ), \\
\text{ \ H}_{B}^{\text{odd}+}  & \cong\text{H}_{B}^{\text{odd}-}\cong
\oplus_{t\text{ odd, }t<d}\text{H}_{B}^{t} (  X (  \mathbf{C} )
, \mathbf{R} )  .
\end{align*}

\begin{theorem}
\bigskip With the above notation and hypotheses, ${\rm H}_{B}^{\bullet+}$ and
${\rm H}_{B}^{\bullet-}$ are both free virtual $\mathbf{R} [  G ]  $-modules.
\end{theorem}

\noindent{\sl Proof.} Since $G$ acts freely on $X (  \mathbf{C} )  $, by the
Lefschetz Fixed Point theorem (see for instance [V]) for each $g\in G$, $g\neq 1$,
the virtual character associated to ${\rm H}_{B}^{\bullet} (  X (
\mathbf{C} )  ,\mathbf{R} )  $ is zero when evaluated on such $g$;
thus ${\rm H}_{B}^{\bullet}={\rm H}_{B}^{\bullet+}+{\rm H}_{B}^{\bullet-}$ is a free virtual
$\mathbf{R} [  G ]  $-module.

Similarly we shall show that ${\rm H}_{B}^{\bullet+}-{\rm H}_{B}^{\bullet-}$ is a free
virtual $\mathbf{R} [  G ]  $-module; this will then establish the
theorem. To see that ${\rm H}_{B}^{\bullet+}-{\rm H}_{B}^{\bullet-}$ is free, we recall
that by the $G$-Signature Theorem 6.12 in [AS] (see also V.18 in [S]), for
each non-trivial element $g\in G$, the value of the virtual character of
${\rm H}_{B}^{\bullet+}-{\rm H}_{B}^{\bullet-}$ evaluated on $g$ is presented in terms
of data associated to the fixed point set $X (  \mathbf{C} )  ^{g}$.
Since $g$ acts without fixed points, it then follows that this virtual
character is zero on all such $g$.\hfill $\square$

\subsubsection{De Rham cohomology and hypercohomology.}

In this paragraph, to ease the notation, we use the symbol $X$ to denote either 
${\mathcal X}_\Q={\mathcal X}\times_{\Spec(\Z)}\Spec(\Q)$ or  ${\mathcal X}_{\R}=
{\mathcal X}\times_{\Spec(\Z)}\Spec(\R)$ depending on the context. 
We follow the terminology of Grothendieck (see [G]); for a given integer
$t$, we consider the (shifted) $t$-th Hodge cohomology group
\[
\text{H}_{\rm Hod}^{t} (  X )   [  d ]  =\text{H}^{t} (
X,\oplus_{n}\Omega_{X }^{n} [  d-n ]   )  =\oplus
_{n}\text{H}^{t+d-n} (  X,\Omega_{X }^{n} )  =\oplus
_{m:\;m+n=t+d}\text{H}^{m} (  X,\Omega_{X }^{n} )
\]
and similarly we put
\[
\text{H}_{\rm Hod}^{\text{ev}} (  X )   [  d ]  =\oplus
_{t\text{ even}}\text{H}_{\rm Hod}^{t} (  X )   [  d ]
,\quad \text{H}_{\rm Hod}^{\text{odd}} (  X )  [  d ]
=\oplus_{t\text{ odd}}\text{H}_{\rm Hod}^{t} (  X )   [  d ]  .
\]
We then let $ (  \text{H}_{\rm Hod}^{\text{ev}} (  X )   [
d ]  ,\sigma^{\text{ev}} )  $ denote ${\rm H}_{\rm Hod}^{\text{ev}} (
X )   [  d ]  $ endowed with the $G$-invariant symmetric form
$\sigma^{\text{ev}}$ and similarly we have $(  \text{H}_{\rm Hod}%
^{\text{odd}} (  X )   [  d ]  ,\sigma^{\text{odd}} )$. In what follows,  
we take $X=X_{\bf R}$.
 Applying Proposition 25 we have a decomposition of $\mathbf{R}%
 [  G ]  $-modules into positive and negative spaces
\[
\text{H}_{\rm Hod}^{\text{ev}}(X)[d]%
=\text{H}_{\rm Hod}^{\text{ev},+} [  d ]  \oplus\text{H}_{\rm Hod}%
^{\text{ev},-} [  d ]  ,\quad \text{H}_{\rm Hod}^{\text{odd}} (  X )
[d]=\text{H}_{\rm Hod}^{\text{\rm odd},+}[d]
\oplus\text{H}_{\rm Hod}^{\text{\rm odd},-}[d] .
\]
In order to obtain detailed information about these decompositions, we shall
need to compare $ (  \text{H}_{\rm Hod}^{\text{ev}} (  X )[d] , 
\sigma^{\text{ev}} )  $ and $ (  \text{H}_{\rm Hod}%
^{\text{odd}} (  X )  [d]  ,\sigma^{\text{odd}} )$ 
with the de Rham hypercohomology ${\rm H}_{\rm dR}^{\bullet} (  X )[d]={\rm H}^{\bullet} 
(  X,\Omega_{X/{\bf R}}^{\bullet}[d]
 )  $ of $\Omega_{X/{\bf R}}^{\bullet} [  d ]  $ endowed with  the
$G$-invariant forms from duality theory: recall that duality for de Rham
hypercohomology gives a perfect $\mathbf{R}$-bilinear form
\[
t^{p}:\text{H}_{\rm dR}^{p}(X)[d]  \times\text{H}%
_{\rm dR}^{-p} (X)[d]  \rightarrow\text{H}_{\rm dR}%
^{0}(X)[2d]  ={\rm H}_{\rm dR}^{2d}(X)
\rightarrow\mathbf{R.}%
\]
where by the Wirtinger theorem (see page 31 in [GH]) the right hand map is
given by real integration $\omega\mapsto\left|  G\right|  ^{-1}\int_{X}\omega$
for a global real $2d$-form $\omega$. As $d$ is even, the map $t^{0}$ is
symmetric (see below); note also that if $x\in{\rm H}_{\rm dR}^{p}(X)
[d]$, $y\in {\rm H}_{\rm dR}^{-p}(X)[d]$,
then $t^{p}\left(  x,y\right)  =\left(  -1\right)  ^{p}t^{-p}\left(
y,x\right)$ which again of course agrees with the commutation rule (15);
hence, as per the construction in 4.1, we may then form the symmetrised
duality maps $\tau_{\text{ }}^{p}$.

We write $\Omega_{X}^{\bullet<m}$ respectively $\Omega_{X}^{\bullet\geq m}$
for the complex
\[
O_{X}\overset{d}{\rightarrow}\Omega_{X/\mathbf{R}}^{1}\overset{d}{\rightarrow
} \cdots\overset{d}{\rightarrow}\Omega_{X/\mathbf{R}}^{m-1}%
\]%
\[
\Omega_{X/\mathbf{R}}^{m}\overset{d}{\rightarrow}\cdots\overset
{d}{\rightarrow}\Omega_{X/\mathbf{R}}^{d}%
\]
where the term $O_{X}$ (respectively  $\Omega_{X/\mathbf{R}}^{m}$) is placed
in degree zero (respectively degree $m$).   We then consider the exact
sequence of complexes
\[
0\rightarrow\Omega_{X}^{\bullet\geq m} [  d ]
\rightarrow\Omega_{X/\mathbf{R}}^{\bullet} [  d ]  \rightarrow
\Omega_{X}^{\bullet<m} [  d ]  \rightarrow0
\]
and we let $F^{-m+d/2}$ denote the image of ${\rm H}_{\rm dR}^{\bullet} (
X,\Omega_{X}^{\bullet\geq m} [  d ]   )  $ in ${\rm H}^{\bullet
} (  X,\Omega_{X/{\bf R}}^{\bullet} [  d])  ={\rm H}_{\rm dR}^{\bullet}(X)[d]   $ under the natural
map induced by $\Omega_{X}^{\bullet\geq m}\hookrightarrow
\Omega_{X/\mathbf{R}}^{\bullet}$. Note that by the degeneration of the Hodge
spectral sequence we know that in fact ${\rm H}^{\bullet} (  X,\Omega
_{X }^{\bullet\geq m}[d])  $ injects
into ${\rm H}^{\bullet} (  X,\Omega_{X/\mathbf{R}}^{\bullet} [  d ])  $.

\begin{theorem}
The quadratic space $ (  \text{\rm H}_{\rm dR}^{\bullet}(X)[d],\ \tau )$, when endowed with 
the filtration $ \{
F^{i} \}  _{i}$,   is a filtered quadratic space, as defined in Lemma
27. There is an isomorphism of $\mathbf{R} [  G ]  $-quadratic
modules:
\[
 (  \text{\rm H}_{dR}^{\bullet}(X)[d] 
,\tau )  \cong (  \text{\rm H}^{d/2} (X,  \Omega_{X/{\bf R}}^{d/2} )
,\sigma_{d/2,d/2} )  \oplus\text{\rm Hyp} (  \oplus_{i<d/2}\text{\rm H}%
^{i} (X,  \Omega_{X/{\bf R}}^{d/2} ) )  \oplus\text{\rm Hyp} (
\text{\rm H}^{\bullet} (X,  \Omega_{X/{\bf R}}^{\bullet>d/2} )   )  .
\]
\end{theorem}

\noindent{\sl  Proof.} Consider the shifted symmetrised duality maps

\begin{equation*}
\mathrm{H}_{\mathrm{dR}}^{-p+d}\left( X\right) \times \mathrm{H}_{\mathrm{dR}%
}^{p+d}\left( X\right) \rightarrow \mathrm{H}_{\mathrm{dR}}^{2d}\left(
X\right) \rightarrow \mathbf{R.}
\end{equation*}
It then follows immediately that the pairing 
\begin{equation*}
\mathrm{H}^{\bullet } ( X,\Omega _{X}^{\bullet \geq m} ) \times 
\mathrm{H}^{\bullet } ( X,\Omega _{X}^{\bullet \geq d-m} )
\rightarrow \mathbf{R}
\end{equation*}
factors through 
\begin{equation*}
\tau _{m}^{{\prime }}:\frac{\mathrm{H}^{\bullet } ( X,\Omega
_{X}^{\bullet \geq m} ) }{\mathrm{H}^{\bullet } ( X,\Omega
_{X}^{\bullet \geq m+1} ) }\times \frac{\mathrm{H}^{\bullet } (
X,\Omega _{X}^{\bullet \geq d-m} ) }{\mathrm{H}^{\bullet } (
X,\Omega _{X}^{\bullet \geq d-m+1} ) }\rightarrow \mathbf{R.}
\end{equation*}
To see that this pairing is perfect, we again appeal to the decomposition of
the Hodge spectral sequence to deduce that for all $m$, $n$
\begin{equation*}
\frac{\text{\textrm{H}}^{n} ( X,\Omega _{X}^{\bullet \geq m} ) }{%
\text{\textrm{H}}^{n} ( X,\Omega _{X}^{\bullet \geq m+1} ) }\cong 
\text{\textrm{H}}^{n} ( X,\Omega _{X}^{m} ) .
\end{equation*}
Thus $\tau _{m}^{{\prime }}$ induces forms 
\begin{equation*}
\tau _{m}^{n}:\text{\textrm{H}}^{n} ( X,\Omega _{X}^{m} ) \times 
\text{\textrm{H}}^{d-n} ( X,\Omega _{X}^{d-m} ) \rightarrow \mathbf{R}.
\end{equation*}
We now claim that under the above isomorphisms the forms $\tau _{m}^{n}$
agree with the pairings $\sigma _{n,m}$ up to the sign $ ( -1 )
^{ ( m+n ) n}$, and in particular agree exactly when $m=d/2=n$. To
show this it is enough to apply $\otimes _{\mathbf{R}}\mathbf{C}$; identify
the left-hand terms with the Dolbeault cohomology groups \textrm{H}$_{%
\overline{\partial }}^{m,n} ( X ) $ and show

\begin{proposition}
The following diagram commutes: 
\begin{equation}
\begin{array}{ccccc}
 {\rm H}^{q} ( X_{\mathbf{C}},\Omega _{X_{\mathbf{C}%
}}^{p} )  & \times  & {\rm H}^{d-q} ( X_{\mathbf{C}%
},\Omega _{X_{\mathbf{C}}}^{d-p} )  & \overset{\cup }{\rightarrow } & 
{\rm H}^{d} ( X_{\mathbf{C}},\Omega _{X_{\mathbf{C}%
}}^{d} )  \\ 
\downarrow  &  & \downarrow  &  & \downarrow  ( -1 ) ^{ (
p+q ) q} \\ 
{\rm H}_{\overline{\partial }}^{p,q} ( X_{\mathbf{C}} ) 
& \times  & {\rm H}_{\overline{\partial }}^{d-p,d-q} ( X_{%
\mathbf{C}} )  & \overset{\wedge }{\rightarrow } & {\rm H}_{%
\overline{\partial }}^{d,d} ( X_{\mathbf{C}} ) 
\end{array}
\end{equation}
where the two left-hand vertical maps are Dolbeault isomorphisms and
right-hand vertical map is the Dolbeault isomorphism multiplied by $ (
-1 ) ^{ ( p+q ) q}$. In particular if $p=q$ then  $ (
-1 ) ^{ ( p+q ) q}=1$ and so in this case $\cup $ and $\wedge $
agree under the Dolbeault isomorphism.
\end{proposition}

\bigskip Before proving the proposition, we first note that it will complete
the proof of the theorem. Indeed, by Lemma 27, we know that 
\begin{eqnarray*}
( \text{\textrm{H}}_{dR}^{\bullet } ( X )  [ d ]
,\tau  )  &\cong & ( \oplus _{i}\text{\textrm{H}}^{i} ( \Omega
_{X}^{d/2} ) ,\tau _{d/2}^{\prime } ) \oplus \text{\textrm{Hyp}}%
 ( \text{\textrm{H}}^{\bullet }( \Omega _{X}^{\bullet >d/2} )
 )  \\
&\cong & ( \text{\textrm{H}}^{d/2} ( \Omega _{X}^{d/2} ) ,\tau
_{d/2}^{d/2} ) \oplus \;\text{\textrm{Hyp}} ( \oplus _{i<d/2}\text{%
\textrm{H}}^{i} ( \Omega _{X}^{d/2} )) \oplus \text{\textrm{%
Hyp}} ( \text{\textrm{H}}^{\bullet } ( \Omega _{X}^{\bullet
>d/2})) ;
\end{eqnarray*}
however, by the above discussion together with the proposition, we know that 
$\tau _{d/2}^{d/2}$ is equal to $\sigma _{d/2,d/2}$ and the result will now
follow.

\bigskip

\textbf{Proof of Proposition 38.} In unraveling the Dolbeault isomorphisms
we shall follow the conventions given in\ Section 3 of Chapter 0 in [GH]. We
first need some notation: let $\frak{A}^{p,q}$ denote the sheaf of $%
C^{\infty }$-forms on $X_{\mathbf{C}}$ of type $ ( p,q ) $ we let $%
\frak{Z}^{p,q} $ denote the sheaf of $\overline{\partial }$%
-closed $C^{\infty}$ forms of type $(p,q)$ and we write $\Omega ^{p}$ 
for the sheaf of holomorphic $p$-forms $%
\frak{Z}^{p,0}$. By the $\overline{\partial }$-Poincar\'{e} Lemma we have
exact sequences of sheaves 
\begin{equation}
0\rightarrow \frak{Z}^{p,q}\rightarrow \frak{A}^{p,q}\overset{\overline{%
\partial }}{\rightarrow }\frak{Z}^{p,q+1}\rightarrow 0.
\end{equation}
For each $n,\;0\leq n<q,$ we then consider the exact sequence 
\begin{equation}
0\rightarrow \frak{Z}^{p,q}\otimes \frak{Z}^{d-p,d-q-n-1}\rightarrow \frak{Z}%
^{p,q}\otimes \frak{A}^{d-p,d-q-n-1}\overset{1\otimes \overline{\partial }}{%
\rightarrow }\frak{Z}^{p,q}\otimes \frak{Z}^{d-p,d-q-n}\rightarrow 0.
\end{equation}
For brevity we write ${\rm H}^{a} ( \frak{Z}^{p,q} ) $ for 
${\rm H}^{a} ( X_{\mathbf{C}},\frak{Z}^{p,q} )$ etc. and we
let $\delta _{1}^{n}$ resp. $\delta _{2}^{n}$ denote the $n$-th coboundary
map associated to the cohomology of the exact sequence (23) resp. (24). Then
for  $x\in {\rm H}^{0} ( \frak{Z}^{p,q} )$, 
$y\in {\rm H}^{n} ( \frak{Z}^{d-p,d-q-n} )$, by using the 
cocycle description of the boundary maps, we obtain $x\wedge
\delta _{1}^{n} ( y ) =\delta _{2}^{n} ( x\wedge y )$.
Hence for each such $n$, $0\leq n<d-q$, we get a commutative diagram 
\begin{equation*}
\begin{array}{ccccc}
\text{\textrm{H}}^{0} ( \frak{Z}^{p,q} )  & \times  & \text{\textrm{%
H}}^{n} ( \frak{Z}^{d-p,d-q-n} )  & \overset{\wedge }{\rightarrow } & 
\text{\textrm{H}}^{n} ( \frak{Z}^{p,q}\frak{\wedge Z}^{d-p,d-q-n} ) 
\\ 
\downarrow 1 &  & \downarrow \delta _{1}^{n} &  & \downarrow \delta _{2}^{n}
\\ 
\text{\textrm{H}}^{0} ( \frak{Z}^{p,q} )  & \times  & \text{\textrm{%
H}}^{n+1} ( \frak{Z}^{d-p,d-q-n-1} )  & \overset{\wedge }{%
\rightarrow } & \text{\textrm{H}}^{n+1} ( \frak{Z}^{p,q}\frak{\wedge Z}%
^{d-p,d-q-n-1} ) 
\end{array}
\end{equation*}
and we claim that $\delta _{2}^{n}$ agrees with $\left( -1\right) ^{p+q}$
times the Dolbeault map $D^{n}:{\rm H}^{n} ( \frak{Z}%
^{d,d-n} ) \rightarrow {\rm H}^{n+1} ( \frak{Z}%
^{d,d-n-1} )$. Recall that $ D^{n}$ is an isomorphism for positive $n
$ and that $D^{0 }$ is a surjection which induces the isomorphism 
\begin{equation*}
\text{\textrm{H}}_{\overline{\partial }}^{p,q} ( X_{\mathbf{C}} ) :=%
\text{\textrm{H}}^{0} ( \frak{Z}^{p,q} ) /\overline{\partial }%
 ( \text{\textrm{H}}^{0} ( \frak{A}^{p,q-1} )  )
\rightarrow \text{\textrm{H}}^{1} ( \frak{Z}^{p,q-1} ) .
\end{equation*}
In proving the claim it will be useful to have the following cocycle formula
for the Dolbeault maps $D^{n}: {\rm H}^{n} ( \frak{Z}^{a,b} )
\rightarrow {\rm H}^{n+1} ( \frak{Z}^{a,b-1} ) $. Once and
for all we fix a sufficiently fine cover $\mathcal{U}=\left\{ U_{i}\right\}
_{i\in I}$ of $X({\mathbf{C}})$. For a $\frak{Z}^{a,b}$-valued $n$-cocycle $%
\omega ^{n}$ we write $\omega \left( i_{0},\ldots ,i_{n}\right) $ for $\omega
\left( U_{i_{0}}\cap \cdots\cap U_{i_{n}}\right) $. Then 
\begin{equation*}
D^{n}\left( \omega ^{n}\right) \left( i_{0},\ldots,i_{n+1}\right)
=
\sum_{k=0}^{n+1} ( -1 ) ^{k}\overline{\partial }^{-1} (
\omega ^{n} ( \widehat{i}_{k} )  ) 
\end{equation*}
where $\omega ^{n} ( \widehat{i}_{k} ) $ means $\omega ^{n} (
i_{0},\ldots ,\widehat{i}_{k},\ldots,i_{n+1} )$. Now let $\gamma \in {\rm H}
^{0}( \frak{Z}^{p,q} )$ , $\omega ^{n}\in {\rm H}^{n} ( 
\frak{Z}^{d-p,d-q-n} )$.  Then 
\begin{equation*}
D^{n}\left( \gamma \wedge \omega ^{n}\right) \left( i_{0},\ldots, i_{n+1}\right)
=
\sum_{k=0}^{n+1} ( -1 ) ^{k}\overline{\partial }^{-1} (
\gamma \wedge \omega ^{n} ( \widehat{i}_{k} ) ) 
\end{equation*}
while 
\begin{equation*}
\delta _{2}^{n} ( \gamma \wedge \omega ^{n} )  (
i_{0},\ldots , i_{n+1} ) =
\sum_{k=0}^{n+1} ( -1 ) ^{k} ( (
1\otimes \overline{\partial } ) ^{-1}\gamma \wedge \omega ^{n} ( 
\widehat{i}_{k})) .
\end{equation*}
The claim then follows by the Leibniz rule for $\overline{\partial }$ and
using the fact the $\overline{\partial } ( \gamma  ( i_{0} )
 ) =0$: indeed, for $k>0$ 
\begin{eqnarray*}
\overline{\partial }^{-1} ( \gamma  ( i_{0} ) \wedge \omega
^{n} ( \widehat{i}_{k} )  )  &=& ( -1 ) ^{p+q}\gamma
 ( i_{0} ) \wedge \overline{\partial }^{-1}\omega ^{n} ( 
\widehat{i}_{k} )  \\
&=& ( -1 ) ^{p+q} ( 1\otimes \overline{\partial } )
^{-1} ( \gamma  ( i_{0} ) \wedge \omega ^{n} ( \widehat{i}%
_{k} )  ) 
\end{eqnarray*}
and similarly when $k=0$ 
\begin{equation*}
\overline{\partial }^{-1} ( \gamma  ( i_{1} ) \wedge \omega
^{n} ( \widehat{i}_{0} )  ) = ( -1 ) ^{p+q} (
1\otimes \overline{\partial } ) ^{-1} ( \gamma  ( i_{1} )
\wedge  \omega ^{n} ( \widehat{i}_{0} )) .
\end{equation*}

Next, for $m$, $0\leq m<q,$ we consider the exact sequences 
\begin{equation}
0\rightarrow \frak{Z}^{p,m-1}\otimes \Omega ^{d-p}\rightarrow \frak{A}%
^{p,m-1}\otimes \Omega ^{d-p}\overset{\overline{\partial }=\overline{%
\partial }\otimes 1}{\rightarrow }\frak{Z}^{p,m}\otimes \Omega
^{d-p}\rightarrow 0
\end{equation}
and we write $\delta _{3}^{m}$ for the $m$-th coboundary map associated to
the cohomology of this exact sequence. For $x\in {\rm H}^{m} ( 
\frak{Z}^{p,m} )$ , $y\in  {\rm H}^{q} (  {\Omega }%
^{d-p} )$, we know that $\delta _{1}^{m} ( x ) \wedge y=\delta
_{3}^{m} ( x\wedge y )$ ([Go] Ch. 6.5, p. 255). Therefore for each such $m$ we obtain a
commutative square 
\begin{equation*}
\begin{array}{ccccc}
\text{\textrm{H}}^{m} ( \frak{Z}^{p,m} )  & \times  & \text{\textrm{%
H}}^{d-q} (  \Omega ^{d-p} )  & \overset{\wedge }{\rightarrow }
& \text{\textrm{H}}^{d-q+m} ( \frak{Z}^{p,m}\wedge  \Omega %
^{d-p} )  \\ 
\downarrow \delta _{1}^{m} &  & \downarrow 1 &  & \downarrow \delta _{3}^{m}
\\ 
\text{\textrm{H}}^{m+1} ( \frak{Z}^{p,m-1})  & \times  & \text{%
\textrm{H}}^{d-q} (  \Omega^{d-p} )  & \overset{\wedge }{%
\rightarrow } & \text{\textrm{H}}^{d-q+m+1} ( \frak{Z}^{p,m-1}\wedge 
{\Omega }^{d-p} ) 
\end{array}
\end{equation*}
and we now claim that $\delta _{3}^{m}$ agrees with the Dolbeault map $D^{d-q+m}$. 
For the sake of brevity we put $p^{\prime }=d-p$, $q^{\prime
}=d-q$ and we let $\omega ^{m}\in {\rm H}^{m} ( \frak{Z}%
^{p,m} )$, $\nu ^{q^{\prime }}\in {\rm H}^{q^{\prime }} ( 
\Omega^{p^{\prime }} )$. Then 
\begin{equation*}
D^{m+q^{\prime }} ( \omega ^{m}\wedge \nu ^{q^{\prime }} )  (
i_{0},\ldots, i_{m+q^{\prime }+1} ) =
\sum_{k=0}^{m+q^{\prime }+1} (
-1 ) ^{k}\overline{\partial }^{-1} ( \omega ^{m}\wedge \nu
^{q^{\prime }} )  ( \widehat{i}_{k} ) 
\end{equation*}
while 
\begin{equation*}
\delta _{3}^{m} ( \omega ^{m}\wedge \nu ^{q^{\prime }} )  (
i_{0},\ldots, i_{m+q^{\prime }+1} ) =
\sum_{k=0}^{m+q^{\prime }+1} (
-1 ) ^{k} ( \overline{\partial }\otimes 1 ) ^{-1} ( \omega
^{m}\wedge \nu ^{q^{\prime }} )  ( \widehat{i}_{k} ) .
\end{equation*}
Since $\overline{\partial }\otimes 1=\overline{\partial }$ 
as we have seen in (25) we obtain that $\delta_m^3$ agrees with $D^{d-q+m}$.  
Recall now that the Dolbeault isomorphisms ([GH] p. 45) are obtained as 
a composition of a succession of  Dolbeault maps (there are $q$ of these maps)
$$
{\rm H}_{\overline{\partial }}^{p,q} ( X_{\mathbf{C}} )={\rm H}^0({\mathfrak Z}^{p,q})/\bar\partial 
({\rm H}^0 ( {\mathfrak A}^{p,q-1} ))\xrightarrow{D} 
\cdots\xrightarrow{D} {\rm H}^{q-1}({\mathfrak Z}^{p,1})
\xrightarrow{D}{\rm H}^q(\Omega^p).
$$
By combining the above results with an inductive argument, we can now 
see that the sign discrepancy between the two pairings is equal to 
the product of $q$ copies of $(-1)^{p+q}$, therefore equal to $(-1)^{(p+q)q}$.
 \ \ \ \ $ \square $

\bigskip

\ \

\begin{corollary}
\bigskip There is a non-canonical $\mathbf{R}[G] $-isometry
\[
\left(  \text{\rm H}_{\rm dR}^{\bullet}(X)[d] 
,\tau\right)  \otimes_{\mathbf{Q}}\mathbf{R}\cong \left(  \text{\rm H}%
_{\rm Hod}^{\bullet}(X)[d]   ,\sigma \right)
\otimes_{\mathbf{Q}}\mathbf{R}.
\]
\end{corollary}

\noindent{\sl Proof.} By the above proposition we know that each of the above quadratic
spaces is isometric to the orthogonal sum of $ (  \text{H}^{d/2} (X,
\Omega_{X/{\bf R}}^{d/2} )  ,\sigma_{d/2,d/2} )  $ and a hyperbolic space.
On the other hand by the degeneration of the Hodge spectral sequence we know
that ${\rm H}_{dR}^{\bullet} (X)  \otimes_{\mathbf{R}}\mathbf{C}$ and
${\rm H}_{\rm Hod}^{\bullet} (X)  \otimes_{\mathbf{R}}\mathbf{C}$ are
isomorphic $\mathbf{C}[G]  $-modules; hence ${\rm H}_{\rm dR}^{\bullet
}\left(  X\right)  $ and ${\rm H}_{\rm Hod}^{\bullet}\left(  X\right)$ are isomorphic
$\mathbf{R} [  G ]  $-modules; therefore we may conclude that the two
hyperbolic spaces are isometric, as required. \hfill$\square $
\medskip

Next we recall that from Proposition 1.4 on p 319 in [D]:

\begin{proposition}
The   comparison isomorphism  
${\rm H}_{\rm dR}^{\bullet} (  X )
\otimes_{\mathbf{Q}}\mathbf{C}\cong {\rm H}_{B}^{\bullet}(X)
\otimes_{\mathbf{R}}\mathbf{C}$ identifies ${\rm H}_{\rm dR}^{\bullet}(X)
\otimes_{\mathbf{Q}}\mathbf{R}$  with ${\rm H}_{B+}^{\bullet}%
\oplus (  \text{\rm H}_{B-}^{\bullet}\otimes_{\mathbf{R}}i\mathbf{R} )$.
\end{proposition}

\bigskip Writing
\[
\text{H}_{B\pm}^{\text{ev},\pm}=\text{H}_{B\pm}^{\text{ev}}\cap\text{H}%
_{B}^{\text{ev},\pm},\quad \text{ \ H}_{\rm dR}^{\text{ev}}=\text{H}_{\rm dR}%
^{\text{ev}} (  X)  , 
\]
and similarly with odd in place of even, we conclude from the above proposition that:

\begin{corollary}%
\[
\dim(\text{\rm H}_{\rm dR}^{\text{\rm ev},+})=\dim(\text{\rm H}_{B+}^{\text{\rm ev},+})+\dim
(\text{\rm H}_{B-}^{\text{\rm ev},-}),\ \ \dim(\text{\rm H}_{\rm dR}^{\text{\rm ev},-})=\dim
({\rm H}_{B+}^{\text{\rm ev},-})+\dim({\rm  H}_{B-}^{\text{\rm ev},+}),%
\]%
\[
\dim(\text{\rm H}_{\rm dR}^{\text{\rm odd},+})=\dim(\text{\rm H}_{B+}^{\text{\rm odd},+})+\dim
(\text{\rm H}_{B-}^{\text{\rm odd},-}),\ \ \dim(\text{\rm H}_{\rm dR}^{\text{\rm odd},-})%
=\dim(\text{\rm H}_{B+}^{\text{\rm odd},-})+\dim(\text{\rm H}_{B-}^{\text{\rm odd},+}).
\]
Moreover by Corollary 39 the same statements hold with ${\rm H}_{\rm Hod}^{{\rm ev},+}$ in
place of ${\rm H}_{\rm dR}^{{\rm ev},+}$ etc. on the left-hand side.
\end{corollary}

\subsection{Proof of Theorem 1.}

Before embarking on the proof of the theorem, we first need:

\begin{proposition}
${\rm H}_{\rm Hod}^{\text{\rm ev}}(X)  -{\rm H}_{\rm Hod}^{\text{\rm odd}}(X)$
is a (virtually) free $\mathbf{Q}[G]$-module.
\end{proposition}

\noindent{\sl Proof.} This follows   from Theorem 36 and the comparison
isomorphism ${\rm H}_{\rm dR}^{\bullet}(X)  \otimes_{\mathbf{Q}%
}\C\cong {\rm H}_{B}^{\bullet}(X) \otimes_{\mathbf{R}%
}\mathbf{C}$. \hfill $\square $
\medskip

Throughout this sub-section we assume $W$ to be a virtual complex symplectic
representation of $G$. Recall from (21) that in order to prove the theorem it
is sufficient to show
\begin{equation}
i^{\delta(Y)  \dim(W) }\varepsilon_{\infty} (
\mathcal{X},W )  =i^{n_{W}^{-} (  \sigma )  }.
\end{equation}
Initially we shall suppose that $\dim (  W )  =0$; then, by linearity
in both sides, we shall conclude the proof by dealing with the case where
$W$ is two copies of the trivial representation.

By the above and Lemma 35, in order to prove Theorem 1 when $\dim (
W )  =0$, we are required to establish the congruence modulo 4
\[
n_{W}^{-}\left(  \sigma\right)  \equiv\left\{
\begin{array}
[c]{c}%
0,\text{ \ \ if \ }d\text{\ is odd,}\\
\chi_{-}\left(  X\otimes_{G}W\right),  \;\;\;\;\text{if \ }d\text{\ is even.}%
\end{array}
\right.
\]

\medskip\qquad\qquad\ \ \qquad\qquad

Case 1. $d$\ \textit{is odd.} On the one hand we know that $\sigma$ is
hyperbolic; on the other hand by the above proposition we know that
${\rm H}_{\rm Hod}^{\bullet}(X)  _{\mathbf{R}}$ is a free virtual
$\mathbf{R} [  G ]  $-module; thus by Lemma 26 we know that
${\rm H}_{\rm Hod}^{\text{ev}-}(X)_{\mathbf{R}}-{\rm H}_{\rm Hod}^{\text{odd}%
-}(X)_{\mathbf{R}}$ is also a free $\mathbf{R} [  G ]
$-module; therefore, because $W$ has dimension zero, it follows at once that
$n_{W}^{-} (  \sigma )  =0$, as required.\medskip

\ Case 2. $d$\ \textit{is even.} Since $d$ is even we may freely use the
notation of 4.3. Then by Corollary 41
\begin{equation}
n_{W}^{-} (  \sigma )  =\dim(  \text{H}_{B+}^{\bullet-}\otimes
W)  ^{G}+\dim(  \text{H}_{B-}^{\bullet+}\otimes W)  ^{G}.
\end{equation}
By Theorem 36 ${\rm H}_{B}^{\bullet-}={\rm H}_{B+}^{\bullet-}+{\rm H}_{B-}^{\bullet-}$
is a free $G$-module and again since $\dim (  W )  =0,\;$we see that
\[
0=\dim (  \text{H}_{B}^{\bullet-}\otimes W )  ^{G}=\dim (
\text{H}_{B+}^{\bullet-}\otimes W )  ^{G}+\dim (  \text{H}%
_{B-}^{\bullet-}\otimes W )  ^{G}%
\]
and\ therefore by (27)
\[
n_{W}^{-} (  \sigma )  =-\dim (  \text{H}_{B-}^{\bullet-}\otimes
W )  ^{G}+\dim (  \text{H}_{B-}^{\bullet+}\otimes W )  ^{G}.
\]
By Proposition 25  we know that each ${\rm H}_{B}^{\bullet\pm}(  X (
\mathbf{C} )  ,\mathbf{R} )  $ is an $\mathbf{R} [  G ]
$-module, and so, reasoning as after the proof of Lemma 15, we see that
$\dim (  \text{H}_{B-}^{\bullet-}\otimes W )  ^{G}$ is even, and it
therefore follows that we have the congruence mod 4
\[
n_{W}^{-} (  \sigma )  \equiv\dim (  \text{H}_{B-}^{\bullet
-}\otimes W )  ^{G}+\dim (  \text{H}_{B-}^{\bullet+}\otimes W )
^{G}=\chi_{-} (  X\otimes_{G}W )  \;\;
\]
as required. 

To conclude we now suppose that $W $ is two copies of the trivial
representation; so that ${\rm H}_{B}^{\bullet}(X)\otimes_{G}W$ is now just two copies
of ${\rm H}_{B}^{\bullet}(Y)$. We are now required to establish the
congruence $\operatorname{mod}4$
\begin{equation}
n_{W}^{-}\left(  \sigma\right)  \equiv\left\{
\begin{array}
[c]{c}%
  \chi\left(  Y\right), \ \ \   \text{
\ \ \ \ \ \ \ \ \ \ \ \  \ \ \ \ \ \ \   if
\ }d\text{\ is odd.}\ \ \ \ \ \ \ \ \ \ \ \ \ \ \\
2\chi^{\pm} (  Y )  +2\chi_{\pm}(  Y), \  \text{if
 $d$  is even and $\pm$  denotes the sign of }\left(  -1\right)
^{d/2+1}.
\end{array}
\right.
\end{equation}

Suppose first that $d$ is odd. Since $\sigma$ is hyperbolic, by Lemma 26
$n_{W}^{-}\left(  \sigma\right)  =\chi\left(  Y\right)  $, as required.

Suppose next that $d$ is even, so that again we may use the notation of 4.3.
Recall that $ (  \text{H}_{B}^{\bullet\pm}(X)\otimes W )  ^{G}$ is two
copies of ${\rm H}_{B}^{\bullet\pm}(Y)  ,$ and, as previously, put
$\chi^{\pm}\left(  Y\right)  =\dim (  \text{H}_{B}^{\bullet\pm} (
Y )   ) $ etc.; more generally, we shall write $\chi_{\pm}^{\pm
}(  Y)$ for $\dim(  \text{H}_{B,\pm}^{\text{ev},\pm
}(  Y)  )  -\dim(  \text{H}_{B,\pm}^{\text{odd},\pm
}\left(  Y\right)  )  $. Observe that by Corollary 41 we again have
congruences $\operatorname{mod}4$%
\[
n_{W}^{-}\left(  \sigma\right)  =2\chi_{+}^{-}\left(  Y\right)  +2\chi_{-}%
^{+}\left(  Y\right)
\]%
\[
=2\chi^{-}\left(  Y\right)  -2\chi_{-}^{-}\left(  Y\right)  +2\chi_{-}%
^{+}\left(  Y\right)
\]%
\[
\equiv2\chi^{-}\left(  Y\right)  -2\chi_{-}^{-}\left(  Y\right)  -2\chi
_{-}^{+}\left(  Y\right)
\]%
\[
\equiv2\chi^{-}\left(  Y\right)  -2\chi_{-}\left(  Y\right)  .
\]

Case 1. $d\equiv2\operatorname{mod}4.\;\;$ In this case by (28) we have to
show the congruence $\operatorname{mod}4$%
\[
2\chi^{-}\left(  Y\right)  -2\chi_{-}\left(  Y\right)  \equiv2\chi^{+}\left(
Y\right)  +2\chi_{+}\left(  Y\right)
\]
which is clear since
\[
\chi_{+}\left(  Y\right)  +\chi_{-}\left(  Y\right)  =\chi\left(  Y\right)
=\chi^{+}\left(  Y\right)  +\chi^{-}\left(  Y\right)  .
\]

Case 2. $d\equiv0\operatorname{mod}4.\;$This follows at once since we have to
show the congruence $\operatorname{mod}4$
\[
2\chi^{-}\left(  Y\right)  +2\chi_{-}\left(  Y\right)  \equiv2\chi^{-}\left(
Y\right)  -2\chi_{-}\left(  Y\right)
\]
which is immediate.\ \ \ \ \ $\square$

\bigskip

\section{Appendix: \  Comparison of definitions.}

The symplectic hermitian class group that we have used, namely ${\rm H}^{\text{s}%
} (  \mathbf{Z} [  G ] )  $, is very well suited to
comparison with Arakelov invariants; indeed, from (5) we see that it is the
natural vehicle for carrying discriminantal signs associated to Arakelov
discriminants. In this Appendix we briefly indicate how the class group
${\rm H}^{\text{s}} (  \mathbf{Z} [  G ]   )$, and hermitian
classes formed in this group, relate to the previous hermitian classes and
classgroups, such as those used for instance in [F1] and [CPT1].

\subsection{Hermitian class groups.}

\bigskip Recall that in Definition 11 we defined the symplectic hermitian
class group H$^{\text{s}}\left(  \mathbf{Z}\left[  G\right]  \right)  $ to be
\begin{equation}
\text{H}^{\text{s}} (  \mathbf{Z} [  G ]   )  =\frac
{\text{Hom}_{\Omega_{\mathbf{Q}}} (  R_{G}^{\text{s}},J_{f} )
\times\text{Hom} (  R_{G}^{\text{s}},\mathbf{R}^{\times} )
}{\operatorname{\rm Im}(\Delta^{\text{s}})\cdot (  \text{Det}^{\text{s}} (
\widehat{\mathbf{Z}}\left[  G\right]  ^{\times} )  \times1 )  }.
\end{equation}
By contrast in [F1] and [CPT1] the hermitian class group ${\rm HCl} (
\mathbf{Z} [  G])  $ is used, which is described in terms of
character functions as
\begin{equation}
\text{HCl} (  \mathbf{Z} [  G ]  )  =\frac{\text{Hom}%
_{\Omega_{\mathbf{Q}}} (  R_{G},J_{f} )  \times\text{Det}(
\mathbf{R} [  G ]  ^{\times} )  \times\text{Hom}_{\Omega
_{\mathbb{Q}}} (  R_{G}^{\text{s}},\overline{\mathbf{Q}}^{\times} )
}{\operatorname{\rm Im}(\widetilde{\Delta})\cdot(  \text{Det} (  \widehat
{\mathbf{Z}} [  G ]  ^{\times}\times\mathbf{R} [  G ]
^{\times} )  \times1 )  }%
\end{equation}
where $\widetilde{\Delta}$ is the twisted diagonal map
\[
\widetilde{\Delta}:\text{Det} (  \mathbf{Q}[G] ^{\times
} )  \rightarrow\text{Hom}_{\Omega_{\mathbb{Q}}} (  R_{G}%
,J_{f} )  \times\text{Det} (  \mathbf{R}[G]  ^{\times
} )  \times\text{Hom}_{\Omega_{\mathbf{Q}}} (  R_{G}^{\text{s}%
},\overline{\mathbf{Q}}^{\times} )
\]
given by $\widetilde{\Delta} (  \text{Det} (  a )   )
={\rm Det} (  a )  \times{\rm Det} (  a )  \times{\rm Det}^{\text{s}%
} (  a )  ^{-1}$. We therefore have a natural map from
\begin{equation}
\phi:\text{HCl} (  \mathbf{Z}[G])  \rightarrow
\text{H}^{\text{s}} (  \mathbf{Z}[G])
\end{equation}
induced by the map
\[
\text{Hom}_{\Omega_{\mathbf{Q}}} (  R_{G},J_{f} )  \times
\text{Det} (  \mathbf{R}[G] ^{\times} )  \times
\text{Hom}_{\Omega_{\mathbf{Q}}}(  R_{G}^{\text{s}},\overline{\mathbf{Q}%
}^{\times} )  \rightarrow\text{Hom}_{\Omega_{\mathbb{Q}}} (
R_{G}^{\text{s}},J_{f} )  \times\text{Hom} (  R_{G}^{\text{s}%
},\mathbf{R}^{\times} )
\]
which: takes the first left-hand factor into the first right-hand factor by
restriction from $R_{G}$ to $R_{G}^{\text{s}}$; which is trivial on the second
left-hand factor; and which maps the third left-hand factor to the second
right-hand factor by inverting the natural map induced by the inclusion
$\overline{\Q} \hookrightarrow \C$.

\subsection{\bigskip Hermitian classes.}

Next we recall the construction of \ hermitian Euler characteristics used in
[F1] and [CPT1]; we compare this definition with that given in 3.2.1, and then
see how they match under the comparison map $\phi$ above.

So here we consider a perfect $\mathbf{Z} [  G ]  $-complex
$P^{\bullet}$
\[
\cdots\rightarrow P^{i}\overset{\partial^{i}}{\rightarrow}P^{i+1}%
\rightarrow\cdots
\]
which supports non-degenerate $G$-invariant $\mathbf{Q}$-valued
forms
\[
\sigma^{i}:\text{H}^{i} (  P^{\bullet} )  _{\mathbf{Q}}\times
\text{H}^{-i} (  P^{\bullet} )  _{\mathbf{Q}}\rightarrow\mathbf{Q}%
\]
which are symmetric in the sense that $\sigma^{i}\left(  x,y\right)
=\sigma^{-i}\left(  y,x\right)  $. In Proposition 2.7 of [CPT1] we show that,
after adding an acyclic complex  to $P^{\bullet}$ if necessary, each
$P_{\mathbf{Q}}^{i}=\mathbf{Q}\otimes P^{i}$ admits a $G$-decomposition
$P_{\mathbf{Q}}^{i}=B^{i}\oplus H^{i}\oplus U^{i}$ with $B^{i}%
= {\rm Im}(\partial^{i-1})$ and with $U^{i}$ mapped isomorphically onto
$B^{i+1}$ by $\partial^{i},$ and there exist non-degenerate $G$-invariant
pairings
\[
\overline{p}_{H}^{i}:H^{i}\times H^{-i}\rightarrow\mathbf{Q,\;\;\;}%
\overline{p}_{B}^{i}:B^{i}\times U^{-i}\rightarrow\mathbf{Q,\;\;\;}%
\overline{p}_{U}^{i}:U^{i}\times B^{-i}\rightarrow\mathbf{Q}%
\]
which lift the $\sigma^{i}$ in the following sense: under the identification
$H^{i}={\rm H}^{i} (  P_{\mathbf{Q}}^{\bullet} )  $
\begin{equation}
\overline{p}_{H}^{i}=\sigma^{i}:H^{i}\times H^{-i}\;\rightarrow\mathbf{Q};
\end{equation}
and for $b\in B^{i},u\in U^{-i}$
\begin{equation}
\overline{p}_{B}^{i} (  b,u )  = \overline{p}_{U}^{-i} (
u,b )  ;
\end{equation}
furthermore these pairings have the crucial property that for all $i$%
\begin{equation}
\overline{p}_{U}^{i-1}(  u^{i-1},\partial^{-i}u^{-i})
=\overline{p}_{B}^{i}(  \partial^{i-1}(  u^{i-1})
,u^{-i})  .
\end{equation}
The above $\overline{p}_{H}^{i},\;\overline{p}_{B}^{i},\;\overline{p}_{U}^{i}
$ then induce pairings
\[
p^{0}:P_{\mathbf{Q}}^{0}\times P_{\mathbf{Q}}^{0}\rightarrow\mathbf{Q}%
\]%
\begin{equation}
p^{i}: (  P_{\mathbf{Q}}^{i}\oplus P_{\mathbf{Q}}^{-i} )
\times (  P_{\mathbf{Q}}^{i}\oplus P_{\mathbf{Q}}^{-i} )
\rightarrow\mathbf{Q,}\text{\ \ for }i>0
\end{equation}
which are non-degenerate, symmetric and $G$-invariant. By construction, we see
that each $p^{i}$ is an orthogonal sum of forms $p\mid_{H^{i}\oplus H^{-i}}$
on $H^{i}\oplus H^{-i}$, $p\mid_{B^{i}\oplus U^{-i}}$ on $B^{i}\oplus
U^{-i}$, $p\mid_{U^{i}\oplus B^{-i}}$ on $U^{i}\oplus B^{-i}$ when $i>0$; and
when $i=0$   we have $p\mid_{H^{0}}=\overline{p}_{H}^{0}$, $p\mid
_{B^{0}\oplus U^{0}}=\overline{p}_{B}^{0}$.\medskip

With the same notation as used in 3.2.1, the hermitian class in [CPT1]
associated to the pair $ (  P^{\bullet},\sigma )$ is denoted
 $d (  P^{\bullet},\sigma )$; then $\phi (  d (
P^{\bullet},\sigma )  )  ^{-1}$ is that class in ${\rm H}^{\text{s}%
} (  \mathbf{Z}[G])  $ represented by the character
map which takes $\theta_{m}\in R_{G}^{s}$ to the value
\begin{equation}%
{\textstyle\prod_{p<\infty}}
\text{Det} (  \lambda_{p}^{i} )   (  \theta_{m} )  ^{ (
-1 )  ^{i}}\times%
{\textstyle\prod_{i\geq0}}
\mathbf{pf} \left(  T_{m} (  \widetilde{p}^{i} (  a_{0}^{ij}%
,a_{0}^{-ij^{\prime}} )   )  _{j,j^{\prime}} \right)  ^{ (
-1 )  ^{i}}%
\end{equation}
where $T_{m}$ is a symplectic representation with character $\theta_{m}$ and
where for $i>0$ the above Pfaffian term is formed with respect to the
$ \{  a_{0}^{ij},a_{0}^{-ij} \}  _{j}$ and $ \{  a_{0}%
^{ij} \}_{j}$ again denotes a chosen $\mathbf{Q}[G]
$-basis of $P^{i}\otimes_{\mathbf{Z}}\mathbf{Q}$.

(To see why we need to invert the class $\phi (  d (  P^{\bullet
},\sigma )   )  $ to get the above representative, note: firstly,
the map $\phi$ involves inversion of the archimedean coordinate (see (31));
secondly, the formula for the finite coordinate in Definition 4.5 in [CPT1] is
the inverse of the finite coordinate that we use.)

In Theorem 2.9 of [CPT1] it is shown that the class $\phi (  d (
P^{\bullet},\sigma )   )  $ is independent of choices, and so only
depends on the complex $P^{\bullet}$ and the cohomology pairings $\sigma_{i}$.
Recall that for brevity we put
\[
b_{jn}^{im}=r_{G} (  a_{0}^{ij}\otimes w_{mn} )  .
\]
With the notation of the previous sub-section we shall show

\bigskip

\begin{proposition}
(a)\ For each integer $m$%
\[%
{ \sum_{i\geq0}}
\left(  -1\right)  ^{i}n_{m}^{-} (  p^{i} )  =%
{ \sum_{i\geq 0}}
\left(  -1\right)  ^{i}n_{m}^{-} (  \sigma^{i} )  .
\]\smallskip

(b) For each integer $m$%
\[
\text{\rm Pf}_{ (  \sigma\otimes\kappa_{m} )  ^{G}} \left(  \xi_{m} (
\otimes_{i} (  \wedge_{jn}b_{jn}^{im} )  ^{ (  -1 )  ^{i}%
} )  \right)=\left|  G\right|  ^{\chi (  P_{\mathbf{Q}}^{\bullet} )  \theta
_{m} (  1 )  /2}%
{ \prod_{i}}
\mathbf{pf}\left(  \widetilde{p}^{i} (  a_{0}^{ij},a_{0}^{-ij^{\prime}%
} )  _{j,j^{\prime}} \right)  ^{ (  -1 )  ^{i}}%
\]
where
\[
\chi (  P_{\mathbf{Q}}^{\bullet} )  =\left|  G\right|  ^{-1}%
{\textstyle\sum_{i}}
\left(  -1\right)  ^{i}\dim_{\mathbf{Q}} (  P_{\mathbf{Q}}^{i} )
=\left|  G\right|  ^{-1}%
{\textstyle\sum_{i}}
\left(  -1\right)  ^{i}\dim_{\mathbf{Q}} (  \text{\rm H}^{i} (
P_{\mathbf{Q}}^{\bullet} )   )  .
\]\smallskip

(c) Let $h_{m}^{i}$ denote the metric on $\det (  P_{m}^{0} )  $
resp. on $\det (  P_{m}^{i}\oplus P_{m}^{-i} )  $ for $i=0$ resp.
$i>0$ afforded by $p^{0}$ resp. $p^{i}$, and let $x_{i}\in\det (
P_{m}^{i} )  $ all be non-zero. Then
\[
h_{m}^{0} (  x_{0} )
{\textstyle\prod_{i>0}}
h_{m}^{i} (  x_{i}\otimes x_{-i} )  ^{(-1)^{i}}=h_{\sigma,m} (
\xi_{m} (  \otimes_{i}x_{i}^{\left(  -1\right)  ^{i}} )  )  .
\]
\end{proposition}

Before proving this proposition, we first note that part (b) together with
(36) has the following important implication

\begin{theorem}%
\[
\chi_{\text{\rm H}}^{\text{\rm s}} (  P^{\bullet},\sigma )  =\phi (
d (  P^{\bullet},\left|  G\right|  \sigma))  ^{-1}.
\]
\end{theorem}

\bigskip

\begin{remark}
{\rm To understand conceptually why we are obliged to renormalise the forms
$\sigma^{\text{ev}},\sigma^{\text{odd}}$ by a factor $\left|  G\right|$, it
is helpful to consider the special case of a tame Galois extension
$N/\mathbf{Q}$ with Galois group $G$. In [F1] and [CPT1] one works with the
hermitian pair $ (  O_{N}, Tr_{N/\mathbf{Q}} )$; however, in [CPT2]
we are obliged to work with $ (  O_{N},\left|  G\right|  ^{-1}%
Tr_{N/\mathbf{Q}} )  $ for the following reason: the metric associated to
the trace map is of course the pullback of the trivial (or standard) metric on
$\mathbf{Q}$; however, as explained in 4.1 (after Lemma 32), we need to
normalise the pullback metric by a factor $ |  G |  ^{-1}$ to ensure
that it agrees with the original metric on $G$-fixed sections.} 
\end{remark}

\begin{example}
{\rm Let ${\rm HP}={\rm Hyp} (  \mathbf{Z}[G]))  $ denote the free
hyperbolic plane, that is to say the module $\mathbf{Z}[G]
\oplus {\rm Hom}_{\mathbf{Z}} (  \mathbf{Z}[G]  ,\mathbf{Z}%
 )  $ with the evaluation pairing $\sigma$ (see 3.2.5). Then from [F1]
pages 42-43, we know that under the decomposition (26) $d (
\text{HP}, \sigma )  $ is represented by the character function which sends the
symplectic character $\theta_{m}$ to the value
\[
1\times1\times\left(  -1\right)  ^{\theta_{m}\left(  -1\right)  /2}.
\]
Thus we see that by the above theorem the class $\chi_{\text{H}}^{\text{s}%
} (  \text{HP},\left|  G\right|  ^{-1}\sigma )  $ respectively
$\chi_{\text{H}}^{\text{s}} (  \text{HP},\sigma )  $ in 
${\rm H}^{\text{s}} (  \mathbf{Z}[G])  $ is represented by
the character function which maps $\theta_{m}$ to
\[
1\times\left(  -1\right)  ^{\theta_{m}\left(  -1\right)  /2}%
,\;\;\text{resp.\ }1\times\left(  -\left|  G\right|  \right)  ^{\theta
_{m}\left(  -1\right)  /2}.
\]}
\end{example}

\textit{Proof of proposition}. To prove (a) we note that for $i>0$, $p\mid
_{H^{i}\oplus H^{-i}}$,  $ p\mid_{B^{i}\oplus U^{-i}}$, and
$p\mid_{B^{-i}\oplus U^{i}}$ are all hyperbolic; hence
\[
n_{m}^{-} (  p^{i} )  =n_{m}^{-} (  \sigma^{i} )
+\dim (  B^{i} )  +\dim (  U^{i} )
\]
and the result follows at once since
\begin{align*}
\dim (  U^{0} )  +%
{\textstyle\sum_{i>0}}
\left(  -1\right)  ^{i} (  \dim (  B^{i} )  +\dim (
U^{i} )   )   & =\\
\dim (  U^{0} )  +%
{\textstyle\sum_{i>0}}
\left(  -1\right)  ^{i} (  \dim (  U^{i-1} )  +\dim (
U^{i} )   )   & =0.
\end{align*}
In proving (b) we shall use Deligne's ``Koszul rule of signs"' in reordering wedge
products and tensor products; however, we shall see that all terms used have
even grade since they arise as determinants of symplectic isotypic parts of
real representations. Note also that for non-zero elements $l$, $l^{\prime}$ of a
complex line $L$, we shall write $l^{\prime}l^{-1}$ for the complex number
such that $l^{\prime}=\left(  l^{\prime}l^{-1}\right)  l.$

For fixed $m,\;$ for brevity we shall put
\[
r_{i}=   \wedge_{jn}b_{jn}^{im}   \in\det (P_{m}^{i}).
\]
We then use the decompositions
\[
P_{m}^{i}=B_{m}^{i}\oplus H_{m}^{i}\oplus U_{m}^{i}%
\]%
\[
r_{i}=x_{B^{i}}\otimes x_{H^{i}}\otimes x_{U^{i}}%
\]
where for a vector space $V,$ $x_{V}$ denotes an element of $\det V$. As
indicated previously, we note that each of the above $x$-terms is a wedge
product with even grade, and so in the sequel, when we reorder such
multi-tensors with the Deligne sign convention, there is no change of sign.
After multiplying the terms $x_{H^{i}}$ by suitable scalars we may assume that
each $x_{U^{i}}$ is mapped to $x_{B^{i+1}}$ by $\det (  \partial
^{i} )$. Thus
\begin{equation}
\xi_{m} (  \otimes_{i}r_{i}^{ (  -1 )  ^{i}} )  =\otimes
_{i}x_{H^{i}}^{\left(  -1\right)  ^{i}}=x_{H^{0}}\otimes (  \otimes
_{i>0} (  x_{H^{2i}}\otimes x_{H^{-2i}} )   ) \otimes  (
\otimes_{i\geq0} (  x_{H^{2i+1}}\otimes x_{H^{-2i-1}} )
^{-1} )  .
\end{equation}

Next we recall the orthogonal decompositions of the $p^{i}$ given prior to the
statement of the proposition and we note that by (34) we know that
\begin{align*}
\text{Pf}_{p\mid_{B^{-i}\oplus U^{i}}}\left(  x_{B^{-i}}\otimes x_{U^{i}%
}\right)   & =\text{Pf}_{p\mid_{B^{-i}\oplus U^{i}}}\left(  \det\left(
\partial\right)  x_{U^{-i-1}}\otimes x_{U^{i}}\right) \\
& =\text{Pf}_{p\mid_{U^{-i-1}\oplus B^{i+1}}}\left(  x_{U^{-i-1}}\otimes
\det\left(  \partial\right)  x_{U^{i}}\right) \\
& =\text{Pf}_{p\mid_{U^{-i-1}\oplus B^{i+1}}}\left(  x_{U^{-i-1}}\otimes
x_{B^{i+1}}\right) \\
& =\text{Pf}_{p\mid_{B^{i+1}\oplus U^{-i-1}}}\left(  x_{B^{i+1}}\otimes
x_{U^{-i-1}}\right)  .
\end{align*}
We may therefore use Lemma 3 to deduce that the product
\[
\text{Pf}_{p^{0}\otimes\kappa_{m}} (  r_{0} )
{\textstyle\prod_{i>0}}
\text{Pf}_{p^{2i}\otimes\kappa_{m}} (  r_{2i} )
{\textstyle\prod_{i\geq0}}
\text{Pf}_{p^{2i+1}\otimes\kappa_{m}} (  r_{2i+1} )  ^{-1}%
\]
telescopes to
\[
\text{Pf}_{\sigma^{0}\otimes\kappa_{m}} (  x_{H^{0}} )
{\textstyle\prod_{i>0}}
\text{Pf}_{\sigma^{2i}\otimes\kappa_{m}} (  x_{H^{2i}}\otimes x_{H^{-2i}%
} )
{\textstyle\prod_{i\geq0}}
\text{Pf}_{\sigma^{2i+1}\otimes\kappa_{m}}\left(  x_{H^{2i+1}}\otimes
x_{H^{-2i-1}}\right)  ^{-1}%
\]
and by (38) this is equal to
\begin{equation}
\text{Pf}_{\sigma\otimes\kappa_{m}} (  \xi_{m} (  \otimes_{i}%
r_{i}^{ (  -1 )  ^{i}} )   )  .
\end{equation}
On the other hand, by Proposition 23, for each $i$%

\begin{equation}
\mathbf{pf}\left(   |  G |  T_{W_{m}}^{ (  q )  } (
\widetilde{p}^{i} (  a_{0}^{ij},a_{0}^{ij^{\prime}} )   )
\right)  =\text{Pf}_{p^{i}\otimes\kappa_{m}} (  r_{i} )
\end{equation}
and so comparing with (39) above (b) is proved.\medskip

The proof of (c) is very similar to that of (b); we therefore only sketch the
proof. Again we decompose each $x_{i}$ as
\[
x_{i}=x_{B^{i}}\otimes x_{H^{i}}\otimes x_{U^{i}}%
\]
and again, after multiplying the $x_{H^{i}}$ by suitable scalars, we may
suppose that each $x_{U^{i}}$ is mapped to $x_{B^{i+1}}$ by $\det (
\partial^{i} )$; hence, as previously,
\begin{multline}
\ \ \ \ \ \ \ \xi_{m} (  \otimes_{i} (  x_{B^{i}}\otimes x_{H^{i}}\otimes x_{U^{i}%
} )  ^{ (  -1 )  ^{i}} ) = \\ \ \ \   =x_{0}\otimes\left(
\otimes_{i>0}\left(  x_{H^{2i}}\otimes x_{H^{-2i}}\right)  \right)
\otimes (  \otimes_{i\geq0}\left(  x_{H^{2i+1}}\otimes x_{H^{-2i-1}%
}\right)  ^{-1} )  .
\end{multline}
Furthermore by (34) we know that
\begin{align*}
h (  B^{-i}\oplus U^{i} )   (  x_{B^{-i}}\otimes x_{U^{i}%
} )   & =h (  B^{-i}\oplus U^{i} )   (  \det (
\partial )  x_{U^{-i-1}}\otimes x_{U^{i}} ) \\
& =h (  B^{i+1}\oplus U^{-i-1} )   (  x_{U^{-i-1}}\otimes
\det (  \partial )  x_{U^{i}} ) \\
& =h (  B^{i+1}\oplus U^{-i-1} )   (  x_{U^{-i-1}}\otimes
x_{B^{i+1}} )  .
\end{align*}
Hence the product $h_{m}^{0} (  x_{0} )
{\textstyle\prod_{i>0}}
h_{m}^{i} (  x_{i}\otimes x_{-i} )  ^{(-1)^{i}}$ telescopes to
\[
h_{H^{0}} (  x_{H^{0}} )
{\textstyle\prod_{i>0}}
h_{H^{2i}\oplus H^{-2i}} (  x_{H^{2i}}\otimes x_{H^{-2i}} )
{\textstyle\prod_{i\geq0}}
h_{H^{2i+1}\oplus H^{-2i-1}} (  x_{H^{2i+1}}\otimes x_{H^{-2i-1}} )
^{-1}%
\]
and by (41) this is equal to $h_{\sigma,m} (  \xi_{m} (  \otimes
_{i}x_{i}^{ (  -1 )  ^{i}} )   )  $, as required.
\hfill $\square$

\bigskip

\subsection{Proof of Theorem 30.}

We again adopt the notation of the theorem, so that $P^{\bullet}$ is a
perfect $\mathbf{Z} [  G ]  $-complex and $\sigma^{\text{ev}}$,
$\sigma^{\text{odd}}$ are $G$-invariant non-degenerate symmetric forms on the
even and odd parts of the cohomology of $P^{\bullet}$. As in (5.2), \ after
extending $P^{\bullet}$ by an acyclic perfect complex if necessary (which
leaves $\chi_{\text{H}}^{\text{s}} (  P^{\bullet},\sigma )  $
unchanged by Theorem 19), we obtain ``lifts'' $p^{i}$ as in (32)-(35),
and hence $\mathbf{Q} [  G ]  $-valued pairings
\[
\widetilde{p}^{0}:P_{\mathbf{Q}}^{0}\times P_{\mathbf{Q}}^{0}\rightarrow
\mathbf{Q} [  G ]
\]%
\[
\widetilde{p}^{i}: (  P_{\mathbf{Q}}^{i}\oplus P_{\mathbf{Q}}^{-i} )
\times (  P_{\mathbf{Q}}^{i}\oplus P_{\mathbf{Q}}^{-i} )
\rightarrow\mathbf{Q} [  G ]  \mathbf{,}\text{\ \ for }i>0.
\]
Using (5) and Theorem 44 we obtain a decomposition
\[
\chi_{\text{H}}^{\text{s}} (  P^{\bullet},\sigma )  =\phi (
d (  P^{\bullet},\left\{  \left|  G\right|  p^{i}\right\}   )
 )  ^{-1}=x_{1}\times x_{2}\in\text{A}^{\text{s}} (  \mathbf{Z}%
 [  G ]   )  \times\text{S}_{\infty} (  \mathbf{Z} [
G ]  )  .
\]
Here,\ with the notation of (36), the class $x_{1}\times x_{2}$ is represented
by the character map which takes the symplectic character to the value:
\[%
{\textstyle\prod_{p<\infty}}
\text{Det} (  \lambda_{p}^{i} )   (  \theta_{m} )  ^{ (
-1 )  ^{i}}\times%
{\textstyle\prod_{i\geq0}}
\left|  \mathbf{pf} (  T_{m} (   |  G |  \widetilde{p}%
^{i} (  a_{0}^{ij},a_{0}^{-ij^{\prime}} )   )  _{j,j^{\prime}%
} )  ^{ (  -1 )  ^{i}}\right|
\]
respectively
\[
\text{sign}\left(
{\textstyle\prod_{i\geq0}}
\mathbf{pf} (  T_{m} (   |  G |  \widetilde{p}^{i} (
a_{0}^{ij},a_{0}^{-ij^{\prime}} )   )  _{j,j^{\prime}} )
^{ (  -1 )  ^{i}}\right)  .
\]
First we consider  the Arakelov class $x_{1}$. By linearity we may suppose
that $W_{m}$ is a left ideal of $\mathbf{C} [  G ]  $; then we denote
by $ \{  v_{ms} \}  $ an orthonormal basis of $W_{m}$ with respect to
the hermitian form $\nu_{W_{m}}$, and we put
\[
c_{js}^{im}=r_{G} (  a_{0}^{ij}\otimes v_{ms} )  .
\]
Then by Proposition 24 we have
\begin{align*}%
{\textstyle\prod_{i\geq0}}
\left|  \mathbf{pf} (  T_{m} (  |  G |  \widetilde{p}%
^{i} (  a_{0}^{ij},a_{0}^{-ij^{\prime}} )   )  _{j,j^{\prime}%
} )  ^{ (  -1 )  ^{i}}\right|   & =%
{\textstyle\prod_{i\geq0}}
\left|  \det (  T_{m} (   |  G |  \widetilde{p}^{i} (
a_{0}^{ij},a_{0}^{-ij^{\prime}} )  )  _{j,j^{\prime}} )
^{ (  -1 )  ^{i}}\right|  ^{1/2}\\
& =%
{\textstyle\prod_{i\geq0}}
h_{p^{i},m}\left(  \wedge_{js}c_{js}^{im}\right)  ^{\left(  -1\right)  ^{i}}%
\end{align*}
and by Proposition 43(c) this latter real number is equal to
\[
h_{\sigma,m} (  \xi_{m} (  \otimes_{i} (  \wedge_{js}c_{js}%
^{im} )  ^{ (  -1 )  ^{i}} )   )  .
\]
We now consider the signature class $x_{2}$. By Proposition 28 we know that
\[
\text{sign}\left(
{\textstyle\prod_{i\geq0}}
\mathbf{pf} (   |  G |  T_{m} (  \widetilde{p}^{i} (
a_{0}^{ij},a_{0}^{-ij^{\prime}} )  )  )  _{j,j^{\prime}%
}^{\ (  -1 )  ^{i}} )\right)  =%
{\textstyle\prod_{i\geq0}}
(\sqrt{-1})^{n_{m}^{-}\left(  p^{i}\right)  }%
\]
and by Proposition 43(a) the right-hand sign is equal to $\left(  \sqrt
{-1}\right)  ^{n_{m}^{-}\left(  \sigma\right)  },$ as required.\ On
considering the definition of $\chi_{\text{H}}^{\text{s}} (  P^{\bullet
},h_{\sigma} )  $ we see that this then implies Theorem
30.\ \ \ \ \ $\square$

\bigskip

\bigskip

\section{REFERENCES}

\bigskip

$\rm{\ \ \ \ }$[AS] M. F. Atiyah, I. Singer, The index of elliptic operators III, Annals
of Math. 87 (1968), 564-604.\smallskip

[B] J. M. Bismut, Equivariant immersions and Quillen metrics, J. Diff. Geom.
41 (1995), 53-157.\smallskip

[Bl] S. Bloch, De Rham cohomology and conductors of curves, Duke Math. Journal
54 no. 2 (1987), 295-308.\smallskip

[CEPT1] T. Chinburg, B. Erez, G. Pappas, M. J. Taylor, Tame actions of group
schemes: integrals and slices, Duke Math. Journal 82 no. 2 (1996),
269-308.\smallskip\ 

[CEPT2] T. Chinburg, B. Erez, G. Pappas, M. J. Taylor, $\varepsilon$-constants
and Galois structure of de Rham cohomology, Annals of Math. 146 (1997),
411-473.\smallskip\smallskip

[CNT] Ph. Cassou-Nogu\`{e}s, M. J. Taylor, Local root numbers and hermitian
Galois structure of rings of integers, Math. Ann. 263 (1983), 251-261.\smallskip

[CPT1] T. Chinburg, G. Pappas, M. J. Taylor, Duality and hermitian Galois
module structure, Proc. L.M.S. (3) 87 (2003), 54-108.\smallskip

[CPT2] T. Chinburg, G. Pappas, M. J. Taylor, $\varepsilon$-constants and
equivariant Arakelov Euler characteristics, Ann.\ Sci. Ecole Norm. Sup. (4) 35
(2002), no.3, 307-352.\smallskip

[D1] P. Deligne, Valeurs de fonctions L et p\'{e}riodes d'int\'{e}grales,
Proc. Symp. Pure Math., 33 pt 2 (1979), 313-346.\smallskip

[D2] P. Deligne, Les constantes des \'{e}quations fonctionelles de la fonction
L, Springer LNM. 349 (1974), 501-597, Springer-Verlag.\smallskip

[DP] A. Dold, D. Puppe, Homologie nicht-additiver Funktoren, Anwendungen, Ann.
Inst. Fourier 11 (1961), 201-312.\smallskip

[F1] A. Fr\"{o}hlich, Classgroups and hermitian modules, Progress in Maths 48,
Birkh\"{a}user, 1984.\smallskip

[F2] A. Fr\"{o}hlich, Galois module structure of algebraic integers, Springer
Ergebnisse, 3 Folge, Band 1, 1983.\smallskip

[GH] P. Griffiths, J. Harris, Principles of algebraic geometry, Wiley, 1978.\smallskip

[G] A. Grothendieck, On the de Rham cohomology of algebraic varieties, Publ.
Math. I.H.E.S., 29, 95-103. \smallskip

[Go] R. Godement, Topologie alg\'ebrique et th\'eorie des faisceaux. 
Publications de l'Institut de Math\'ematique de l'Universit\'e de Strasbourg, XIII. 
Actualit\'es Scientifiques et Industrielles, No. 1252. Hermann, Paris, 1973. 

[GSZ] H. Gillet, C. Soul\'{e}, with an Appendix by D. Zagier, Analytic torsion
and the arithmetic Todd genus, Topology 30 (1991), 21-54.\smallskip

[H] R. Hartshorne, Algebraic geometry, Springer GTM 52.\smallskip

[I] L. Illusie, Complexe cotangent et d\'{e}formations, Springer Lecture Notes
239 (1971).\smallskip

[KS] K. Kato, T. Saito, Conductor formula of Bloch, preprint 2001.\smallskip

[KM] F. Knusden, D. Mumford, The projectivity of the moduli space of stable
curves I, preliminaries on ``det'' and ``div'', Math. Scand. 39 1 (1976), 19-55.\smallskip

[M] J. S. Milne, \'{E}tale cohomology, Princeton University Press, 1980.\smallskip

[Ms] J. Millson, Closed geodesics and the $\eta$-invariant, Annals of Math.
108 (1974), 1-39.\smallskip

[RS] D. B. Ray, I. M. Singer, Analytic torsion for complex manifolds, Annals
of Math. (2) 98 (1973), 154-177.\smallskip

[R] S. Rosenberg, The Laplacian on a Riemannian manifold, L.M.S. Student Texts
31 (1997), C.U.P.\smallskip

[Sa] T. Saito, Jacobi sum Hecke characters, de Rham discriminant, and the
determinant of $\ell$-adic cohomologies, J. Algebraic Geometry 3 (1994), 411-434.\smallskip

[S] P. Shanahan, The Atiyah-Singer Index theorem, Springer LNM 638, (1978).\smallskip

[So] \ C. Soul\'{e}, D. Abramovich, J.-F. Burnol, J. Kramer, Lectures on
Arakelov Geometry, Cambridge Studies in advanced mathematics 33 (1992), C.U.P.\smallskip

[T] M. J. Taylor, Classgroups of group rings, LMS Lecture Notes 91, 1984, C.U.P.\smallskip

[V] J. L. Verdier, Caract\'{e}ristique d'Euler-Poincar\'{e}, Bull. Soc. Math.
France 101 (1973), 447-448.

\bigskip

Ted Chinburg,

University of Pennsylvania,

Philadelphia PA 19104, USA

ted@math.upenn.edu\bigskip

Georgios Pappas,

Michigan State University,

East Lansing MI 48824, USA

pappas@math.msu.edu\bigskip

Martin Taylor,

UMIST,

Manchester M60 1QD, UK

martin.taylor@umist.ac.uk

\bigskip

\bigskip
\end{document}